\documentclass{amsart}  

\usepackage{amsmath, amsthm, amssymb} 
\usepackage[margin=1in]{geometry} 
\usepackage[colorlinks, allcolors=orange!90!black, pagebackref]{hyperref} 
\usepackage{tikz-cd} 
\usepackage{todonotes} 
\usepackage[nameinlink,capitalize]{cleveref} 
\usepackage{mathtools}
\usepackage[shortlabels]{enumitem} 
\usepackage{microtype} 

\newcommand\setItemnumber[1]{\setcounter{enumi}{\numexpr#1-1\relax}}

\renewcommand*{\backrefalt}[4]{%
\ifcase #1 %
No citations.%
\or
(Cited on page #2).%
\else
(Cited on pages #2).%
\fi
}

\newcommand{\mf}{\mathfrak}

\newcommand{\A}{\mathbb{A}} 
\newcommand{\fa}{\mathfrak{a}} 
\newcommand{\fb}{\mathfrak{b}} 
\newcommand{\uA}{\underline{\mathrm{A}}} 
\newcommand{\F}{\mathbb{F}} 
\newcommand{\N}{\mathbb{N}} 
\newcommand{\fp}{\mathfrak{p}} 
\newcommand{\QQ}{\mathbb{Q}} 
\newcommand{\Q}{\mathtt{Q}} 
\newcommand{\Z}{\mathbb{Z}} 

\newcommand{\yo}{y}

\DeclareMathOperator{\Hom}{Hom} 
\DeclareMathOperator{\Fun}{Fun} 
\DeclareMathOperator{\Spec}{Spec} 
\DeclareMathOperator{\res}{res} 
\DeclareMathOperator{\nm}{nm} 
\DeclareMathOperator{\tr}{tr} 
\DeclareMathOperator{\conj}{c} 
\DeclareMathOperator{\cores}{cores} 
\DeclareMathOperator{\conm}{conm} 
\DeclareMathOperator{\cotr}{cotr} 
\DeclareMathOperator{\coc}{coc} 
\DeclareMathOperator{\id}{id} 
\DeclareMathOperator{\FP}{FP} 
\DeclareMathOperator{\height}{ht} 
\DeclareMathOperator{\coht}{coht} 
\DeclareMathOperator{\RU}{RU} 

\newcommand{\Set}{\mathsf{Set}} 
\newcommand{\set}{\mathsf{set}} 
\newcommand{\Poly}{\mathcal{P}} 
\newcommand{\STamb}{\mathsf{STamb}} 
\newcommand{\Tamb}{\mathsf{Tamb}} 
\newcommand{\CRing}{\mathsf{CRing}} 

\DeclareMathSymbol{\mhyphen}{\mathord}{AMSa}{"39}

\newlength\bbheight
\newcommand{\ghost}{
	\text{%
		\normalfont
		\bbheight=\fontcharht\font`0 
		\raisebox{-0.2\bbheight}[\bbheight][0pt]{%
			\includegraphics[height=1.4\bbheight]{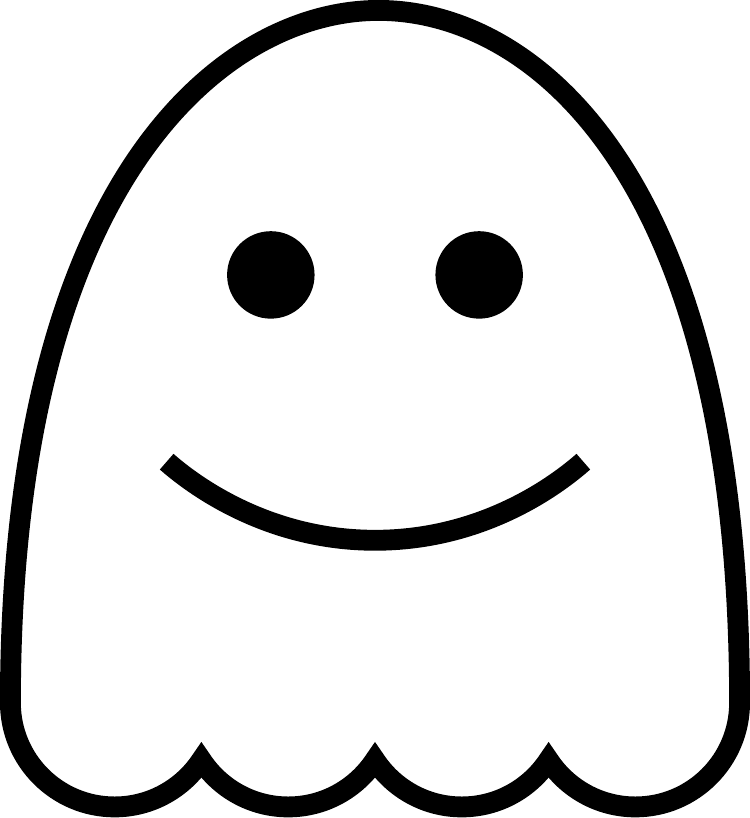}%
		}%
	}%
}
\newcommand{\ghostmap}{\chi}

\numberwithin{equation}{section} 

\newtheorem{letterthm}{Theorem}

\newtheorem{thm}[equation]{Theorem}
\AddToHook{env/thm/begin}{\crefalias{equation}{theorem}}
\newtheorem{proposition}[equation]{Proposition}
\AddToHook{env/proposition/begin}{\crefalias{equation}{proposition}}
\newtheorem{lemma}[equation]{Lemma}
\AddToHook{env/lemma/begin}{\crefalias{equation}{lemma}}
\newtheorem{corollary}[equation]{Corollary}
\AddToHook{env/corollary/begin}{\crefalias{equation}{corollary}}

\theoremstyle{definition}
\newtheorem{definition}[equation]{Definition}
\AddToHook{env/definition/begin}{\crefalias{equation}{definition}}
\newtheorem{notation}[equation]{Notation}
\AddToHook{env/notation/begin}{\crefalias{equation}{notation}}
\newtheorem{example}[equation]{Example}
\AddToHook{env/example/begin}{\crefalias{equation}{example}}
\newtheorem{remark}[equation]{Remark}
\AddToHook{env/remark/begin}{\crefalias{equation}{remark}}
\newtheorem{inductivehypothesis}[equation]{Inductive Hypothesis}
\AddToHook{env/inductivehypothesis/begin}{\crefalias{equation}{inductivehypothesis}}
\crefname{inductivehypothesis}{inductive hypothesis}{inductive hypotheses}

\renewcommand{\ghostmap}{\chi}

\begin{document}

\title{On the Tambara affine line}
\author{David Chan}
\address{Department of Mathematics, Michigan State University, East Lansing, Michigan}
\email{\href{mailto:chandav2@msu.edu}{chandav2@msu.edu}}
\author{David Mehrle}
\address{Department of Mathematics, University of Kentucky, Lexington, Kentucky}
\email{\href{mailto:davidm@uky.edu}{davidm@uky.edu}}
\author{J.D. Quigley}
\address{Department of Mathematics, University of Virginia, Charlottesville, Virginia}
\email{\href{mailto:mbp6pj@virginia.edu}{mpb6pj@virginia.edu}}
\author{Ben Spitz}
\address{Department of Mathematics, Indiana University, Bloomington, Indiana}
\email{\href{mailto:bespitz@iu.edu}{bespitz@iu.edu}}
\author{Danika Van Niel}
\address{Department of Mathematics and Statistics, Binghamton University, Binghamton, New York}
\email{\href{mailto:dvanniel@binghamton.edu}{dvanniel@binghamton.edu}}

\subjclass[2020]{55P91, 19A22, 13A50, 13B21}

\begin{abstract}
	Tambara functors are the analogue of commutative rings in equivariant algebra. Nakaoka defined ideals in Tambara functors, leading to the definition of the Nakaoka spectrum of prime ideals in a Tambara functor. In this work, we continue the study of the Nakoaka spectra of Tambara functors.

	We describe, in terms of the Zariski spectra of ordinary commutative rings, the Nakaoka spectra of many Tambara functors. In particular: we identify the Nakaoka spectrum of the fixed-point Tambara functor of any $G$-ring with the GIT quotient of its classical Zariski spectrum; we describe the Nakaoka spectrum of the complex representation ring Tambara functor over a cyclic group of prime order $p$; we describe the affine line (the Nakaoka spectra of free Tambara functors on one generator) over a cyclic group of prime order $p$ in terms of the Zariski spectra of $\mathbb{Z}[x]$, $\mathbb{Z}[x,y]$, and the ring of cyclic polynomials $\mathbb{Z}[x_0,\ldots,x_{p-1}]^{C_p}$.

	To obtain these results, we introduce a ``ghost construction" which produces an integral extension of any $C_p$-Tambara functor, the Nakaoka spectrum of which is describable. To relate the Nakaoka spectrum of a Tambara functor to that of its ghost, we prove several new results in equivariant commutative algebra, including a weak form of the Hilbert Basis Theorem, going up, lying over, and levelwise radicality of prime ideals in Tambara functors. These results also allow us to compute the Krull dimensions of many Tambara functors.
\end{abstract}

\maketitle

\tableofcontents

\section{Introduction}

Equivariant algebra is the study of equivariant analogues of classical algebraic objects.
In this setting, the role of commutative rings is played by objects called Tambara functors.
Tambara functors have been central to describing multiplicative structures in equivariant stable homotopy theory \cite{BH2015,ElmantoHaugseng}.
Our goal in this paper is to explore some algebro-geometric notions in this context with a view towards potential future applications in equivariant tensor-triangular geometry.

Nakaoka defined \emph{prime ideals} (\cref{def:prime}) and \emph{prime ideal spectra} (\cref{definition: nakaoka spectrum}) of Tambara functors, by analogy to commutative algebra \cite{Nak2012}. To disambiguate the term \emph{spectrum}, we refer to the prime ideal spectrum of a Tambara functor as its \emph{Nakaoka spectrum}. We view the study of Nakaoka spectra as a stepping stone towards the development of highly structured equivariant tensor-triangular geometry.

Nonequivariantly, the Balmer spectrum of a tensor-triangulated category is a space whose geometry contains information about nilpotence and localization in that category \cite{Bal05,Bal10}. For example, the celebrated thick subcategory and nilpotence theorems of Devinatz, Hopkins, and Smith \cite{DHS,HS98} can be reinterpreted in terms of the Balmer spectrum of the homotopy category of finite spectra. There has been a tremendous amount of recent work on extending these ideas to the equivariant setting \cite{BS17,6A19,BGH20}, but these works do not fully incorporate the norm structures present in the equivariant stable homotopy category.  A key insight of Hill and Hopkins \cite{HH16} is that the norm structures are necessary to even define equivariant commutative ring spectra. We posit that the norms should be taken into account when studying the tensor-triangular geometry of equivariant categories.

The norm structures present in equivariant homotopy theory are captured on the level of homotopy by Tambara functors. Thus, a full understanding of equivariant Balmer spectra would benefit from a more complete understanding of the algebraic geometry of Tambara functors. In this paper we expand the theory and available computations of Nakaoka spectra, with the goal of developing intuition and context for later work on equivariant tensor-triangular geometry.

Our first theorem demonstrates that Nakaoka spectra share desirable topological properties with prime ideal spectra of commutative rings.

\begin{letterthm}[{Topology of the Nakaoka spectrum, \Cref{SS:Topology}}]\label{letterthm:topology of spec}
	Let $T$ be a Tambara functor.
	\begin{enumerate}[(a)]
		\item The space $\Spec(T)$ is quasi-compact and sober.
		\item If $T$ is Noetherian, then $\Spec(T)$ is Noetherian.
		\item If $T$ is Noetherian, then a subset $Z \subseteq \Spec(T)$ is closed if and only if it is the closure of a finite set.
	\end{enumerate}
\end{letterthm}

Our next main result relates classical and Nakaoka spectra. For any ring $R$ with action by a finite group $G$ there is an associated fixed-point Tambara functor $\FP(R)$ built out of the fixed point rings $R^H$ for $H$ a subgroup of $G$, and we may consider its Nakaoka spectrum.
On the other hand, the $G$-action on $R$ induces a $G$-action on its prime ideal spectrum $\Spec(R)$. There is a natural relation between these two constructions.

\begin{letterthm}[{Nakaoka spectra as GIT quotients, \Cref{theorem: homeomorphisms of spec}}]\label{letterthm:all the equispecs are the same}
	Let $G$ be a finite group and let $R$ be a commutative ring with a $G$-action. The Nakaoka spectrum of the fixed-point Tambara functor of $R$ is naturally homeomorphic to the  quotient $\Spec(R){/}G \cong \Spec(R^G)$.
\end{letterthm}

\begin{remark}
	The quotient $\Spec(R)/G$ appears in geometric invariant theory (GIT) as the GIT quotient of $\Spec(R)$ by $G$. In general, the GIT quotient of a $G$-scheme $X$ need not be the point-set quotient $X/G$. This suggests the possibility of using Nakaoka spectra to produce models for GIT quotients of more general $G$-schemes; however, this would require the development of non-affine Tambara algebraic geometry, which is beyond the scope of this work.
\end{remark}

\subsection{Computations}

The first fundamental computation in algebraic geometry is $\Spec(\Z)$. The appropriate generalization of this to equivariant algebra is the Nakaoka spectrum of the \emph{Burnside $G$-Tambara functor} $\uA$, the initial object in the category of $G$-Tambara functors.
In \cite{Nak2014a}, Nakaoka computes the spectrum of the Burnside Tambara functor $\uA$ in the case that $G$ is a cyclic $p$-group. This computation was generalized to all finite cyclic groups by Calle and Ginnett in \cite{CG2023}.
To our knowledge, these are the only computations of Nakaoka spectra which have previously appeared in the literature.
Since the Burnside Tambara functor $\uA$ is not a fixed-point functor unless $G$ is trivial, \cref{letterthm:all the equispecs are the same} is orthogonal to the calculations of Nakaoka and Calle--Ginnett and provides new calculations of Nakaoka spectra.
In work of Calle and the authors, we apply the tools developed below to compute the Nakaoka spectrum of the Burnside Tambara functor for any finite group \cite{6A25}.

The next fundamental computation in algebraic geometry is the \emph{affine line} $\A^{\!1} \coloneqq \Spec(\mathbb{Z}[x])$.
This scheme represents the forgetful functor $U$ from commutative rings to sets: there is a bijection
\begin{equation}
	\label{equation: affine line universal property}
	U(R) \cong \mathsf{Sch}(\Spec(R), \A^1),
\end{equation}
natural in the commutative ring $R$. As a formal consequence, $\Z[x]$ has a natural co-(commutative ring) object\footnote{A co-(commutative ring) object (equivalently, a co-commutative co-ring object) in a cocartesian monoidal category $\mathsf{C}$ is an object $X$ with morphisms $0,1 \colon X \to \varnothing$ and $+,\times \colon X \to X \amalg X$ which make $\Hom(X,Y)$ into a commutative ring for all objects $Y$. This is also the same thing as an object $X$ together with a choice of lift of $\Hom(X,{-}) \colon \mathsf{C} \to \Set$ along $U \colon \CRing \to \Set$.} structure in the category of commutative rings.

Blumberg and Hill \cite[Proposition 2.7]{BH2019} show there is a co-Tambara object in the category of $G$-Tambara functors which plays the role of $\Z[x]$, which we call $\uA[{-}]$.
The full data of this co-Tambara object $\uA[{-}]$  consists of a family of Tambara functors $\uA[G/H]$, one for each subgroup $H \leq G$, together with structure maps, satisfying axioms dual to the axioms of a Tambara functor. The individual Tambara functors $\uA[G/H]$ are free objects that together satisfy a universal property analogous to \eqref{equation: affine line universal property}.
In the case that $G$ is a cyclic group of prime order, Blumberg and Hill \cite{BH2019} fully describe the data of $\uA[{-}]$.

As a step towards Tambara algebraic geometry, we study the Nakaoka spectra of the Tambara functors $\uA[{G/H}]$, which together we call the \emph{Tambara affine line} $\underline{\A}^{\!1}$. Let $C_p$ be a cyclic group of prime order.

\begin{letterthm}[{The Tambara affine line over $C_p$, \Cref{SS:A1}}]\label{letterthm:affine line}\
	\begin{enumerate}[(a)]
		\item The Nakaoka spectrum of $\uA[C_p/C_p]$ is an explicit quotient of $\Spec(\Z[x]) \, \amalg \, \Spec(\Z[x,y])$ (\cref{thm:set structure C_p/C_p}).

		\item Let $R = \Z[x_g \mid g \in C_p]^{C_p}$ be the ring of cyclic polynomials. The Nakaoka spectrum of $\uA[C_p/e]$ is an explicit quotient of $\Spec(R) \, \amalg \, \Spec(\Z[x])$ (\cref{thm:set structure C_p/e}).

		\item The co-Tambara structure maps between $\uA[C_p/C_p]$ and $\uA[C_p/e]$ are described in \cref{prop: corestriction,prop: cotransfer,prop: conorm,prop: coconjugation}.
	\end{enumerate}
\end{letterthm}

\begin{remark}
	When studying the prime ideals of a commutative ring $R$, it is often helpful to consider a morphism $R \to S$ such that $\Spec(S) \to \Spec(R)$ is surjective, where the primes of $S$ are easier to describe. We use the same idea to study $\Spec(\uA[G/H])$, where the simpler Tambara functor is obtained through the \emph{ghost} construction described below. To see why mapping to simpler Tambara functors is necessary, the reader is invited to look ahead at the explicit description of $\uA[C_p/e]$ in \cref{Ex:AxGe}.
\end{remark}

We also compute the Nakaoka spectrum of the complex representation Tambara functor $\RU$ when $G = C_p$. The Tambara functor $\RU$ contains the data of the complex representation rings of all subgroups of $G$, together with restriction, induction, and tensor-induction data. Thus, $\RU$ is of central importance in representation theory; it also appears as the zeroth homotopy Tambara functor of equivariant topological $K$-theory.

\begin{letterthm}[{\Cref{thm:RU}}]\label{letterthm: Spec RU}
	Let $G = C_p$ and let $\RU$ denote the complex representation ring Tambara functor. The linearization map $\uA \to \RU$ induces a homeomorphism of Nakaoka spectra $\Spec(\uA) \cong \Spec(\RU)$.
\end{letterthm}

\subsection{Methods}
One of the most powerful tools for studying prime ideal spectra in commutative algebra is the going up theorem. A homomorphism of commutative rings $f \colon R \to S$ satisfies \emph{going up} if containments of primes in $S$ lift to $R$.
This definition translates directly to Tambara functors (\cref{defn: going up}); we prove a going up theorem for equivariant algebra.

\begin{letterthm}[{Going up, \Cref{theorem: going up} and  \Cref{cor: levelwise nil integral map gives surjection on spec}}]
	Let $f \colon R \to S$ be a levelwise-integral morphism of $G$-Tambara functors. Then $f$ satisfies going up.
	If $f$ is moreover levelwise injective, then $\Spec(S) \to \Spec(R)$ is surjective.
\end{letterthm}

To apply the going up theorem, we introduce a construction on $C_p$-Tambara functors which we call the \emph{ghost} (\cref{def: ghost}). The construction and the name are inspired by the ghost ring of a Burnside ring and the ghost Tambara functor of the Burnside Tambara functor that appears in \cite{CG2023}. This also echoes work of Th\'evenaz \cite[Sec.~4]{Thevenaz:SomeRemarks}.

\begin{letterthm}[{The ghost construction, \Cref{section: ghost}}]
	There is a functor $\ghost\colon \Tamb_{C_p} \to \Tamb_{C_p}$ with the following properties:
	\begin{enumerate}[(a)]
		\item For each $C_p$-Tambara functor $T$, we have a natural \emph{ghost map} $\ghostmap_T \colon T \to \ghost(T)$.
		\item The ghost map satisfies going up.
		\item The ghost map induces a surjection $\Spec(\ghost(T)) \twoheadrightarrow \Spec(T)$.
		\item There is a bijection between $\Spec(\ghost(T))$ and $\Spec(T(C_p/e)^{C_p}) \amalg \Spec(\Phi^{C_p}(T))$, where $\Phi^{C_p}(T)$ is a commutative ring called the geometric fixed points of $T$ (\cref{defn: geometric fixed points}).
	\end{enumerate}
\end{letterthm}

\noindent We exploit this bijection of (d) in \cref{section: the ghost and the spectrum} to understand the topology of $\Spec(\ghost(T))$.

The main advantage of this construction is that $\ghost(T)$ is often simpler to describe than $T$, and its spectrum can be described entirely in terms of ordinary commutative algebra (point (d) above). Therefore, combining the going-up theorem with the ghost construction provides a method of computing Nakaoka spectra in essentially every example of interest when $G = C_p$.  It is via this approach that we complete the computations of $\Spec(\uA[C_p/C_p])$ and $\Spec(\uA[C_p/e])$ in \Cref{letterthm:affine line}, at least at the level of posets. Determining the topologies of these Nakaoka spectra takes additional work, and makes use of a version of the Hilbert Basis Theorem which we prove in \cref{cor:weak hilbert basis} below.

Finally, we note that the ghost construction allows us to bound the Krull dimensions of $C_p$-Tambara functors in terms of the Krull dimensions of ordinary commutative rings.

\begin{letterthm}[{\Cref{theorem: krull dimension bounds for integral norm map}}]
	Let $T$ be a $C_p$-Tambara functor such that $\nm\colon T(C_p/e) \to \Phi^{C_p}(T)$ is an integral map. Then
	\[
		\dim(T) \leq
		\dim(\ghost(T)) \leq
		\max\big(\!\dim(T(C_p/e)), \, \dim(\Phi^{C_p}(T))+1\big).
	\]
	If, in addition, the norm is injective, then $\dim(T(C_p/e)) = \dim(\Phi^{C_p}(T))$ and we obtain
	\[
		\dim(T) \leq \dim(T(C_p/e))+1.
	\]
\end{letterthm}
We use this to give computations of Krull dimensions of various Tambara functors of interest in \cref{sec: Krull}.  Aside from measuring the topological complexity of Nakaoka spectra, we note that our dimension computations play an essential role in the proof of \cref{letterthm: Spec RU}.

%
%
%

\subsection{Notation}

We use $G$ to denote a finite group. All rings are assumed commutative and unital. We use $\gamma$ to denote a fixed generator of the cyclic group $C_p$ (we will never have two different cyclic groups in play at the same time). We use symbols like $\subset$, $<$, and $\triangleleft$ to denote proper containments, while $\subseteq$, $\leq$, and $\trianglelefteq$ denote not-necessarily-proper containments. We use $\langle a_1, \ldots, a_n \rangle$ to denote the ideal generated by $a_1, \ldots, a_n$.

\subsection{Acknowledgements}

The authors would like to thank Mike Hill, Clark Lyons, Nat Stapleton, and Yuri Sulyma for helpful discussions. We also thank Will Sawin for comments leading to \Cref{remark: A2 is hard to compute}. Finally, we are grateful to an anonymous referee for their valuable comments and helpful suggestions. 

Chan was partially supported by NSF grant DMS-2135960.
Mehrle, Quigley, and Van Niel were partially supported by NSF grant DMS-2135884.
Quigley was also partially supported by NSF grant DMS-2414922 (formerly DMS-2203785 and DMS-2314082).
Spitz was partially supported by NSF grant DMS-2136090.

\section{Prime Ideals of Tambara Functors}

Let $G$ be a finite group. Tambara functors \cite{Tam1993} are an equivariant analogue of commutative rings -- for each finite group $G$, there is a notion of \emph{$G$-Tambara functor}, and when $G$ is the trivial group, $G$-Tambara functors coincide exactly with commutative rings. There is, likewise, an equivariant analogue of abelian groups, which go by the name of \emph{Mackey functors}. We refer the curious reader to \cite{Str2012} for a thorough introduction to both Mackey and Tambara functors. In this paper, we will discuss only Tambara functors.

\subsection{The Polynomial Category and Tambara functors}

As the name implies, Tambara functors are functors. The domain of any Tambara functor is the category $\Poly_G$, called the \emph{category of polynomials valued in finite $G$-sets}.

The objects of $\Poly_G$ are finite $G$-sets.  The morphisms in $\Poly_G(X,Y)$ are isomorphism classes of \emph{bispans} (also called \emph{polynomials})
\[
	X\leftarrow A\to B\to Y
\]
where $A$ and $B$ are finite $G$-sets and all maps are $G$-equivariant.  The bispan above is isomorphic to a bispan
\[
	X\leftarrow A'\to B'\to Y
\]
if there is a commutative diagram
\[
	\begin{tikzcd}
		X \ar[d,"="]& \ar[l] A \ar[r]\ar[d] & B \ar[d] \ar[r]& Y\ar[d,"="]\\
		X & \ar[l] A' \ar[r] & B'\ar[r]& Y
	\end{tikzcd}
\]
where the two middle vertical maps are $G$-equivariant bijections. We write $[X \leftarrow A \to B \to Y]$ for the isomorphism class of the bispan $X \leftarrow A \to B \to Y$. We refer the reader to \cite[Def. 5.8]{Str2012} for the definition of composition of bispans.

The polynomial category has products, given on objects by disjoint union of $G$-sets. Let $S\colon \Poly_G\to \Set$ be a product preserving functor. In particular, this forces $S(\varnothing)$ to be a singleton. Moreover, given any finite $G$-set $X$, $S(X)$ is a commutative semiring\footnote{A semiring is exactly like a ring, except that it need not have additive inverses. These are also sometimes known as \emph{rigs}.} by \cite[Prop. 5.4]{Str2012}.

\begin{definition}
	A \emph{$G$-semi-Tambara functor} is a product preserving functor $S\colon \Poly_G\to \Set$.  A \emph{$G$-Tambara functor} is a semi-Tambara functor $T$ such that for all $X$ the commutative semiring $T(X)$ is a commutative ring.  A morphism of (semi-)Tambara functors is a natural transformation.  We write $\Tamb_G$ for the category of $G$-Tambara functors.
\end{definition}

By unpacking the structure of $\Poly_G$, one can also give a more explicit definition of Tambara functors. For any equivariant map $f\colon X\to Y$ between finite $G$-sets, there are three corresponding distinguished morphisms in $\Poly_G$:
\begin{equation}\label{eq: distinguished bispans}
	\begin{aligned}
		R_f & = [Y \xleftarrow{f} X \xrightarrow{=} X \xrightarrow{=} X] \\
		T_f & = [X \xleftarrow{=} X \xrightarrow{=} X \xrightarrow{f} Y] \\
		N_f & = [X \xleftarrow{=} X \xrightarrow{f} Y \xrightarrow{=} Y]
	\end{aligned}
\end{equation}
which we call the restriction, transfer, and norm along $f$, respectively. Composition in $\Poly_G$ is defined so that any bispan decomposes as
\[
	[X\xleftarrow{f} A\xrightarrow{g} B\xrightarrow{h} Y ] = T_{h}N_{g}R_f.
\]

The remaining composition laws in $\Poly_G$ amount to a description of how to interchange the restrictions, norms, and transfers to rewrite a composite in this form -- this is originally due to \cite[Sec. 2]{Tam1993} (see also \cite[Sec. 6]{Str2012}).

Let $T$ be a $G$-Tambara functor. Recall that products in $\Poly_G$ are given by disjoint union. Since every finite $G$-set, say $X$, is the disjoint union of its orbits and $T$ is product-preserving, the commutative rings $T(X)$ are determined up to isomorphism by the rings $T(G/H)$ obtained by evaluating $T$ on transitive $G$-sets.  Moreover, since every equivariant function is the composition of folds maps with the union of functions between transitive $G$-sets,\footnote{and maps from the empty $G$-set to transitive $G$-sets, but these do not complicate things because (as previously discussed) $T(\varnothing)$ is the zero ring,} the full data of $T$ is determined by the commutative rings $T(G/H)$ and the functions $T(R_f)$, $T(T_f)$ and $T(N_f)$ when $f\colon G/H\to G/K$ is an equivariant function between transitive $G$-sets.

Let $H$ be a subgroup of $K$.  We write $q_H^K\colon G/H\to G/K$ for the equivariant function uniquely determined by $q_H^K(eH) = eK$. Let $f\colon G/H\to G/K$ be any equivariant function.  If $f(eH) = gK$ then $H\subset gKg^{-1} $ and we can decompose $f$ as the composite
\[
	G/H \xrightarrow{q_H^{gKg^{-1}}} G/gKg^{-1} \xrightarrow{c_{g}} G/K
\]
where $c_g$ is the equivariant bijection determined by sending the identity coset of $G/gKg^{-1}$ to the coset $gK$.

\begin{notation}
	Let $T$ be a Tambara functor. We adopt the following notation:
	\begin{align*}
		\res^K_H & \coloneqq T(R_{q^K_H}) \colon T(G/K)\to T(G/H) \\
		\tr^K_H  & \coloneqq T(T_{q^K_H})\colon T(G/H)\to T(G/K)  \\
		\nm^K_H  & \coloneqq T(N_{q^K_H})\colon T(G/H)\to T(G/K)
	\end{align*}
	And drop the decorations when $H$ and $K$ are understood from context. We call these morphisms the \emph{restriction}, \emph{transfer}, and \emph{norm}, respectively. Additionally, for each $g \in G$, there are \emph{conjugation isomorphisms}
	\[
		\conj_{K,g} \coloneqq T(T_{c_g}) \colon T(G/K)\to T(G/gKg^{-1}).
	\]
	Note that in the case of the conjugation isomorphisms one could equivalently pick $\conj_{K,g} = T(N_{c_g})$ or $\conj_{K,g} = T(R_{c_{g}^{-1}})$, because the morphisms $T_{c_g}, N_{c_g}$, and $R_{c_g^{-1}}$ are equivalent in the polynomial category. We drop the subgroup $K$ from $\conj_{K,g}$ when it is understood from context.

	For each subgroup $K$ and $g \in N_G(K)$, the maps $\conj_{K,g}$ act on $T(G/K)$. It follows from the definitions that if $g \in K$, then $\conj_{K,g}$ acts trivially on $T(G/K)$. Thus we obtain an action of the \emph{Weyl group} $W_G(K) \coloneqq N_G(K) /K$ on $T(G/K)$.
\end{notation}

For a function $\varphi\colon R\to S$ between commutative rings we say that $\varphi$ is \emph{additive} if it is a homomorphism of the underlying additive groups.  We say that $\varphi$ is \emph{multiplicative} if it is a homomorphism of the underlying multiplicative monoids.

\begin{lemma}
	Let $T$ be a $G$-Tambara functor.  For any $H\leq K\leq G$ the maps $\res^H_K$ and $\conj_{K,g}$ are commutative ring homomorphisms. The transfer $\tr^K_H$ is additive, and the map $\nm^K_H$ is multiplicative.
\end{lemma}

\begin{proposition}
	A $G$-Tambara functor $T$ is determined, up to isomorphism, by the commutative rings $T(G/H)$ for all $H\leq G$ and the maps $\res^K_H$, $\tr^K_H$, $\nm^K_H$, and $\conj_{K,g}$ where $H\leq K$ and $g\in G$.
\end{proposition}

The restriction, transfer, norm, and conjugation maps in a Tambara functor are subject to various relations imposed by the composition rules in the category $\Poly_G$.  A full list can be found in \cite{Tam1993}.  We review a few of these below.

\begin{proposition}[Frobenius reciprocity] \label{prop: Frobenius reciprocity}
	Let $T$ be a $G$-Tambara functor, let $H\leq K$ and let $x\in T(G/H)$ and $y\in T(G/K)$.  Then
	\[
		\tr_H^K(x)\cdot y = \tr^K_H(x\cdot \res^K_H(y)).
	\]
\end{proposition}

\begin{proposition}[Tambara reciprocity {\cite[Theorem 2.4]{HM2019}}]\label{proposition: Tambara reciprocity}
	Let $T$ be a $G$-Tambara functor, let $H\leq K$ and let $x,y\in T(G/H)$.  Then
	\[
		\nm_H^K(x+y) = \nm_H^K(x)+\nm_H^K(y)+ \tau
	\] where $\tau$ is a sum of terms, each of which is in the image of a transfer map $\tr^K_{L}$ for some proper subgroup $L<K$.
\end{proposition}

\begin{proposition}[Double coset formulae]
	For any $H,L\leq K \leq G$ we let $g$ run over a set of representatives of the double cosets $L\backslash K \slash H$.  For any $G$-Tambara functor $T$ and $x\in T(G/H)$ we have
	\[
		\res^K_L\circ \tr^K_H(x) = \sum\limits_{LgH \in  L\backslash K \slash H} \tr^L_{L\cap g H g^{-1}}\res^{g H g^{-1}}_{L\cap g H g^{-1}} \conj_{g}(x),
	\]
	\[
		\res^K_L\circ \nm^K_H(x) = \prod\limits_{L g H\in  L\backslash K \slash H} \nm^L_{L\cap g H g^{-1}}\res^{g H g^{-1}}_{L\cap g H g^{-1}} \conj_{g}(x).
	\]
\end{proposition}

\begin{notation}
	When $G = C_p$ (a cyclic group of prime order), we write simply $\res$, $\tr$, $\nm$, $\conj_g$ for $\res_e^{C_p}, \tr_e^{C_p}, \nm_e^{C_p}, \conj_{e,g}$. In this case, the full data of a $C_p$-Tambara functor $T$ is described by the commutative rings $T(C_p/e)$ and $T(C_p/C_p)$ together with the four functions $\res,\ \tr,\  \nm$, and $\conj_\gamma$, where $\gamma$ is a chosen generator of $C_p$. We can write all of this in a \emph{Lewis diagram}:
	\[
		\begin{tikzcd}[row sep=large]
			{T(C_p/C_p)}
			\ar[d, "\res" description]
			\\
			{T(C_p/e)}
			\ar[u, "\tr", bend left=50]
			\ar[u, "\nm"', bend right=50]
			\arrow[from=2-1, to=2-1, loop, in=300, out=240, distance=5mm, "\conj_\gamma"']
		\end{tikzcd}
	\]
	Note that we use the convention that $\tr$ goes on the left and $\nm$ goes on the right. The self-map $\conj_\gamma$ describes the Weyl group action of $C_p$ on $T(C_p/e)$; when this action is trivial, we simply write the word ``trivial.''
\end{notation}
We end this subsection with some examples of Tambara functors.

\begin{example}\label{Ex:FP}
	Let $R$ be a commutative ring with $G$-action.  We write $\FP(R)$ for the \emph{fixed-point Tambara functor} $\FP(R)(G/H) \coloneqq R^H$ with restrictions given by the inclusions and conjugation given by the $G$-action. The transfers and norms are given by
	\[
		\tr_K^H(x)  = \sum_{g \in W_H(K)} g \cdot x, \quad \mathrm{and}\quad   \nm_K^H(x)  = \prod_{g \in W_H(K)} g \cdot x.
	\]
\end{example}

\begin{example}
	For a concrete example of a fixed-point Tambara functor, consider the ring $R = \Z[x,y]$ with the $C_2$-action given by exchanging the variables $x$ and $y$. The Lewis diagram of $\FP(R)$ is:
	\[
		\begin{tikzcd}[row sep=large]
			{\Z[x,y]^{C_2}}
			\ar[d, "\res" description]
			\\
			{\Z[x,y]}
			\ar[u, "\tr", bend left=50]
			\ar[u, "\nm"', bend right=50]
			\arrow["{\mathrm{swap}}" below, from=2-1, to=2-1, loop, in=305, out=235, distance=5mm]
		\end{tikzcd}
	\]
	Here, $\Z[x,y]^{C_2}$ is the ring of the fixed points of the $C_2$-action -- the ring of symmetric polynomials in two variables. Restriction is the inclusion of fixed points, and transfer and norm are given by
	\[
		\tr(f(x,y))  = f(x,y) + f(y,x), \quad \mathrm{and}\quad   \nm(f(x,y))  = f(x,y) f(y,x).
	\]
\end{example}

\subsection{Tambara Ideals}
An ideal of a Tambara functor $T$ is the kernel of a morphism with domain $T$. One can unwind precisely what this means to obtain an internal characterization of ideals.

\begin{definition}[{\cite[Definition 2.1]{Nak2012}}]
	\label{definition: tambara ideals}
	Let $T$ be a $G$-Tambara functor. An \emph{ideal} $I$ of $T$ is a collection of ideals $I(G/H) \subseteq T(G/H)$ closed under transfer, norm, restriction, and conjugation.
	If $I$ is an ideal of $T$, we write $I \trianglelefteq T$.
\end{definition}

An ideal in a ring is proper if and only if it contains no units.  The same is true of ideals in Tambara functors.

\begin{lemma}\label{lemma: ideal is proper if it is levelwise proper}
	Let $I$ be an ideal of a $G$-Tambara functor $T$.  The following are equivalent:
	\begin{enumerate}[(a)]
		\item $I = T$,
		\item there is a subgroup $H\leq G$ such that $I(G/H) = T(G/H)$, and
		\item there is a subgroup $H\leq G$ such that $I(G/H)$ contains a unit.
	\end{enumerate}
\end{lemma}
\begin{proof}
	The first implies the second and the second implies the third by the definition of the containment, so it suffices to prove that the third implies the first. Suppose that $I$ is an ideal with a unit $u \in I(G/H)$.  Since for any $K\leq G$ we have that $I(G/K)$ is an ideal of $T(G/K)$, it suffices to show that for all $K\leq G$ there is a unit in $I(G/K)$. This is straightforward: define $v =\nm^K_e\res^H_e(u)$ and note that this must be a unit because the multiplicativity of norm and restriction implies that both operations preserve units.
\end{proof}

Recall that a commutative ring $R$ is Noetherian if it satisfies the \emph{ascending chain condition} on ideals: any ascending chain of ideals
\[
	I_1\subseteq I_2\subseteq I_3\subseteq \dots\subseteq R
\]
has $I_n=I_{n+1}$ for all sufficiently large $n$.  This definition is easily adapted to Tambara functors. For ideals $I, J \trianglelefteq T$, say that $I \subseteq J$ if $I(G/H) \subseteq J(G/H)$ for all $H \leq G$, and $I \subset J$ if $I \subseteq J$ and $I \neq J$.

\begin{definition}\label{def: noetherian}
	A Tambara functor $T$ is \emph{Noetherian} if any ascending chain of ideals
	\[
		I_1\subseteq I_2\subseteq I_3\subseteq \dots\subseteq T
	\]
	has $I_n=I_{n+1}$ for all sufficiently large $n$.
\end{definition}
Being Noetherian is a finiteness condition that guarantees each ideal of a ring $R$ is finitely generated. We tackle the question of Noetherianity of certain Tambara functors in \cref{sec: hilbert basis} below.

\subsection{Prime Ideals}

The collection of ideals of a Tambara functor $T$ is closed under arbitrary intersection, and so one can speak of the ideal \emph{generated by} a levelwise collection of subsets of $T$.  This allows us to define the product of ideals:

\begin{definition}
	Let $I$ and $J$ be left ideals of a $G$-Tambara functor $T$. The \emph{product} of $I$ and $J$, denoted $IJ$, is defined to be the ideal of $T$ generated by
	\[\{I(G/H)J(G/H)\}_{H \leq G}.\]
\end{definition}

With this definition in hand, one obtains the notion of a prime ideal, which exactly parallels the definition of prime ideal in the theory of rings.

\begin{definition}\label{def:prime}
	Let $P$ be an ideal of a Tambara functor $T$. We say that $P$ is \emph{prime} if it is a proper ideal and for all ideals $I, J \trianglelefteq T$, $IJ \subseteq P$ implies $I \subseteq P$ or $J \subseteq P$.
\end{definition}

This definition is abstractly quite straightforward, but in practice is very difficult to check. We instead make use of an equivalent characterization due to Nakaoka.

\begin{definition}[{\cite[Corollary 4.5]{Nak2012},\cite[Def. 2.6]{CG2023}}]\label{Def:PropositionQ}
	Let $I$ be an ideal of a Tambara functor $T$ and let $H_1,H_2 \leq G$. Let $x \in T(G/H_1)$ and $y \in T(G/H_2)$. Define the proposition 

	\smallskip
	\begin{description}
		\label{qpab}
		\item[$\Q(I,x,y)$] For all subgroups $L, K_1, K_2, H_1, H_2 \leq G$ and $g_1, g_2 \in G$ satisfying the following inclusions:
		      \[
			      g_1 K_1 g_1^{-1} \leq L, \quad
			      g_2 K_2 g_2^{-1} \leq L, \quad
			      K_1 \leq H_1, \quad \text{ and } \quad
			      K_2 \leq H_2,
		      \]
		      \indent we have
		      \[
			      \left(
			      \nm_{g_1 K_1 g_1^{-1}}^L
			      \circ
			      \conj_{K_1,g_1}
			      \circ
			      \res_{K_1}^{H_1}(x)
			      \right)
			      \cdot
			      \left(
			      \nm_{g_2 K_2 g_2^{-1}}^L
			      \circ
			      \conj_{K_2,g_2}
			      \circ
			      \res_{K_2}^{H_2}(y)
			      \right)
			      \in I(G/L).
		      \]
	\end{description}
\end{definition}

\medskip

\begin{proposition}[{\cite[Corollary 4.5]{Nak2012}}]\label{Prop:Nak}
	An ideal $\fp$ of a Tambara functor $T$ is prime if and only if for any $a \in T(G/H)$ and $b \in T(G/H')$, the statement $\Q(\fp,a,b)$ implies $a \in \fp(G/H)$ or $b \in \fp(G/H')$.
\end{proposition}

Terms of the form $\nm_{g^{-1}Kg}^L \circ \conj_{g,K} \circ\res_{K}^{H}(x)$ appear frequently when using statements $\Q(I,x,y)$ to check that an ideal is prime. We introduce some terminology to help discuss them.

\begin{definition}
	\label{def: multiplicative translate}
	A \emph{multiplicative translate of $x$} into $T(G/L)$ is any term of the form
	\[
		\nm_{gKg^{-1}}^L \circ \conj_{K,g} \circ\res_{K}^{H}(x).
	\]
	A \emph{generalized product} of $x \in T(G/H_1)$ and $y \in T(G/H_2)$ is a product of a multiplicative translate of $x$ into $T(G/L)$ and a multiplicative translate of $y$ into $T(G/L)$ (for some $L$).
\end{definition}

With this terminology, the proposition $\Q(I,x,y)$ says simply that all generalized products of $x$ and $y$ lie in $I$. It can be helpful to see the following visual representation of the structure of a generalized product:
\[
	\begin{tikzcd}
		&
		&
		T(G/L)
		\\
		x \in T(G/H_1)
		\ar[d, "\res^{H_1}_{K_1}"]
		&
		&
		&
		&
		T(G/H_2) \ni y
		\ar[d, "\res^{H_2}_{K_2}"]
		\\
		\phantom{a \in }T(G/K_1)
		\ar[r, "\conj_{K_1,g_1}"]
		&
		T(G/g_1 K_1 g_1^{-1})
		\ar[uur, "\nm_{g_1 K_1 g_1^{-1}}^L" description]
		&
		&
		T(G/g_2 K_2 g_2^{-1})
		\ar[uul, "\nm_{g_2 K_2 g_2^{-1}}^L" description]
		&
		T(G/K_2) \phantom{\ni b}
		\ar[l, "\conj_{K_2,g_2}"']
	\end{tikzcd}
\]

\begin{remark}
	Let $I \leq J$ be ideals of a Tambara functor $T$. For any elements $a \in T(G/H)$ and $b \in T(G/H')$, we have $\Q(I,a,b) \implies \Q(J,a,b)$.
\end{remark}

We include the next proposition to illustrate how to work with the statements $\Q(\mf{p},a,b)$.

\begin{proposition}
	Let $T$ be a nonzero Tambara functor. Then $T$ has at least one minimal prime ideal.
\end{proposition}
\begin{proof}
	Since $T$ is nonzero, it has a maximal ideal, which is necessarily prime \cite[Proposition 4.3]{Nak2012}. Thus, by Zorn's lemma, it suffices to prove that an intersection of a nonempty chain of prime ideals is prime.

	Let $\{\mf{p}_i\}_{i \in I}$ be an indexed collection of prime ideals of $T$, where $I$ is a nonempty totally ordered set. Suppose $\mf{p}_i \subseteq \mf{p}_j$ whenever $i \leq j$. We claim that $\mf{p} \coloneqq \bigcap_{i \in I} \mf{p}_i$ is a prime ideal of $T$.

	It is clear that $\mf{p}$ is an ideal. Next, since $I$ is nonempty, pick some $i \in I$. Then $\mf{p} \subseteq \mf{p}_i \subsetneq T$ shows that $\mf{p}$ is a proper ideal of $T$.

	Finally, assume $\Q(\mf{p},x,y)$ holds for some $x \in T(G/H)$ and $y \in T(G/K)$. Towards a contradiction, suppose that neither $x$ nor $y$ lies in $\mf{p}$. Since $x \notin \mf{p}$, there is some $i \in I$ such that $x \notin \mf{p}_i$. Likewise, there is some $j \in I$ such that $y \notin \mf{p}_j$. Since $I$ is totally ordered, we have without loss of generality that $i \leq j$. Then $\mf{p}_i \subseteq \mf{p}_j$, so neither $x$ nor $y$ lies in $\mf{p}_i$. However, $\Q(\mf{p},x,y)$ implies $\Q(\mf{p}_i,x,y)$, so this contradicts the fact that $\mf{p}_i$ is prime.
\end{proof}

\subsection{The Zariski Topology}\label{SS:Topology}

For a Tambara functor $T$, we denote the set of all Tambara prime ideals by $\Spec(T)$.  This forms a poset under inclusion of prime ideals.  We can also give this space a version of the Zariski topology defined as follows (\cite[Definition 4.6]{Nak2012}):  for any ideal $I\subseteq T$, define the set
\[
	V(I) \coloneqq \{\mf{p}\in \Spec(T) \mid I\subseteq \mf{p} \}
\]
and note that the collection of all $V(I)$  as $I$ runs over all ideals of $T$ is closed under finite unions and arbitrary intersections (\cite[Proposition 4.8]{Nak2012}).

\begin{definition}\label{definition: nakaoka spectrum}
	The \textit{Nakaoka spectrum} of a Tambara functor $T$ is the topological space $\Spec(T)$ with the Zariski topology defined by letting the closed sets be those of the form $V(I)$ where $I$ runs over all ideals of $T$.
\end{definition}

Just as is the case for commutative rings, the Nakaoka spectrum of a Tambara functor $T$ is always a quasi-compact space\footnote{Definition: a space is \emph{quasi-compact} if every open cover has a finite subcover. This is the same as what topologists usually mean by \emph{compact}, while the term ``quasi-compact'' is used more commonly used in algebraic geometry.}. The proof is identical to the proof for commutative rings.

\begin{proposition}\label{prop: Spec is qc}
	For any Tambara functor $T$, the space $\Spec(T)$ is quasi-compact.
\end{proposition}
%
%

The Zariski topology on $\Spec(T)$ has a closed basis
\[
	\{V_H(\langle f \rangle) \colon H \leq G, f \in T(G/H)\}
\]
given by
\[
	V_H(\langle f \rangle) = \{\mf{p} \in \Spec(T) \colon f \in \mf{p}(G/H)\}.
\]
However, it can be somewhat cumbersome to think in terms of this closed basis, and in practice one only wants to know if a given set is closed. In commutative algebra, the topology on $\Spec(R)$ is determined by its poset structure when $R$ is Noetherian. We will prove the same statement for Noetherian Tambara functors. Recall that a Tambara functor is Noetherian if its poset of ideals satisfies the ascending chain condition (\Cref{def: noetherian}). We aim to prove the following:

\begin{thm}\label{thm: topology in noetherian case}
	Let $T$ be a Noetherian Tambara functor. A subset $Z \subseteq \Spec(T)$ is closed if and only if it is the closure of a finite set.
\end{thm}

The upshot of this theorem is that the closure of a finite set $\{x_1, \dots, x_n\}$ is the union of the closures of $\{x_1\}, \dots, \{x_n\}$, and the closure of a point $\{\mf{p}\}$ in $\Spec(T)$ is easy to understand: it consists precisely of the prime ideals containing $\mf{p}$. Thus, in the Noetherian setting, there is no need to understand the full ideal structure of $T$ in order to understand the topology of $\Spec(T)$. To prove \Cref{thm: topology in noetherian case}, we mimic the classical story for spectra of commutative rings, and use the theory of sober spaces.

\begin{definition}
	Let $X$ be a topological space. A closed subset $Z \subseteq X$ is \emph{irreducible} if $Z \neq \varnothing$ and $Z$ cannot be written as the union of two proper closed subsets.
\end{definition}

\begin{definition} Let $X$ be a topological space. We say that:
	\begin{enumerate}[(a)]
		\item $X$ is \emph{sober} if every irreducible closed subset $Z \subseteq X$ has exactly one element $z$ such that $Z = \overline{\{z\}}$;
		\item $X$ is \emph{Noetherian} if the poset of open subsets of $X$ satisfies the ascending chain condition (equivalently, the poset of closed subsets of $X$ satisfies the descending chain condition).
	\end{enumerate}
\end{definition}

\begin{lemma}\label{lem: spec is sober}
	Let $T$ be a Tambara functor. Then $\Spec(T)$ is sober.
\end{lemma}
\begin{proof}
	Let $Z$ be an irreducible closed subset of $\Spec(T)$. Then $Z = V(I)$ for some ideal $I \leq T$. We may also suppose, without loss of generality, that $I$ equals the intersection of the prime ideals containing $I$ (this determines the same closed subset as $I$).

	Now let $x \in T(G/H)$ and $y \in T(G/K)$ be such that $\Q(I,x,y)$ holds. Let $\mf{p}$ be any prime ideal containing $I$. We have $\Q(\mf{p},x,y)$, so $x \in \mf{p}$ or $y \in \mf{p}$. Thus, $\mf{p} \in V(I+\langle x \rangle) \cup V(I+\langle y \rangle)$. We conclude that
	\[Z = V(I) = V(I + \langle x \rangle) \cup V(I + \langle y \rangle).\]

	Since $Z$ is irreducible, we may assume without loss of generality that $V(I) = V(I + \langle x \rangle)$. In other words, every prime ideal containing $I$ also contains $x$. Since we assumed that $I$ is the intersection of the primes containing $I$, we also have that $x \in I$.

	Since $Z = V(I)$ is nonempty, we also have that $I \neq T$, so we conclude that $I$ is prime. Now we have $Z = V(I) = \overline{\{I\}}$, and we need only show uniqueness. So, let $\mf{p} \in Z$ be such that $Z = \overline{\{\mf{p}\}}$. Then $V(I) = V(\mf{p})$, so $\mf{p} \subseteq I$ and $I \subseteq \mf{p}$, so $\mf{p} = I$.
\end{proof}

\begin{lemma}\label{lem: T noeth implies spec T noeth}
	Let $T$ be a Noetherian Tambara functor. Then $\Spec(T)$ is a Noetherian topological space.
\end{lemma}
\begin{proof}
	We will prove that closed subsets of $\Spec(T)$ have the descending chain condition. Let $Z_1 \supseteq Z_2 \supseteq \cdots $ be a descending chain of closed subsets of $\Spec(T)$. By definition of the Zariski topology, there exist ideals $I_1, I_2, \dots$ such that $Z_i = V(I_i)$ for all $i$, and we may assume without loss of generality that each $I_i$ is the intersection of the prime ideals containing it.

	Now for any $i$ and any $\mf{p} \in \Spec(T)$ such that $I_{i+1} \subseteq \mf{p}$, we have $\mf{p} \in Z_{i+1}$, so $\mf{p} \in Z_i$, and thus $I_i \subseteq \mf{p}$. Since $\mf{p}$ was arbitrary, we conclude that $I_i$ is contained in all primes which contain $I_{i+1}$. Since $I_{i+1}$ is the intersection of the primes containing it, we conclude that $I_i \subseteq I_{i+1}$. Since $i$ was arbitrary, we have an ascending chain $I_1 \subseteq I_2 \subseteq \cdots$ of ideals of $T$.  Since $T$ is Noetherian, this chain stabilizes, and thus the chain $Z_1 \supseteq Z_2 \supseteq \cdots $ stabilizes.
\end{proof}

\begin{lemma}\label{lem: noeth + sober implies good closed sets}
	Let $X$ be a Noetherian, sober topological space. Let $Z \subseteq X$. The following are equivalent:
	\begin{enumerate}[(a)]
		\item $Z$ is closed;
		\item There is a finite set $Z_0 \subseteq Z$ such that $Z = \overline{Z_0}$.
	\end{enumerate}
\end{lemma}
\begin{proof}
	By definition (b) implies (a). Suppose for contradiction that (a) does not imply (b), i.e. that there exists a closed set $Z$ which is not the closure of any finite set.

	Since $X$ is Noetherian, Zorn's lemma implies that we may suppose $Z$ is minimal among such sets. Since $X$ is sober, $Z$ cannot be irreducible, since irreducible closed sets are the closures of singleton sets. Thus, $Z = A \cup B$ for some closed sets $A, B$ with $A, B \subsetneq Z$. By minimality of $Z$, we have that $A = \overline{A_0}$ and $B = \overline{B_0}$ for some finite sets $A_0, B_0$. But then
	\[Z = \overline{A_0} \cup \overline{B_0} = \overline{A_0 \cup B_0}\]
	is a contradiction.
\end{proof}

The previous three lemmas combine to prove \Cref{thm: topology in noetherian case}.


\begin{proof}[Proof of \Cref{thm: topology in noetherian case}]
	Let $T$ be a Noetherian Tambara functor. Then $\Spec(T)$ is Noetherian and sober by \Cref{lem: spec is sober,lem: T noeth implies spec T noeth}, and the claim follows by \Cref{lem: noeth + sober implies good closed sets}.
\end{proof}

\section{A Weak Hilbert Basis Theorem}
\label{sec: hilbert basis}

In light of \Cref{thm: topology in noetherian case}, the topology on $\Spec(T)$ is much easier to compute when $T$ is a Noetherian Tambara functor, since we need only understand its poset of prime ideals (as opposed to the poset of all ideals). In commutative algebra, the Hilbert Basis Theorem is a central result stating that finite-type algebras over a Noetherian ring are again Noetherian -- in particular, the coordinate rings of classical algebraic varieties are Noetherian.

In the Tambara setting, we will prove a weakened version of this theorem, which still has sufficient power to prove the Noetherianity of the free Tambara functors we consider in this paper. As a first step, we note a sufficient condition for a Tambara functor to be Noetherian.

\begin{lemma}\label{lemma: levelwise finitely generated implies Noetherian}
	If a Tambara functor $T$ is levelwise finitely generated as a commutative ring (i.e. levelwise finite type over $\Z$), then $T$ is Noetherian.
\end{lemma}
\begin{proof}
	If $T$ is levelwise finitely generated over $\mathbb{Z}$ then it is also levelwise Noetherian by the classical Hilbert Basis Theorem. Therefore, any infinite ascending chain of ideals consists of levelwise infinite ascending chains of ideals of Noetherian rings. Such a chain is eventually constant at every level, and there are finitely many levels, so such a chain is eventually constant.
\end{proof}

Our weakened version of the Hilbert Basis Theorem (\Cref{cor:weak hilbert basis}) replaces ``Noetherian'' with ``levelwise finitely generated''; i.e.\ we will show that finite-type algebras over levelwise finitely generated $C_p$-Tambara functors are levelwise finitely generated. The Burnside Tambara functor (the initial object in the category of $G$-Tambara functors) is levelwise finitely generated, thus we will conclude in particular that finite-type free $C_p$-Tambara functors are Noetherian. To begin, we must introduce these finite-type free algebras.

\subsection{Free Tambara Functors}
\label{sec: free Tambara functors}

There are free objects in the category of $G$-Tambara functors, analogous to the free commutative rings (a.k.a.\ polynomial rings) $\Z[X]$ for a set $X$. More generally, for $R$ any commutative ring and $X$ any set, there is an $R$-algebra $R[X]$ (the polynomial algebra) which is free on the set $X$. To recall, the universal property of these free algebras is that there is a bijection
\[
	\mathsf{Alg}_R(R[X], S) \cong \Set(X,S)
\]
which is natural in both the $R$-algebra $S$ and the set $X$. The elements of $X$ are the \emph{generators} of the free algebra $R[X]$. In particular, the free $R$-algebra $R[x]$ on a single generator $x$ has the property that morphisms $R[x] \to S$ are in natural bijective correspondence with elements of $S$. In other words,
\[
	\mathsf{Alg}_R(R[x], S) \cong S.
\]

The situation with Tambara functors is slightly more complicated because a Tambara functor can be free on any \emph{$G$-set} $X$.  The \emph{orbits} of $X$ are the generators of the free algebra on $X$, with each generator living at a level corresponding to the isomorphism type of the orbit.  In particular, if $T$ is a $G$-Tambara functor we can consider the free $T$-algebra $T[G/H]$ which is generated, as a $T$-algebra, by a single generator which lives in $T[G/H](G/H)$. That is, we have
\[
	\mathsf{Alg}_T(T[G/H],S) \cong S(G/G)
\]
naturally in $T$-algebras $S$. The first examples of theses Tambara functors $T[G/H]$ (in the case $H=G$) are introduced by Nakaoka in \cite[Proposition 4.2]{Nak2013}. He also introduces a related construction $T[\mathbf{x}]$, which corepresents $S \mapsto S(G/e)^G$.

In this section, we will describe these free algebras over Tambara functors. Just as $R[X] \cong \Z[X] \otimes R$ naturally in commutative rings $R$, we have $T[X] \cong \uA[X] \boxtimes T$ naturally in Tambara functors $T$, where $\uA[X]$ is the free Tambara functor on $X$ and $\boxtimes$ is the coproduct in the category of Tambara functors.

\begin{proposition}[{\cite[Sec. 9]{Str2012}}]
	The category of $G$-Tambara functors admits binary coproducts, with the coproduct operation denoted $\boxtimes$, called the \emph{box product}.
\end{proposition}

Since $\Tamb_G$ is a category of functors, the Yoneda lemma gives a convenient description of these free Tambara functors (and in particular the initial Tambara functor $\uA = \uA[\varnothing]$). The Yoneda lemma also allows us to see that the free $G$-Tambara functors form a co-Tambara object in $\Tamb_G$.

\begin{definition}
	Let $\mathsf{C}$ be a category admitting all coproducts, and let $G$ be a finite group. A \emph{co-($G$-Tambara) object} in $\mathsf{C}$ is a $G$-Tambara functor valued in $\mathsf{C}^\mathrm{op}$, i.e. a product-preserving functor
	\[F \colon \Poly_G \to \mathsf{C}^\mathrm{op}\]
	such that each commutative semiring object $F(X)$ is in fact a commutative ring object.\footnote{It is indeed a \emph{property} for a commutative semiring object in a cartesian monoidal category to be a commutative ring object -- i.e. any commutative monoid object in a cartesian monoidal category admits at most one antipode map making it into an abelian group object. The proof is the same as the argument that two-sided inverses in a monoid are unique.}
\end{definition}

Note, for example, that a co-($G$-Tambara) object in $\mathsf{Set}^\mathrm{op}$ is exactly the same thing as a $G$-Tambara functor, just as a co-(commutative ring) object in $\mathsf{Set}^\mathrm{op}$ is the same thing as a commutative ring.

With this perspective, we now summarize a key idea in \cite{BH2015}. The Yoneda embedding \[\yo \colon \Poly_G^\mathrm{op} \to \Fun(\Poly_G, \Set)\] preserves coproducts and the image is contained in the category $\STamb_G$ of $G$-semi-Tambara functors (i.e. $\yo(X)$ preserves products for all $X$). We then group-complete each $\yo(X)(Y) = \Poly_G(X,Y)$ to yield a coproduct-preserving functor $\yo^+ \colon \Poly_G^\mathrm{op} \to \Tamb_G$ \cite[Theorem 6.1]{Tam1993}. Thus, we obtain a co-($G$-Tambara) object $\uA[{-}]$ in the category of $G$-Tambara functors, which is simply $\yo^+$ viewed as a functor $\Poly_G \to \Tamb_G^\mathrm{op}$.
The Yoneda lemma together with the universal property of group completion determines a universal property for $\uA[-]$.

\begin{proposition}[{\cite[Proposition 2.7]{BH2019}}]
	\label{prop:free-construction-is-yoneda}
	The functor $\uA[{-}]$ corepresents the identity functor on the category of $G$-Tambara functors, in the sense that
	\[\Tamb_G(\uA[{-}], T) \cong T\]
	naturally in Tambara functors $T$.  We call $\uA[X]$ the \emph{free Tambara functor} on the finite $G$-set $X$.
\end{proposition}

\begin{example}
	The \emph{Burnside Tambara functor}, denoted $\uA$, is the free $G$-Tambara functor on $\varnothing$, i.e. it is the initial object of $\Tamb_G$. In other words, $\uA = \uA[\varnothing]$ in the notation from \Cref{prop:free-construction-is-yoneda}. We note that a bispan \[\varnothing \leftarrow \bullet \to \bullet \to X\] must in fact have the form \[\varnothing \leftarrow \varnothing \to Y \xrightarrow{f} X\] for some morphism $f$, so $\Poly_G(\varnothing,X)$ is in natural bijection with the set of isomorphism classes of objects in the slice category $G\mhyphen\set/X$. The semiring operations on these isomorphism classes (coming from the semiring operations of $\Poly_G(\varnothing,X)$) are
	\[[Y_1 \xrightarrow f X] + [Y_2 \xrightarrow g X] =  [Y_1 \amalg Y_2 \xrightarrow{(f,g)} X]\]
	and
	\[[Y_1 \xrightarrow f X] + [Y_2 \xrightarrow g X] = [Y_1 \times_X Y_2 \rightarrow X].\]
	This gives a description of $\uA(X)$ as the ring of virtual isomorphism classes of objects in $G\mhyphen\set/X$ under the operations described above.

	Now we make a further observation: the functor $G\mhyphen\set/(G/H) \to H\mhyphen\set$ given by sending a morphism $f \colon Y \to G/H$ to its fiber $f^{-1}(eH)$ over the trivial coset is an equivalence of categories. Under this equivalence, we obtain a more familiar description of $\uA$, namely, $\uA(G/H) = A(H)$ is the \emph{Burnside ring} of $H$, the ring of virtual isomorphism classes of finite $H$-sets with addition performed by disjoint union and multiplication performed by cartesian product. The operation $\res_H^K \colon \uA(G/K) \to \uA(G/H)$ comes from the forgetful functor from $K$-sets to $H$-sets, the operation $\tr_H^K \colon \uA(G/H) \to \uA(G/K)$ is given by induction $X \mapsto K \times_H X$, and the operation $\nm_H^K \colon \uA(G/H) \to \uA(G/K)$ is given by coinduction $X \mapsto \operatorname{Map}_H(K,X)$. The conjugation action on $\uA(G/H)$ is given by sending transitive $H$-sets $H/L$ to their conjugates $H/gLg^{-1}$.
\end{example}

\begin{example}
	Let $G = C_p$. Then $\uA$ has the following Lewis diagram
	\[
		\uA \colon \quad
		\begin{tikzcd}[row sep=huge]
			\Z[t]/\langle t^2 - pt \rangle
			\ar[d, "\res"{description}]
			\\
			\Z,
			\ar[u, bend right=50, "\nm"{right}]
			\ar[u, bend left=50, "\tr"{left}]
			\arrow[from=2-1, to=2-1, loop, in=300, out=240, distance=5mm, "{\mathrm{trivial}}"']
		\end{tikzcd}
	\]
	where
	\[
		\res(t)  = p, \quad
		\tr(x)   = tx, \quad \text{ and } \quad
		\nm(x)   = x + \frac{x^p-x}{p} t.
	\]

\end{example}

\begin{example}\label{Ex:AxGG}
	The free $C_p$-Tambara functor on a fixed generator, $\uA[C_p/C_p]$, is described by the Lewis diagram
	\[
		\uA[C_p/C_p] \colon \quad
		\begin{tikzcd}[row sep=huge]
			\Z[x,t,n]/\langle t^2 - pt, tn - tx^p\rangle
			\ar[d, "\res"{description}]
			\\
			\Z[x],
			\ar[u, bend right=50, "\nm"{right}]
			\ar[u, bend left=50, "\tr"{left}]
			\arrow[from=2-1, to=2-1, loop, in=300, out=240, distance=5mm, "{\mathrm{trivial}}"']
		\end{tikzcd}
	\]
	The restriction and transfer are determined by
	\[
		\res(t) = p, \quad
		\res(x) = x, \quad
		\res(n) = x^p, \quad \text{ and } \quad
		\tr(f) = tf.
	\]
	The norm is determined by the norm of $\uA$ ($\nm(k) = k + \frac{k^p - k}{p}t$), the condition $\nm(x) = n$, and the formula for the norm of a sum from \cref{proposition: Tambara reciprocity}.
\end{example}

The free $C_2$-Tambara functor on an underlying generator is described in \cite{BH2019}. A generalization of the calculations therein to $C_{p^n}$ can be found in \cite[Sec. 3]{MQS2024}.

\begin{example}\label{Ex:AxGe}
	The free $C_p$-Tambara functor on an underlying generator is described by the Lewis diagram
	\[
		\uA[C_p/e] \colon \quad
		\begin{tikzcd}[row sep=huge]
			\Z[n]\big[\{t_{\vec v}\}_{\vec v \in \N^{\times p}}\big]/I
			\ar[d, "\res"{description}]
			\\
			\Z[x_0, x_1, \ldots, x_{p-1}],
			\ar[u, bend right=50, "\nm"{right}]
			\ar[u, bend left=50, "\tr"{left}]
			\arrow[from=2-1, to=2-1, loop, in=300, out=240, distance=5mm, "\gamma"']
		\end{tikzcd}
	\]
	where $I$ is the ideal of $\Z[n]\big[\{t_{\vec v}\}_{\vec v \in \N^{\times p}}\big]$ generated by relations
	\[
		t_{\vec 0}^2 - p t_{\vec 0}, \qquad
		t_{\vec v} - t_{\gamma \cdot \vec v}, \qquad
		t_{\vec v} \cdot t_{\vec w} - \sum_{i = 0}^{p-1}t_{\vec v + \gamma^i \cdot \vec w}, \qquad  \text{ and } \qquad
		n t_{\vec v} - t_{\vec v + \vec 1},
	\]
	where $\gamma \cdot \vec v$ is shorthand for the cyclic permutation of indices:
	\[
		\gamma \cdot (v_0, \ldots, v_{p-1}) = (v_1, \ldots, v_{p-1}, v_0),
	\]
	and $\vec 1 = (1,1,\ldots,1) \in \N^{\times p}$.

	The conjugation action on $\uA[C_p/e](C_p/e)$ is given by permuting the variables:
	\[
		\gamma^i \cdot x_j = x_{i + j}
	\]
	with indices taken mod $p$.

	For $\vec v = (v_0, \ldots, v_{p-1}) \in \N^{\times p}$, write $x^{\vec v}$ as shorthand for
	\[
		x^{\vec v} \coloneqq x_0^{v_0}x_1^{v_1} \cdots x_{p-1}^{v_{p-1}}.
	\]
	The restriction and transfer of $\uA[C_p/e]$ are then determined by the formulae:
	\[
		\res(t_{\vec v}) = \sum_{i = 0}^{p-1} x^{\gamma^i \cdot \vec v},   \quad
		\res(n) = x^{\vec 1} = x_0x_1\cdots x^{p-1}, \quad \text{ and } \quad
		\tr(x^{\vec v}) = t_{\vec v}.
	\]
	The norm is determined by the norm of $\uA$, the condition that $\nm(x_i) = n$ for all $i \in \{0,1,\ldots,p-1\}$, and the formula for the norm of a sum from \cref{proposition: Tambara reciprocity}.
\end{example}

\subsection{Noetherianity}
In this section we prove a weak version of the Hilbert Basis Theorem for some Tambara functors. In particular, we will show that $\uA[C_p/C_p]$ and $\uA[C_p/e]$ are Noetherian Tambara functors. The former claim is immediate, but the latter will take a careful argument.

\begin{proposition}\label{prop: AxGG is levelwise finitely generated}
	$\uA[C_p/C_p]$ is levelwise finitely generated, and hence Noetherian.
\end{proposition}
\begin{proof}
	Direct from \Cref{Ex:AxGG} and \Cref{lemma: levelwise finitely generated implies Noetherian}.
\end{proof}

\begin{proposition}
	\label{prop: free Cp is levelwise finitely generated}
	The free Tambara functor $\uA[C_p/e]$ is levelwise finitely generated as a $\Z$-algebra.
\end{proposition}

\begin{proof}

	We show that the ring $R = \uA[C_p/e](C_p/C_p)$ is a finitely generated $\mathbb{Z}$-algebra.  Since $ \uA[C_p/e](C_p/e) = \mathbb{Z}[x_0,\dots,x_{p-1}]$ is also finitely generated as a $\Z$-algebra, this will complete the proof.

	We begin with an overview of the key ideas of the proof.  Let $\vec{v} = (v_0,\dots,v_{p-1})\in \mathbb{N}^{\times p}$ be some vector, and let $K =K(\vec{v})= \sum_{i=0}^{p-1} v_i$.  We consider the subring $R_0\subseteq R$ generated by $t_{\vec{v}}$ where $K\leq p+1$.  Since there are only finitely many such $t_{\vec{v}}$, the ring $R_0$ is finitely generated and hence Noetherian.  One shows that $t_{\vec{v}}$ is in $R_0$ for all $\vec{v}$ by induction on $K$, hence $R_0 =R$.  The inductive step is proved by inducting on the number of zeros present in $\vec{v}$.  The base case is straightforward, but the inductive step requires a further induction on the number of ones in $\vec{v}$ and so on; the nested inductions eventually stop for size reasons.  We now give the complete argument.

	Let $R_0\subseteq R$ be the subring generated by $n$ and all $t_{\vec{v}}$ such that
	\[
		|\vec{v}| \coloneqq \sum\limits_{i=0}^{p-1}v_i \leq p
	\]
	We will show that $R_0$ contains $t_{\vec{v}}$ for all $\vec{v}$ and therefore $R_0=R$ is a finitely generated $\mathbb{Z}$-algebra.

	We begin by establishing some notation. For all $i\in \mathbb{N}$ and $\vec{v} = (v_0,\ldots,v_{p-1})$, let
	\[
		S_i(\vec{v}) = \{j\in \{0,\dots,p-1\}\mid v_j=i\}
	\]
	and let $f_i(\vec{v}) = |S_i(\vec{v})|$.  Let
	\[
		T_i(\vec{v}) = \bigcup_{j\geq i} S_j(\vec{v})
	\]
	be the set of indices $k$ where $v_k\geq i$. If the vector $\vec{v}$ is fixed, we will suppress it from the notation. Finally, for any subset $A\subseteq \{0,\dots,p-1\}$ we write $e(A)\in \mathbb{N}^p$ for the tuple given by
	\[
		e(A)_i = \begin{cases}
			1 & i\in A,    \\
			0 & i\notin A.
		\end{cases}
	\]

	Let us write $K = |\vec{v}|$, we proceed by induction on $K$ to show that $t_{\vec v} \in R_0$. The base case, when $K\leq p$ we have $t_{\vec{v}}\in R_0$, holds by definition of $R_0$.  We adopt the inductive hypothesis:
	\begin{enumerate}
		\setItemnumber{-1}
		\item Suppose $t_{\vec{v}}\in R_0$ whenever $|\vec{v}|<K$.
	\end{enumerate}

	Let $\vec{v}$ be a $p$-tuple with $|\vec{v}|=K$.  We will show that $t_{\vec{v}}\in R_0$ by induction on $f_0(\vec{v})$.

	For the base case $f_0(\vec{v})=0$ we have
	\[
		t_{\vec{v}} = n\cdot t_{\vec{v}-\vec{1}}.
	\]
	Since $n\in R_0$ by definition, and $t_{\vec{v}-\vec{1}}\in R_0$ by inductive hypothesis (-1), we have $t_{\vec{v}}\in R_0$ whenever $f_0(\vec{v})=0$, establishing the base case.  We adopt the inductive hypothesis:
	\begin{enumerate}
		\setItemnumber{0}
		\item Suppose there is an $m_0>0$ so that $t_{\vec{v}}\in R_0$ whenever $|\vec{v}|=K$ and $f_0(\vec{v})<m_0$.
	\end{enumerate}

	For the inductive step, suppose we have a $\vec{v}$ with $|\vec{v}|=K$ and $f_0(\vec{v})=m_0$.  We show that $t_{\vec{v}}\in R_0$ by induction on $f_1(\vec{v})$.

	For the base case $f_1(\vec{v})=0$ we have
	\[
		t_{\vec{v}} = t_{\vec{v}-e(T_1)}\cdot t_{e(T_1)} - \sum\limits_{1\leq i\leq p-1} t_{\vec{v}-e(T_1)+\gamma^ie(T_1)}.
	\]
	By assumption, $f_0(\vec{v})\neq 0$ so $T_1\neq \{0,\dots,p-1\}$.  The assumption that $|\vec{v}|=K>p$ implies that $T_1\neq \varnothing$.  Since $p$ is prime the only $\gamma^i$-fixed points in the power set of $\{0,\dots,p-1\}$ are the whole set and the empty set, so  $\gamma^iT_1\neq T_1$ for any $1\leq i\leq p-1$. Since $T_1\amalg S_0=\{0,\dots,p-1\}$, it must be that $\gamma^i\cdot T_1\cap S_0\neq \varnothing$ from which it follows that $f_0(\vec{v}-e(T_1)+\gamma^ie(T_1))<f_0(\vec{v})=m_0$, since somewhere a $0$ from $\vec{v}$ became a $1$, and because $f_1(\vec{v})=0$ there are no new zeros in the vector.  Thus by inductive hypothesis (0) every term on the right is in $R_0$ so we have established the base case of the induction on $f_1(\vec{v})$.  We adopt the inductive hypothesis:
	\begin{enumerate}
		\setItemnumber{1}
		\item Suppose there is an $m_1>0$ so that $t_{\vec{v}}\in R_0$ whenever $|\vec{v}|=K$, $f_0(\vec{v})=m_0$, and $f_1(\vec{v})<m_1$.
	\end{enumerate}

	For the inductive step, suppose we have a $\vec{v}$ with $|\vec{v}|=K$, $f_0(\vec{v})=m_0$, and $f_1(\vec{v})=m_1$.  We show that $t_{\vec{v}}\in R_0$ by induction on $f_2(\vec{v})$.  At this point the reader might be beginning to see the pattern, but we write out one more case for concreteness.

	For the base case $f_2(\vec{v})=0$ we have
	\[
		t_{\vec{v}} = t_{\vec{v}-e(T_2)}\cdot t_{e(T_2)} - \sum\limits_{1\leq i\leq p-1} t_{\vec{v}-e(T_2)+\gamma^ie(T_2)}.
	\]
	We have $T_2\amalg S_0\amalg S_1 = \{0,\dots,p-1\}$ and $\gamma^i\cdot T_2\neq T_2$ so we have $\gamma^i\cdot T_2\cap (S_0\amalg S_1)\neq \varnothing$.  If $\gamma^i\cdot T_2\cap S_0\neq \varnothing$ then $f_0(t_{\vec{v}-e(T_2)+\gamma^ie(T_2)})<m_0$ by the same argument as before.  If $\gamma^i\cdot T_2\cap S_0 = \varnothing$ then  $f_0(t_{\vec{v}-e(T_2)+\gamma^ie(T_2)})=m_0$ and $f_1(t_{\vec{v}-e(T_2)+\gamma^ie(T_2)})<f_1(\vec{v})=m_1$ since one of the $1$'s in $\vec{v}$ became a $2$, because  $\gamma^i\cdot T_2\cap S_1\neq \varnothing$, and no $0$'s became $1$'s because $T_2\cap S_0=\varnothing$.  Thus every term on the right is in $R_0$ so we have established the base case for the induction on $f_2(\vec{v})$. We adopt the inductive hypothesis:
	\begin{enumerate}
		\setItemnumber{2}
		\item Suppose there is an $m_2>0$ so that $t_{\vec{v}}\in R_0$ whenever $|\vec{v}|=K$, $f_0(\vec{v})=m_0$, $f_2(\vec{v})=m_1$, and $f_2(\vec{v})<m_2$.
	\end{enumerate}
	One proceeds in this way, finding numbers $m_i>0$ and inductive hypotheses:
	\begin{enumerate}[(i)]
		\item Suppose there is an $m_i>0$ so that $t_{\vec{v}}\in R_0$ whenever $|\vec{v}|=K$, \[f_0(\vec{v})=m_0, \dots, f_{i-1}(\vec{v})=m_{i-1},\] and $f_i(\vec{v})<m_i$.
	\end{enumerate}
	The base case of inductive hypothesis (i) is implied by (-1) through (i-1).  Moreover, if we can prove the inductive step for any (i) this will immediately imply the inductive step for (i-1) and working all the way backwards it establishes the inductive step for (-1) completing the proof.

	Eventually we arrive at inductive hypothesis $(K+1)$, where again we can establish the base case, and we are left to consider the inductive step: we need to show that if $m_{K+1}>0$ and $f_{K+1}(\vec{v})=m_{K+1}$ then $t_{\vec{v}}\in R_0$.  But it is impossible for $|\vec{v}|=K$ and $f_{K+1}(\vec{v})\neq 0$ to both be true, and thus only the base case matters, as the inductive step is vacuously true, and the proof is complete.
\end{proof}

Using \Cref{lemma: levelwise finitely generated implies Noetherian} we have the following.

\begin{corollary}
	The free Tambara functor $\uA[C_p/e]$ is Noetherian.
\end{corollary}

We have established that both free $C_p$-Tambara functors on a single generator are levelwise finitely generated. It follows by the formula for box products \cite[Definition 1.2.1]{Maz2013} that \emph{all} free $C_p$-Tambara functors are levelwise finitely generated, and indeed all free $T$-algebras are levelwise finitely generated whenever $T$ is levelwise finitely generated. This general statement is what we call the \emph{Weak Hilbert Basis Theorem}. We note the following box-product argument was also obtained independently by Noah Wisdom.

\begin{corollary}[Weak Hilbert Basis Theorem]\label{cor:weak hilbert basis}
	Let $X$ be a finite $C_p$-set and let $T$ be a levelwise finitely generated $C_p$-Tambara functor. The free $T$-algebra $T[X]$ is levelwise finitely generated and, in particular, Noetherian.
\end{corollary}

\begin{proof}
	We start with the case $T = \uA$. Decompose $X$ into transitive $C_p$-sets as
	\[
		X \cong (C_p/e)^{\sqcup i} \sqcup (C_p/C_p)^{\sqcup j}.
	\]
	Then
	\[
		\uA[X] \cong \uA[C_p/e]^{\boxtimes i} \boxtimes \uA[C_p/C_p]^{\boxtimes j}.
	\]
	By the formula for box products from \cite[Definition 1.2.1]{Maz2013},
	\[
		\uA[X](C_p/e) \cong \uA[C_p/e](C_p/e)^{\otimes i} \otimes \uA[C_p/C_p](C_p/e)^{\otimes j},
	\]
	which is finitely generated by \Cref{prop: AxGG is levelwise finitely generated,prop: free Cp is levelwise finitely generated}. Similarly, $\uA[X](C_p/C_p)$ is a quotient of
	\[
		\left(  \uA[C_p/e](C_p/C_p)^{\otimes i} \otimes \uA[C_p/C_p](C_p/C_p)^{\otimes j}  \right)
		\oplus
		\left( \uA[C_p/e](C_p/e)^{\otimes i} \otimes \uA[C_p/C_p](C_p/e)^{\otimes j} \right),
	\]
	which is finitely generated by \Cref{prop: AxGG is levelwise finitely generated,prop: free Cp is levelwise finitely generated}. Thus, $\uA[X]$ is levelwise finitely generated, and we apply \Cref{lemma: levelwise finitely generated implies Noetherian} to show that $\uA[X]$ is Noetherian.

	The case where $T$ is a levelwise finitely generated $C_p$-Tambara functor follows from the case $T = \uA$ together with the isomorphism $T[X] \cong T \boxtimes \uA[X]$ and the formulae for box products as above.
\end{proof}

\begin{remark}
	The classical Hilbert Basis Theorem says that if a commutative ring $R$ is Noetherian, then $R[x]$ is Noetherian. To show that $T[X]$ is Noetherian, we assume a much stronger condition (levelwise finitely generated) than the classical Hilbert Basis Theorem demands.
\end{remark}

\begin{remark}
	Now that we know the free $C_p$-Tambara functors are Noetherian, \Cref{thm: topology in noetherian case}  gives a concise understanding of the closed subsets of $\Spec(\uA[C_p/C_p])$ and $\Spec(\uA[C_p/e])$ in terms of the containments of prime ideals.
\end{remark}

\section{Tambara Primes are Levelwise Radical}
\label{sec: levelwise radical}

Although the levels of a prime ideal of a Tambara functor are not necessarily prime, they do exhibit some levelwise structure: each level of a prime ideal of a Tambara functor is a radical ideal.  This section is devoted to the proof of this fact, which is \cref{thm: Tambara prime are radical} below.

\subsection{\texorpdfstring{$G$}{G}-prime Ideals}

Before considering prime ideals of Tambara functors we first consider certain ideals in rings with $G$-action.  If $G$ acts on a ring $R$ through ring homomorphisms and $I$ is an ideal of $R$ we say that $I$ is \emph{$G$-invariant} if $g(I)\subseteq I$ for all $g\in G$.

\begin{definition} \label{defn: G-prime ideals}
	Let $R$ be a commutative ring with $G$-action by ring homomorphisms.
	We call a $G$-invariant ideal $\mf{p}$ a \emph{$G$-prime ideal} if
	\begin{enumerate}[(a)]
		\item $\mf{p} \neq R$, and
		\item for $x,y \in R$, if $x(g\cdot y)\in \mf{p}$ for all $g\in G$, then either $x$ or $y$ is in $\mf{p}$.
	\end{enumerate}
	Let $\Spec_G(R)$ denote the set of all $G$-prime ideals of $R$.
\end{definition}

\begin{remark}
	Note that while condition (b) appears to be asymmetric in $x$ and $y$, it is not.  Indeed, because $I$ is $G$-invariant we have $x(g \cdot y)\in I$ for all $g\in G$ if and only if for all pairs $g_1,g_2\in G$ we have $(g_1\cdot x) (g_2\cdot y)\in I$.
\end{remark}

If $G$ is the trivial group, or more generally if the action of $G$ is trivial on a $G$-prime ideal, this recovers precisely the standard notion of a prime ideal of a commutative ring. Every prime ideal is $G$-prime, but not all $G$-prime ideals are prime. These $G$-prime ideals arise naturally in the study of Tambara functors, as shown in the following lemma.

\begin{lemma}\label{lem: bottom is G-prime}
	Let $\mf{p}$ be a prime ideal of a Tambara functor $T$. Then $\mf{p}(G/e)$ is a $G$-prime ideal of $T(G/e)$.
\end{lemma}
\begin{proof}
	By \cref{definition: tambara ideals}, the ideal $\mf{p}(G/e)$ is a $G$-invariant ideal. It suffices to show that for $x,y \in T(G/e)$, if $x(g \cdot y)\in \mf{p}(G/e)$ for all $g\in G$, then either $x$ or $y$ is in $\mf{p}(G/e)$.  Since there are no proper subgroups of the trivial group the proposition $\Q(\mf{p},x,y)$ (\cref{Def:PropositionQ}) consists of checking whether or not all products of the form
	\[
		\nm^K_e(g_1 \cdot x) \nm^K_e(g_2 \cdot y)
		= \nm^K_e(g_1 \cdot (x(g_1^{-1}g_2 \cdot y)))
	\]
	are in $\mf{p}$ for all $K\leq G$ and $g_1,g_2\in G$.  Since we assumed that $x (g_1^{-1}g_2 \cdot y)\in \mf{p}(G/e)$ these elements are always in $\mf{p}$.  Since $\mf{p}$ is a prime ideal, and $\Q(\mf{p},x,y)$ is satisfied, we must have that $x$ or $y$ is in $\mf{p}(G/e)$, proving the claim.
\end{proof}

\begin{remark}
	It is not true in general that Tambara prime ideals are levelwise $W$-prime for the actions of the Weyl groups $W$. For example, the zero ideal is prime in the Burnside $C_p$-Tambara functor $\uA$, but $\uA(C_p/C_p) = \mathbb{Z}[t]/(t^2-pt)$ is not a domain so the zero ideal is not prime in $\uA(C_p/C_p)$. More generally, for any finite group $G$ the zero ideal of $\uA(G/H)$ is never $W_G(H)$-prime for $H\neq e$, despite the fact that $\underline{0}$ is a prime ideal of $\uA$.
\end{remark}

It turns out that $G$-prime ideals are always radical.

\begin{proposition}\label{prop: G prime is radical}
	Let $\mf{p}\subset R$ be a $G$-prime ideal of a $G$-commutative ring.  Then $\mf{p}$ is a radical ideal.
\end{proposition}

\begin{proof}
	Suppose that $x^n\in \mf{p}$ for some $n$.  Let $A$ be a finite, non-empty multiset composed of elements of $G$.  We write
	\[
		\operatorname{n}_A(x) = \prod\limits_{g\in A} g\cdot x.
	\]
	We will show that $\operatorname{n}_A(x)\in \mf{p}$ for all choices of non-empty multisets $A$; the proposition statement is the specific choice $A = \{1\}$.

	First, note that if $|A|>n|G|$ then there is some $g\in G$ so that $\operatorname{n}_A(x)$ is divisible by $(g\cdot x)^n = g\cdot(x^n)\in \mf{p}$, hence $n_A(x)\in \mf{p}$.  Suppose that for some $N$ we know that $|A|>N$ implies that $\operatorname{n}_A(x)\in \mf{p}$.  The result follows by induction if we show that $\operatorname{n}_A(x)\in \mf{p}$ for all $A$ with $|A|=N$.  Note that for any $g\in G$ we have $g\cdot \operatorname{n}_A(x) = \operatorname{n}_{g \cdot A}(x)$ and so for any $g\in G$ we have
	\[
		\operatorname{n}_A(x)(g\cdot \operatorname{n}_A(x)) = \operatorname{n}_A(x) \operatorname{n}_{g \cdot A}(x) = \operatorname{n}_{A\amalg (g\cdot A)}(x),
	\]
	and this element is in $\mf{p}$ because $|A\amalg (g\cdot A)|=2N>N$.  Since $\mf{p}$ is $G$-prime we must have $\operatorname{n}_A(x)\in \mf{p}$ and so we are done.
\end{proof}

\begin{corollary}\label{cor: G prime determined by fixed points}
	If $\mf{p}$ and $\mf{q}$ are two $G$-prime ideals such that $\mf{p}^G = \mf{q}^G$ then $\mf{p} = \mf{q}$.
\end{corollary}
\begin{proof}
	We will show that $\mf{p}\subseteq \mf{q}$; the reverse inclusion follows by symmetry.  Let $a\in \mf{p}$ and note that by \cref{prop: G prime is radical} we will be done if we show that $a^{|G|}\in \mf{q}$.

	We have that $a$ is a solution to the equation
	\[
		\prod\limits_{g\in G} (x-ga)=0.
	\]
	The coefficients in the polynomial expansion of this product are (up to $\pm1$) the elementary symmetric polynomials on the set $\{ga\}_{g\in G}$. In particular, the coefficient of every term in the polynomial except $x^{|G|}$ is an element of $\mf{p}^{G} = \mf{q}^G$.  Evaluating at $x=a$, it follows that $a^{|G|} \equiv 0\mod \mf{q}$, so $a\in \mf{q}$ by \cref{prop: G prime is radical}.
\end{proof}

\subsection{Tambara Primes are Levelwise Radical}

Our next goal is to prove the following theorem, which shows Tambara prime ideals are composed of radical ideals, we do this by using a variation of the above arguments for $G$-prime ideals.

\begin{thm}\label{thm: Tambara prime are radical}
	Every Tambara prime ideal is levelwise radical.  That is, if $\mf{p}\subset T$ is a Tambara prime ideal then $\mf{p}(G/H)\subset T(G/H)$ is a radical ideal for all $H\leq G$.
\end{thm}

Let $T$ be a $G$-Tambara functor and let $\mf{p}$ be a prime ideal of $T$.
We know that $\mf{p}(G/e)$ is a $G$-prime ideal and thus is radical by \cref{prop: G prime is radical}.
This forms the base case for an inductive argument.

\begin{inductivehypothesis}
	\label{inductive step}
	Let $H$ be a subgroup of $G$ such that $\mf{p}(G/L)$ is a radical ideal of $T(G/L)$ for all proper subgroups $L < H$. Let $x \in T(G/H)$ such that $x^n \in \mf{p}(G/H)$.
\end{inductivehypothesis}

\cref{thm: Tambara prime are radical} is proved once we show that $x \in \mf{p}(G/H)$, i.e. $\mf{p}(G/H)$ is a radical ideal of $T(G/H)$.
With this setup in mind, we make the following definition.

\begin{definition}
	For any $y \in T(G/K)$ for any $K \leq G$, define
	\[
		P(y) \coloneqq
		\sup\left\{
		m \in \Z
		\ \middle\vert\
		y = \prod_{i=1}^m
		\nm_{g_i L_ig_i^{-1}}^K
		\conj_{L_i,g_i}
		\res^H_{L_i}(x)
		\right\}.
	\]
	That is, $P(y)$ is the supremum over the number of terms in a factorizations of $y$ into a product of multiplicative translates of $x$.
	If $y$ cannot be factored into a product of multiplicative translates of $x$ then $P(y) = -\infty$.
\end{definition}

Note that $P(x) \geq 1$.

\begin{lemma}
	\label{lemma: bound implies radical}
	Suppose there is an integer $M \geq 0$ such that for any $K\leq G$ and any $y \in T(G/K)$,  $P(y) > M$ implies that $y \in \mf{p}(G/K)$. Then $M = 0$ suffices: $P(y) > 0$ implies that $y \in \mf{p}(G/K)$.
\end{lemma}

\begin{proof}
	If $M = 0$, then we are done. Suppose then that $M>0$ and assume there is $y \in T(G/K)$ with $P(y) = M$. The proposition $\Q(\mf{p}, y, y)$ (\cref{Def:PropositionQ}) is a statement about whether certain elements $z$ are in $\mf{p}$, where $z$ is the product of two multiplicative translates of $y$. Such a $z$ necessarily satisfies the inequality
	\begin{equation}
		\label{equation: P(z) > M}
		P(z) \geq 2 P(y) = 2M > M.
	\end{equation}
	Thus, $z$ is in $\mf{p}$ because $P(z) > M$. Since $\mf{p}$ is prime it follows that $y \in \mf{p}$. Thus, $P(y) \geq M$ suffices to determine membership in $\mf{p}$; we don't need the strict inequality.

	But if $P(y) > M-1$ suffices, we can repeat the argument to see that $P(y) > M-2$ suffices. Repeating, we obtain the statement: $P(y) > 0$ implies that $y \in \mf{p}$. We may not iterate beyond this because $M = 0$ does not satisfy the inequality \eqref{equation: P(z) > M}.
\end{proof}

With this lemma all that remains to prove \cref{thm: Tambara prime are radical} is to show there is some $M$ such that $P(y) > M$ implies that $y \in \mf{p}$. In fact, it is enough to show that for all $K \leq G$, there is some integer $M(K)$ such that $y \in T(G/K)$ and $P(y) > M(K)$ implies that $y \in \mf{p}(G/K)$. We may then take $M$ to be the largest of the $M(K)$. We divide this into three cases: subgroups not containing a conjugate of $H$, subgroups conjugate to $H$, and subgroups properly containing a conjugate of $H$.

\begin{lemma}
	\label{lemma: existence of bound case 1}
	Under the \cref{inductive step}, let $K$ be a subgroup of $G$ such that $H$ is not subconjugate to $K$.
	Then $P(y) > 0$ implies that $y \in \mf{p}(G/K)$, for any $y \in T(G/K)$.
\end{lemma}

\begin{proof}
	In this case, any multiplicative translate of $x$ into $T(G/K)$ has the form
	\[
		y = \nm_{g Lg^{-1}}^K \conj_g \res^H_L(x),
	\]
	where $L$ must be a proper subgroup of $H$. By the inductive assumptions, we have $\mf{p}(G/L)$ is radical for all $L < H$. Therefore, $\res^H_L(x^n) = \res^H_L(x)^n \in \mf{p}(G/L)$ yields $\res^H_L(x) \in \mf{p}(G/L)$. Therefore, any multiplicative translate of $x$ into $T(G/K)$ is in $\mf{p}$. In other words, if $P(y) > 0$ then $y \in \mf{p}(G/K)$.
\end{proof}

\begin{lemma}
	\label{lemma: existence of bound case 2}
	Under the \cref{inductive step}, let $K = gHg^{-1}$ for some $g \in G$.
	Then there is an integer $M = M(K) \geq 0$ such that $P(y) > M$ implies that $y \in \mf{p}(G/K)$, for any $y \in T(G/K)$.
\end{lemma}

\begin{proof}
	Suppose that $y\in T(G/K)$.  Note that any multiplicative translate of $x$ which is in $T(G/gHg^{-1})$  has the form
	\[
		y = \nm^{gHg^{-1}}_{g_1Lg_1^{-1}}\conj_{g_1}\res^H_{L}(x)
	\]
	for some subgroup $L\leq H$.  If $L<H$ is a proper subgroup then by the same reasoning as the previous lemma we have $\res^H_L(x)\in \mf{p}$.  Thus any multiplicative translate with a proper restriction is in $\mf{p}$, so we need only consider those without a proper restriction. Since $H$ and $gHg^{-1}$ have the same number of elements, a multiplicative translate from $T(G/H)$ to $T(G/gHg^{-1})$ without a restriction also cannot have a norm and thus we need only consider multiplicative translates which are products of conjugations.  Suppose that $y\in T(G/K)$ and $P(y)>n\cdot |G|$.  Then we can write $y$ as a product
	\[
		y = \prod\limits_{i=1}^{A} c_{g_i}(x)
	\]
	where $g_i H g_i^{-1}=K$ and $A>n\cdot |G|$.  The number $A$ is large enough so that some $g_i$ is repeated at least $n$ times. Then $y$ is divisible by $\conj_{g_i}(x)^n = \conj_{g_i}(x^n)\in \mf{p}$.  Thus we can take $M(K) = n\cdot |G|$.
\end{proof}

\begin{lemma}
	\label{lemma: existence of bound case 3}
	Under the \cref{inductive step}, let $K$ be a subgroup of $G$ that properly contains $gHg^{-1}$ for some $g \in G$.
	Then there is an integer $M = M(K) \geq 0$ such that $P(y) > M$ implies that $y \in \mf{p}(G/K)$, for any $y \in T(G/K)$.
\end{lemma}

\begin{proof}
	Since such an $M(K)$ exists if and only if $M(gKg^{-1})$ exists for all $g\in G$ we may assume that $H\leq K$.  First, suppose that $H$ is a maximal subgroup of $K$.  Then every proper subgroup $L<K$ is either conjugate to $H$, or does not contain a conjugate of $H$. By the previous two lemmas, we may suppose we know that $M(L)$ exists for all proper subgroups $L<K$.  We let $U_K$ be the maximum of all the $M(L)$ for $L<K$, let $V_K$ be the number of proper subgroups of $K$, and define $M(K) = U_KV_K$. Let $y\in T(G/K)$ with $P(y)>M(K)$.  Then we can write
	\[
		y = \prod\limits_{i=1}^{N} \nm^{K}_{L_i} \conj_{g_i}\res^H_{g_i^{-1}L_i g_i}(x)
	\]
	for some $N > M(K)$, where $g_i^{-1}L_i g_i\leq H$ and $L_i\leq K$.  Note that because $K$ has strictly more elements than $H$ we must have that $L_i$ is a proper subgroup of $K$.  We can rewrite the product as
	\[
		y = \prod\limits_{L<K} \left(\prod\limits_{L_i=L}\nm^{K}_{L_i} \conj_{g_i}\res^H_{g_i^{-1}L_i g_i}(x)\right) = \prod\limits_{L<K} \nm^K_{L}\left(\prod\limits_{L_i=L} \conj_{g_i}\res^H_{g_i^{-1}L_i g_i}(x)\right)
	\]
	and our choice of $M(K)$ is large enough so that there is some $L$ with $|\{i \mid L_i=L\}|>U_K\geq M(L)$.  It follows that for this $L$ we have
	\[
		\nm^K_{L}\left(\prod\limits_{i,L_i=L} \conj_{g_i}\res^H_{L_i}(x)\right)\in \mf{p},
	\]
	and thus $y \in \mf{p}$.  This proves the case where $H$ is a maximal subgroup of $K$.

	Note that the only place in the above argument where maximality of $H$ in $K$ is used is to know that $M(L)$ exists for all proper $L<K$. To extend the argument to non-maximal subgroups $H < K$, we induct over the maximal length $n$ of chains 
	\[
		H < L_1 < \ldots < L_n < K
	\]
of proper inclusions of subgroups between $H$ and $K$. If $n = 0$, then $H$ is maximal in $K$ and the argument of the previous paragraph holds. If $n = 1$, then every proper subgroup $J$ of $K$ is either 
\begin{enumerate}[(a)]
	\item contained in a conjugate of $H$,
	\item does not contain a conjugate of $H$, or
	\item contains $H$ as a maximal proper subgroup. 
\end{enumerate}
In cases (a) and (b), the previous two lemmas show that $M(J)$ exists. In case (c), we have a maximal chain $H < J < K$ in which $H$ is a maximal proper subgroup of $J$, and the argument of the previous paragraph shows that $M(J)$ exists. Thus, we return to the situation where $M(J)$ exists for all proper subgroups $J < K$, and the argument of the previous paragraph again produces the integer $M(K)$. Proceeding by induction on the length $n$ of a chain proves the claim for all subgroups $K$ containing $H$. 
%
%
\end{proof}

We are now ready to prove \cref{thm: Tambara prime are radical}.

\begin{proof}[Proof of \cref{thm: Tambara prime are radical}]
	The proof proceeds by induction. For the base case, $\mf{p}(G/e)$ is a radical ideal of $T(G/e)$ by \cref{prop: G prime is radical}.

	Now we take on \cref{inductive step}: let $H \leq G$ be a subgroup such that $\mf{p}(G/L)$ is a radical ideal for all $L < H$, and let $x \in T(G/H)$ be such that $x^n \in \mf{p}(G/H)$. We will show $x \in \mf{p}(G/H)$.

	By \cref{lemma: existence of bound case 1,lemma: existence of bound case 2,lemma: existence of bound case 3}, there is an integer $M > 0$ such that $P(y) > M$ implies that $y \in \mf{p}(G/K)$, for $y \in T(G/K)$ for any $K \leq G$; this integer is the maximum over all $K$ of the integers $M(K)$ produced by these lemmas. Then by \cref{lemma: bound implies radical}, $P(y) > 0$ implies that $y \in \mf{p}(G/K)$ for any $K$.

	Finally, $P(x) \geq 1$ since $x$ is a product of one multiplicative translate of itself. Hence, $P(x) > 0$, so $x \in \mf{p}(G/H)$. Thus $\mf{p}(G/H)$ is a radical ideal of $T(G/H)$.
\end{proof}

\section{The Going Up Theorem}

Integral maps play an important role in commutative algebra. We recall that a morphism $f \colon A \to B$ of commutative rings is said to be \emph{integral} if, for all $b \in B$, there exist elements $a_0, \dots, a_{n-1} \in A$ such that
\[b^n = f(a_{n-1}) b^{n-1} + \dots + f(a_1) b + f(a_0),\]
i.e. $b$ is the root of a monic polynomial with coefficients in $A$. Integral maps of commutative rings satisfy the \emph{going up} property, which states that containments of prime ideals in $A$ can be lifted to containments of prime ideals in $B$. Additionally, integral extensions (i.e. injective integral maps) satisfy two additional conditions: \emph{lying over} and \emph{incomparability}. In this section, we will demonstrate that levelwise-integral morphisms of Tambara functors also satisfy going up, and levelwise-integral extensions of Tambara functors also satisfy lying over. It turns out that incomparability does not hold for levelwise-integral extensions of Tambara functors, which will have important implications later.

\begin{definition}\label{defn: going up}
	Given a map of Tambara functors $f\colon S\to T$ we say that $f$ satisfies \emph{going up} if for any chain of prime ideals $\mf{p}_0\subset \mf{p}_1$ in $S$ and prime ideal $\mf{q}_0$ in $T$ such that $f^{-1}(\mf{q}_0)=\mf{p}_0$ there is some $\mf{q}_1$ with $\mf{q}_0\subset \mf{q}_1$ and $f^{-1}(\mf{q}_1) = \mf{p}_1$.
\end{definition}

As noted above, in ordinary commutative algebra an integral map $f\colon R\to S$ always satisfies going up.

\begin{lemma}\label{lemma: going up closed under composition}
	Let $f \colon S \to T$ and $g \colon T \to R$ be morphisms of Tambara functors. If $f$ and $g$ both satisfy going up, then $g \circ f$ satisfies going up.
\end{lemma}

\begin{lemma}\label{lemma: going up for surjections}
	Let $S \to T$ be a levelwise surjective map of Tambara functors. Then $S \to T$ satisfies going up.
\end{lemma}
\begin{proof}
	A levelwise surjection is isomorphic to a quotient map $\pi \colon S \to S/I$. The poset of prime ideals of $S/I$ is isomorphic to the poset of prime ideals of $S$ which contain $I$ \cite[Corollary 4.10]{Nak2012} via $\mf{p} \mapsto \pi^{-1}(\mf{p})$, so $S \to S/I$ satisfies going up.
\end{proof}

Recall that a multiplicative translate of $x \in T(G/H)$ is a term of the form $\nm \conj_g \res(x)$ (\cref{def: multiplicative translate}). Similarly, we call any term of the form
\[
	\tr_L^{L'} \nm_{gKg^{-1}}^L \conj_{K,g} \res^H_K(x)
\]
a \emph{Tambara translate} of $x$.

\begin{lemma}\label{lemma: proposition Q under translates}
	Let $T$ be a Tambara functor, $\mf{p}$ an ideal of $T$, and $x,y$ elements in $T$.  If the proposition $\Q(\mf{p},x,y)$ holds then $\Q(\mf{p},x',y')$ holds where $x'$ and $y'$ are any Tambara translates of $x$ and $y$.
\end{lemma}

\begin{proof}
	By symmetry it suffices to prove the case $y'=y$.  Since every Tambara translate is a composite of transfers, restrictions, norms, and conjugations it suffices to prove each of these cases separately.  All cases except transfer are directly covered by the statement $\Q(\mf{p},x,y)$.  Transfers are covered by the TNR-relations (referring to the notation in \Cref{eq: distinguished bispans} and the relations are discussed after), since a generalized product of a transfer of $x$ with $y$ is a product of a multiplicative translate of $y$ with a transfer of a multiplicative translate of $x$.  By Frobenius reciprocity (\cref{prop: Frobenius reciprocity}) this is a transfer of a generalized product of $x$ and $y$, which is in $\mf{p}$ by the assumption that $\Q(\mf{p},x,y)$ is true.
\end{proof}

\begin{lemma}\label{lemma: max complement is prime}
	Let $S$ and $T$ be Tambara functors. Let $\mf{p}$ be a prime ideal of $S$ and $B$ the levelwise complement of $\mf{p}$ in $S$. If $S\subseteq T$ and $\mf{q}$ is an ideal of $T$ which is maximal with respect to $\mf{q}\cap B = \varnothing$ then $\mf{q}$ is a prime ideal.
\end{lemma}
\begin{proof}
	Let $x\in T(G/H)$ and $y\in T(G/K)$ and suppose we know that all generalized products of $x$ and $y$ are in $\mf{q}$, that is, $\Q(\mf{q},x,y)$ holds. If neither $x$ nor $y$ is in $\mf{q}$ then by maximality of $\mf{q}$ there are subgroups $H'$ and $K'$ of $G$ so that $(B\cap (\mf{q},x))(G/H')\neq \varnothing$ and $(B\cap (\mf{q},y))(G/K')\neq \varnothing$, where $(\mf{q},x)$ is the ideal generated by $\mf{q}$ and $x$.  This implies there are $x'\in T(G/H')$ and $y'\in T(G/K')$, Tambara translates of $x$ and $y$ respectively, and $b_1\in B(G/H')$ and $b_2\in B(G/K')$ so that
	\begin{align*}
		b_1 & = s_1x'+q_1 \\
		b_2 & = s_2y'+q_2
	\end{align*}
	for some $q_i\in \mf{q}$ and $s_i \in S$.  We claim that the proposition $\Q(\mf{q},b_1,b_2)$ holds.  Granting this, since $b_1,b_2\in B\subset S$, we have that the proposition $\Q(\mf{q}\cap S,b_1,b_2) =\Q(\mf{p},b_1,b_2)$ holds. This is impossible, since it would imply that  $b_1$ or $b_2$ is in $B\cap \mf{p}=\varnothing$; with this we will be done by contradiction.

	To see that $\Q(\mf{q},b_1,b_2)$ holds, suppose that $\Phi= \nm_f \conj_g \res_h$ is any multiplicative translate. We will demonstrate that
	\[
		\Phi(b_1) = \Phi(s_1x'+q_1) = \Phi(s_1x')+q_1'
	\]
	where $q_1'$ is some element in $\mf{q}$. This follows when $\Phi$ is a restriction or conjugation since these are ring homomorphisms.  When $\Phi$ is a norm this follows from compatibility of norms and addition as in \cref{proposition: Tambara reciprocity}.  If $\Psi$ is also a multiplicative translate, we have
	\[
		\Phi(b_1)\Psi(b_2) = (\Phi(s_1x')+ q_1')(\Psi(s_2y')+q_2') = \Phi(s_1)\Psi(s_2)\Phi(x')\Psi(y')+q_3
	\]
	where $q_3$ is some element in $\mf{q}$.  By the previous lemma, $\Q(\mf{q},x',y')$ is true, so $\Phi(x')\Psi(y')$ is in $\mf{q}$. Thus $\Phi(b_1)\Psi(b_2)\in \mf{q}$ establishing the claim. This completes the proof.
\end{proof}

\begin{lemma}[{cf. \cite[Theorem 41]{Kaplansky}}]\label{lemma: equivalent condition for going up}
	The following are equivalent for an inclusion of Tambara functors $S \subseteq T$.
	\begin{enumerate}[(a)]
		\item The extension $f\colon S\to T$ satisfies going up.
		\item For any prime $\mf{p}\subset S$, let $\mf{q}$ be an ideal (prime by the last lemma) of $T$ which is maximal with respect to the property that $\mf{q}\cap S\subseteq \mf{p}$, then $\mf{q}\cap S = \mf{p}$.
	\end{enumerate}
\end{lemma}
\begin{proof}
	To see that $(a)\Rightarrow (b)$, suppose that $f$ satisfies going up and let $\mf{p}$ and $\mf{q}$ be as in the statement of $(b)$. Certainly $\mf{q}\cap S$ is a prime ideal of $S$ which is contained in $\mf{p}$.  If this inclusion is strict we can use going up to produce a prime ideal $\mf{q}_1\supset \mf{q}$ and with $\mf{q}_1\cap S=\mf{p}$, which contradicts the maximality of $\mf{q}$.

	To prove $(b)\Rightarrow (a)$ suppose that $\mf{p}_0\subset \mf{p}_1$ are primes in $S$ and $\mf{q}_0$ is a prime in $T$ with $\mf{q}_0\cap S = \mf{p}_0$.  Let $\mf{q}_1$ be maximal among all ideals of $T$ which contain $\mf{q}_0$ and whose intersection with $S$ is contained in $\mf{p}_1$.  Since this ideal satisfies the hypotheses of $(b)$ we must have $\mf{q}_1\cap S = \mf{p}_1$ (this also implies $\mf{q}_0\neq \mf{q}_1$) and so we have proved (a).
\end{proof}

An immediate corollary of characterization (b) is that an extension of Tambara functors which satisfies going up must also satisfy the \emph{lying over} property that the induced maps on $\Spec$ is surjective.

\begin{corollary}[Lying over]\label{cor: lying over}
	If $S\subseteq T$ satisfies going up then the following map is surjective:
	\[
		\Spec(T)\to \Spec(S).
	\]
\end{corollary}

We are now ready to prove a going-up theorem for Tambara functors.  We adapt the proof in Kaplansky's book \cite[Theorem 44]{Kaplansky} for commutative rings. The key step in the proof is that prime ideals are radical, and thus the proof works in this setting because of \cref{thm: Tambara prime are radical}.

\begin{proposition}\label{lemma: going up for extensions}
	Let $S \subseteq T$ be a levelwise integral extension of Tambara functors, i.e. an inclusion of $G$-Tambara functors such that for all $H\leq G$ the ring $T(G/H)$ is integral over $S(G/H)$. Then $S\subseteq T$ satisfies going up.
\end{proposition}
\begin{proof}
	We prove that the Tambara functor inclusion has property $(b)$ from \cref{lemma: equivalent condition for going up}.  Let $\mf{p}$ a prime ideal of $S$ and let $\mf{q}$ be the ideal of $T$ which is maximal with respect to $\mf{q}\cap S\subseteq \mf{p}$.  This ideal is prime by \cref{lemma: max complement is prime}.  Toward a contradiction, suppose there is some $x\in \mf{p}(G/H)\setminus \mf{q}(G/H)$ and let $I$ denote the smallest ideal of $T$ which contains both $\mf{q}$ and $x$.  By the maximality of $\mf{q}$ there must be some $K\leq G$, a $y\in I(G/K)\cap \mf{p}(G/K)$, a Tambara translate of $x$, a $q\in \mf{q}(G/K)$, and a $t\in T(G/K)$ such that $ty+q\in S(G/K)\setminus \mf{p}(G/K)$.

	Since $T(G/K)$ is integral over $S(G/K)$ there is a polynomial relation
	\[
		t^n+s_{n-1}t^{n-1}+\dots+s_1t+s_0=0
	\]
	with $s_i\in S(G/K)$.  Multiplying this relation by $y^n$ we have
	\[
		(ty)^n+s_{n-1}y(ty)^{n-1}+\dots+s_1y^{n-1}(ty)+s_0y^n=0,
	\]
	replacing $ty$ with $ty+q$ we see that
	\[
		(ty+q)^n+s_{n-1}y(ty+q)^{n-1}+\dots+s_1y^{n-1}(ty+q)+s_0y^n\in \mf{q},
	\]
	but this element is also in $S$ since $ty+q$ and $y$ are in $S$.  Thus this sum is in $\mf{q}\cap S \subset \mf{p}$.  Moreover, since $y\in \mf{p}$ it follows that $(ty+q)^n\in \mf{p}$.  But by \cref{thm: Tambara prime are radical} we have that $ty+q\in \mf{p}$, which is a contradiction.  Thus we must have $\mf{q}\cap S = \mf{p}$.
\end{proof}

\begin{thm}[Going up]\label{theorem: going up}
	Let $f \colon S \to T$ be a levelwise integral morphism of Tambara functors. Then $f$ satisfies going up.
\end{thm}
\begin{proof}
	We can factor $f$ as $S \to S/\ker(f) \to T$. The map $S \to S/\ker(f)$ is levelwise surjective, hence satisfies going up by \Cref{lemma: going up for surjections}. The map $S/\ker(f) \to T$ is a levelwise integral extension, hence satisfies going up by \Cref{lemma: going up for extensions}. We are done by \Cref{lemma: going up closed under composition}.
\end{proof}

We obtain the next corollary by combining the going up theorem with the lying over property (\cref{cor: lying over}).

\begin{corollary}\label{cor: levelwise integral extension gives surjection on spec}
	If $f \colon S \to T$ is a levelwise integral extension of Tambara functors then $f^*\colon \Spec(T)\to \Spec(S)$ is surjective.
\end{corollary}

\begin{corollary}\label{cor: levelwise nil integral map gives surjection on spec}
	Let $f \colon S \to T$ be an integral map of Tambara functors.  If
	\[\ker(f) \subseteq \bigcap_{\mf{p} \in \Spec(S)} \mf{p}\]
	then $f^*\colon \Spec(S)\to \Spec(R)$ is surjective.
\end{corollary}
\begin{proof}
	We can factor $f$ as $S \to S/\ker(f) \to T$. Since $S/\ker(f) \to T$ is a levelwise integral extension, $\Spec(T) \to \Spec(S/\ker(f))$ is surjective. Since $\ker(f)$ is contained in every prime ideal of $S$, we also have that $\Spec(S/\ker(f)) \to \Spec(S)$ is a homeomorphism. Thus, $f^* \colon \Spec(T) \to \Spec(S)$ is surjective.
\end{proof}

\begin{remark}
	There is an additional property satisfied by integral extensions in commutative algebra, called \emph{incomparability}. This does \emph{not} hold for integral extensions of Tambara functors, as we will see in \Cref{example: Burnside ghost picture}.
\end{remark}

\section{Fixed Point Functors and GIT Quotients}

Throughout this section, let $R$ be a commutative ring with $G$-action. Recall the definition of the fixed-point Tambara functor $\FP(R)$ from \cref{Ex:FP}.  The goal of this section is to prove \Cref{theorem: homeomorphisms of spec}, which says that $\Spec(\FP(R))$ is homeomorphic to the Zariski spectrum $\Spec(R^G)$. 

\begin{proposition}\label{prop: Spec of FP(R)}
	Let $R$ be a commutative ring with $G$-action.  If $\mf{p}\subseteq \FP(R)$ is a prime ideal then $\mf{p}(G/H) = R^H\cap \mf{p}(G/e)$ for all $H\leq G$.
\end{proposition}
\begin{proof}
	Since the restriction maps of $\FP(R)$ are injective, and Tambara ideals are closed under restrictions, we know that $\mf{p}(G/H) \subseteq \mf{p}(G/e)\cap R^H$ for all $H$.  It remains to prove the reverse inclusion.  Suppose that $x\in \mf{p}(G/e)\cap R^H\subseteq \FP(R)(G/H)$ and note that $\res^H_e(x)\in \mf{p}(G/e)$ so
	\[
		\nm^H_e(\res^H_e(x)) = \prod\limits_{h\in H}h\cdot x = x^{|H|}
	\]
	is an element in $\mf{p}(G/H)$.  By \cref{thm: Tambara prime are radical} this ideal is radical so $x\in \mf{p}(G/H)$ which concludes the proof.
\end{proof}

It follows that the $G$-primes of $R$ and the primes of the fixed point $G$-Tambara functor are in bijection. Recall that $\Spec_G(R)$ denotes the set of all $G$-prime ideals of $R$ (\Cref{defn: G-prime ideals}), and $\Spec(\FP(R))$ is the Nakaoka prime ideal spectrum (\cref{definition: nakaoka spectrum}) of the fixed-point Tambara functor $\FP(R)$.

\begin{corollary}
	\label{cor: SpecG and Spec(FP) poset iso}
	There is a poset isomorphism between $\Spec(\FP(R))$ and $\Spec_G(R)$.
\end{corollary}

This proposition allows us to connect the Tambara spectrum of $\FP(R)$ with the GIT quotient of $\Spec(R)$. We first need a lemma.

\begin{lemma}\label{lemma: intersection of primes is G-prime}
	Let $R$ be a $G$-ring and let $\mf{p}\in \Spec(R)$. The ideal
	\[
		\overline{\mf{p}} = \bigcap\limits_{g\in G} g\cdot \mf{p}
	\]
	is a $G$-prime ideal.
\end{lemma}
\begin{proof}
	Since the action of $G$ permutes the terms in the intersection we see that $\overline{\mf{p}}$ is $G$-invariant.  Suppose that $x,y\in R$ and $x(gy)\in \overline{\mf{p}}$ for all $g\in G$. Towards a contradiction, suppose that neither $x$ nor $y$ is in $\overline{\mf{p}}$.  Then there must be some $g,h\in G$ such that $x\notin g\mf{p}$ and $y\notin h\mf{p}$.  Equivalently, we have that neither $g^{-1}x$ nor $h^{-1}y$ is in $\mf{p}$.  On the other hand, we have
	\[
		g\cdot((g^{-1}x)(h^{-1}y) = x(gh^{-1}y)\in \overline{\mf{p}}
	\]
	by assumption.  Since $\overline{\mf{p}}$ is $G$-invariant we have $(g^{-1}x)(h^{-1}y)\in \overline{\mf{p}}\subset \mf{p}$.  This is impossible because $\mf{p}$ is a prime ideal, so we have our contradiction.
\end{proof}

If $R$ is a $G$-ring we write $\Spec_G(R)$ for the poset of $G$-prime ideals. We write $\Spec(R)/G$ for the orbits of $\Spec(R)$ under the $G$-action. Recall that the \emph{geometric invariant theory (GIT) quotient} of an affine scheme $\Spec(R)$ with respect to an action of $G$ on $R$ is the scheme $\Spec(R^G)$, where $R^G$ is the $G$-fixed points of $R$ \cite[Chapter 1.2]{MumfordFogarty}

\begin{proposition}\label{theorem: all the equispecs are the same}
	There are poset isomorphisms between $\Spec(R)/G$, $\Spec_G(R)$, $\Spec(R^G)$ and $\Spec(\FP(R))$.
\end{proposition}

\begin{proof}
	The isomorphism of posets between $\Spec_G(R)$ and $\Spec(\FP(R))$ is established by \cref{cor: SpecG and Spec(FP) poset iso} so it remains to establish poset isomorphisms between $\Spec(R)/G$, $\Spec_G(R)$, and $\Spec(R^G)$.  If $\mf{p}\in \Spec(R)$ we write $[\mf{p}]$ for the corresponding element in the orbits $\Spec(R)/G$.  We define
	\begin{align*}
		\phi & \colon \Spec(R)/G\to \Spec_G(R) \\
		\psi & \colon \Spec_G(R)\to \Spec(R^G)
	\end{align*}
	by $\phi([\mf{p}]) = \overline{\mf{p}}$, in the notation of \cref{lemma: intersection of primes is G-prime}, and $\psi(\mf{q}) = \mf{q}^G$. Note that both $\phi$ and $\psi$ preserve inclusions.

	Consider the diagram
	\[
		\begin{tikzcd}
			\Spec(R) \ar[d,"i^*"'] \ar[r] & \Spec(R)/G \ar[d,"\phi"] \\
			\Spec(R^G)  & \Spec_G(R) \ar[l,"\psi"]
		\end{tikzcd}
	\]
	where the unlabeled arrow is the canonical quotient and $i\colon R^G\to R$ is the inclusion of fixed points. Since $i$ is an integral extension, $i^*\colon \Spec(R) \to \Spec(R^G)$ is surjective and sends $\mf{p}$ to $R^G \cap \mf{p} = \mf{p}^{G}$. In the other direction on the diagram, $\mf{p}$ is sent to
	\[
		\big(\bigcap_{g \in G} g \cdot \mf{p}\big)^G = \bigcap_{g \in G} (g \cdot \mf{p})^G = \bigcap_{g \in G} (\mf{p})^G = (\mf{p})^G.
	\]
	Therefore this diagram is commutative. Since $i^*$ is surjective, then $\psi$ must also be surjective.  By \cref{cor: G prime determined by fixed points} we know that $\psi$ is also injective so $\psi$ is an order-preserving bijection. Since $i^*$ is a surjection, and $\psi$ is a bijection, then $\phi$ must also be a surjection. Last, we need to show that $\phi$ is injective.

	To see that $\phi$ is injective, suppose that $\mf{p},\mf{q}\in \Spec(R)$ with $\overline{\mf{p}} = \overline{\mf{q}}$.  If $x\in \mf{p}$ then
	\[
		\nm(x) = \prod_{g \in G} gx \in \overline{\mf{p}} = \overline{\mf{q}}.
	\]
	It follows that there is some $g\in G$ so that $gx\in \mf{q}$ or, equivalently, $x\in g^{-1}\mf{q}$.  It follows that
	\[
		\mf{p}\subseteq \bigcup_{g\in G} g\cdot \mf{q}
	\]
	so by the prime avoidance lemma \cite[Lemma 3.3]{Eisenbud} we have $\mf{p}\subseteq g\cdot \mf{q}$ for some $g$.   Repeating the argument we see that $\mf{q}\subseteq h\mf{p}$ for some $h\in G$. So we have $\mf{p}\subseteq g\mf{q}\subseteq gh\mf{p}$. If we show that $(gh)\mf{p} = \mf{p}$ we are done since this gives $g\mf{q} = \mf{p}$ so these two prime ideals are in the same orbit.  To prove that  $(gh)\mf{p} = \mf{p}$, let $n$ be the order of $gh$ and observe that
	\[
		\mf{p}\subseteq (gh)\mf{p}\subseteq (gh)^2\mf{p}\subseteq\dots\subseteq (gh)^{n}\mf{p} =\mf{p}
	\]
	and so $(gh)\mf{p} =\mf{p}$.
\end{proof}

\begin{remark}
	When $R$ has a (left) action by a finite group $G$ we obtain an associated (right) action on  the affine scheme $\Spec(R)$ and it is natural to ask about the categorical quotient $\Spec(R)/\!/G$, defined as the colimit of the functor $BG\to \mathsf{AffScheme}$ which realizes the action.  Since affine schemes are the opposite of the category of commutative rings this is canonically isomorphic to the spectrum of the limit of the functor $BG\to \mathsf{CRing}$ which realizes the $G$-action on $R$.  In particular, $\Spec(R)/\!/G\cong \Spec(R^G)$.  Moreover, the above proposition tells us there is an isomorphism of posets $\Spec(R)/\!/G\cong \Spec(R)/G$, where the latter denotes the usual quotient space.

	The categorical quotient  $\Spec(R)/\!/G$ is known as the \emph{GIT quotient}.  We observe that the result above gives a new proof of the classical fact that the GIT quotient of an affine scheme by a finite group exists and is the actual quotient in spaces \cite[Chapter 1.2]{MumfordFogarty}
\end{remark}

We can upgrade the poset isomorphism above to a homeomorphism.  For this we first need a technical lemma.

\begin{lemma}\label{lemma: a lemma for homeomorphism}
	Let $R$ be a $G$-ring and let $\mf{p}\in \mathrm{Spec}(R^G)$ be any prime ideal.  Let $\overline{\mf{q}}\subseteq R$ be any $G$-prime ideal with $\overline{\mf{q}}^G = \mf{p}$.  Then for any $G$-invariant ideal $I\subseteq R$  we have $I^G\subseteq \mf{p}$ if and only if $I\subseteq \overline{\mf{p}}$.
\end{lemma}
\begin{proof}
	The backward implication is immediate from taking fixed points so it remains to check the forward implication.  Note that by \cref{theorem: all the equispecs are the same} the $G$-prime ideal $\overline{\mf{q}}$ is actually unique.  Unpacking the proof it can be obtained by picking a prime ideal $\mf{q}\in \mathrm{Spec}(R)$ with $\mf{q}^G = \mf{p}$ and taking
	\[
		\overline{\mf{q}} = \bigcap\limits_{g\in G}g\cdot \mf{q}.
	\]
	While the ideal $\mf{q}$ is not unique, its $G$-orbit is unique, again by \cref{theorem: all the equispecs are the same}.

	Suppose that $I$ is some $G$-invariant ideal with $I^G\subseteq \mf{p}$.  Then $I\subseteq \mf{q}$ for some ideal $\mf{q}$ which is maximal with respect to $\mf{q}^G\subseteq \mf{p}$.   Since the extension $R^G\to R$ is integral, it follows from \cite[Theorem 41]{Kaplansky} that any such $\mf{q}$ is prime and $\mf{q}^G = \mf{p}$.  Since $I$ is $G$-invariant we have
	\[
		I\subseteq \bigcap\limits_{g\in G}g\cdot \mf{q} = \overline{\mf{q}}. \qedhere
	\]
\end{proof}

\begin{thm}
	\label{theorem: homeomorphisms of spec}
	The bijection of \cref{theorem: all the equispecs are the same} gives a homeomorphism
	\[
		\Spec(R^G) \cong \Spec(\FP(R))
	\]
	where the left space has the Zariski topology and the right space has the topology from \cref{definition: nakaoka spectrum}.
\end{thm}

\begin{proof}
	Recall that the bijection is induced by the function
	\[
		\alpha \colon \Spec(\FP(R))\to \Spec(R^G)
	\]
	defined by $\alpha(\mf{q}) = \mf{q}(G/e)^G=\mf{q}(G/G)$ where the second equality is \cref{prop: Spec of FP(R)}.  The inverse
	\[
		\beta\colon \Spec(R^G)\to \Spec(\FP(R))
	\]
	sends $\mf{p}\in \Spec(R^G)$ to $\FP(\overline{\mf{p}})$, where $\overline{\mf{p}}$ is a (necessarily unique) $G$-prime ideal of $R$ with $\overline{\mf{p}}^G = \mf{p}$.

	To see that $\alpha$ is continuous we show that the preimage of closed sets is closed.  Any closed subset of $\Spec(R^G)$ has the form
	\[
		V(I) = \{ \mf{p}\in \Spec(R^G)\mid I\subseteq \mf{p}\}
	\]
	where $I\subseteq R^G$ is some ideal.  We will write $\FP(I)\subseteq \FP(R)$ for the constant ideal at $I$.

	We claim that the preimage of $V(I)$ under $\alpha$ is precisely
	\[
		V(\FP(I)) = \{\mf{q}\in \Spec(\FP(R))\mid \FP(I)\subseteq \mf{q} \}
	\]
	which is closed.  To see this, note that $\mf{q}\in \alpha^{-1}(V(I))$ if and only if $I\subseteq \mf{q}(G/G)$.  But by \cref{prop: Spec of FP(R)} we have $I\subseteq \mf{q}(G/G)$ if and only if $\FP(I)\subseteq \mf{q}$ and  thus $\alpha^{-1}(V(I)) = V(\FP(I))$.

	To see that $\beta$ is continuous fix an ideal $I\subseteq \FP(R)$ so that an arbitrary closed set in $\Spec(\FP(R))$ has the form
	\[
		V(I)  = \{\mf{q}\in \Spec(\FP(R))\mid I\subseteq \mf{q} \}.
	\]
	By definition we have
	\[
		\beta^{-1}(V(I)) = \{\mf{p}\in \Spec(R^G)\mid I\subseteq \FP(\overline{\mf{p}}) \}.
	\]
	By \cref{prop: Spec of FP(R)}, we see that $I\subseteq \FP(\overline{\mf{p}})$ if and only if $I(G/e)\subseteq \overline{\mf{p}}$.
	By \cref{lemma: a lemma for homeomorphism} we have that $I(G/e)\subseteq \overline{\mf{p}}$ if and only if $I(G/e)^G\subseteq \mf{p}$ and so $\beta^{-1}(V(I))$ is equal to
	\[
		V(I(G/e)^G) = \{\mf{p}\in \Spec(R^G)\mid I(G/e)^G\subset\mf{p} \}
	\]
	which is closed in $\Spec(R^G)$ by definition of the topology.
\end{proof}

\section{The Ghost of a \texorpdfstring{$C_p$}{Cₚ}-Tambara Functor}
\label{section: ghost}

In this section we specialize to the case $G= C_p$ and define the \emph{ghost} of a $C_p$-Tambara functor $T$, building on ideas of Th\'evenaz \cite{Thevenaz:SomeRemarks} and Calle--Ginnett \cite{CG2023}.  Let $\gamma$ be a generator of $C_p$.

\begin{notation}
	\label{notation:underlying/fixed}
	For a $C_p$-Tambara functor $T$, we refer to $T(C_p/e)$ as the \emph{underlying level of $T$} and $T(C_p/C_p)$ as the \emph{fixed level of $T$}. Similarly, for a morphism $\varphi \colon S \to T$ of $C_p$-Tambara functors, we say that $\varphi(C_p/e) \colon S(C_p/e) \to T(C_p/e)$ is the \emph{underlying level of $\varphi$} and $\varphi(C_p/C_p) \colon S(C_p/C_p) \to T(C_p/C_p)$ is the \emph{fixed level of $\varphi$}. This terminology comes from fixed point Tambara functors, where the fixed level consists entirely of fixed points and the underlying level is the underlying $C_p$-ring from which the fixed points are taken. This also informs the terminology of \cref{Ex:AxGG,Ex:AxGe}.
\end{notation}

\begin{notation}\label{defn: geometric fixed points}  For $T$ a $C_p$-Tambara functor, we use the following notation.
	\begin{enumerate}[(a)]
		\item Let $\tau$ denote the image of the transfer
		      \[
			      \tr \colon T(C_p/e) \to T(C_p/C_p).
		      \]
		      This is an ideal of $T(C_p/C_p)$ by Frobenius reciprocity (\Cref{prop: Frobenius reciprocity}).
		\item Let $\Phi^{C_p} T \coloneqq T(C_p/C_p)/\tau$ denote the \emph{geometric fixed points} of $T$ \cite[Def. 5.10]{BGHL2019}.
	\end{enumerate}
\end{notation}

We note that the geometric fixed points has been called the \emph{Brauer quotient} in work of Th\'evenaz and others \cite{Thevenaz:SomeRemarks}. The term geometric fixed points comes from the viewpoint of equivariant stable homotopy theory \cite{MandellMay,BGHL2019}.

\begin{definition}\label{def: ghost}
	Let $T$ be a $C_p$-Tambara functor.  We define the $\emph{ghost}$ of $T$ (denoted $\ghost(T)$) to be the $C_p$-Tambara functor defined by the Lewis diagram
	\[
		\begin{tikzcd}[row sep=huge]
			T(C_p/e)^{C_p} \times \Phi^{C_p} T
			\ar[d, "{\res}" description]
			\\
			T(C_p/e)
			\ar[u, bend right=50, "\nm"{right}]
			\ar[u, bend left=50, "\tr"{left}]
			\arrow[from=2-1, to=2-1, loop, in=300, out=240, distance=5mm, "\gamma"']
		\end{tikzcd}
	\]
	where the restriction, transfer, and norm are given by
	\[
		\res(x,y) = x,  \quad
		\tr(x) = \left(\sum\limits_{g\in G}g\cdot x,0\right), \quad \text{and} \quad
		\nm(x) = \left(\prod\limits_{g\in G}g\cdot x, \nm(x)+\tau \right).
	\]
\end{definition}

\begin{remark}
	Constructions similar to our ghost construction have previously appeared in work of others. The ghost of the Burnside Tambara functor was defined in \cite[Def.~3.3]{CG2023}. It also generalizes work of Th\'evenaz, who introduced a similar construction for Green functors\footnote{A \emph{Green functor} is like a Tambara functor, but without the norms. Every Tambara functor has an underlying Green functor, but not all Green functors extend to Tambara functors.} \cite[Sec.~4]{Thevenaz:SomeRemarks}. The outline of the ghost construction is also present in \cite[Sec.~19]{Str2012} and in the background of \cite{Bru2005}. A related but different construction also appears in \cite{Sulyma2023}; see \cref{remark:transversal} for more details.
\end{remark}

The ghost of $T$ is certainly a commutative ring at every level.  To check it is a Tambara functor one needs to check that the TNR-relations, meaning the notation used in \Cref{eq: distinguished bispans} along with the relations discussed after this equation, are satisfied.  This is straightforward here because the restriction map is a projection; we leave the details to the reader.

\begin{proposition}
	For a $C_p$-Tambara functor $T$ the ghost $\ghost(T)$ is a $C_p$-Tambara functor.
\end{proposition}

The name for the ghost construction is motivated by the following example.

\begin{example}
	\label{example: Burnside ghost}
	The ghost of the Burnside $G$-Tambara functor $\uA$ is the $G$-Tambara functor which is at level $G/H$ is the \emph{ghost ring} of the Burnside ring $A(H)$. By definition, the ghost ring of $A(H)$ is $\prod_{[K]} \Z$ where $[K]$ ranges over the $H$-conjugacy classes of subgroups of $H$. This ring was known to and used by Burnside~\cite{burnside_table_of_marks}, although this definition was not introduced in this way.
\end{example}

For $G = C_p$, we have only two levels, and the ghost of $\uA$ is given by $\Z$ at level $C_p/e$ and $\Z \times \Z$ at level $C_p/C_p$.

\begin{example}
	\label{example: ghost not inj on spec}
	The ghost of the Burnside $C_p$-Tambara functor $\uA$ is given by the Lewis diagram
	\[
		\ghost(\uA) =
		\begin{tikzcd}[row sep=huge]
			\mathbb{Z} \times \mathbb{Z}
			\ar[d, "\res" description]
			\\
			\mathbb{Z}
			\ar[u, bend right=50, "\nm"{right}]
			\ar[u, bend left=50, "\tr"{left}]
			\arrow[from=2-1, to=2-1, loop, in=300, out=240, distance=5mm, "{\mathrm{trivial}}"']
		\end{tikzcd}
	\]
	where $\res(x,y) = x$, $\tr(x) = (px,0)$, and $\nm(x) = (x^2,x)$.
	This is the ghost Tambara functor defined in \cite[Def. 2.3]{CG2023}.
\end{example}

\begin{example}\label{Ex:Ghost of AxGe}
	We describe $\ghost(T)$ where $T = \uA[C_p/e]$ is the free $C_p$-Tambara functor on an underlying generator. By \cref{Ex:AxGe}, we have
	\[
		T(C_p/e) = \Z[x_0,\ldots,x_{p-1}], \, \text{and}
	\]
	\[
		T(C_p/C_p) = \Z[n][\{ t_{\vec{v}} \}_{\vec{v} \in \N^{\times p}}] / I.
	\]
	The transfer ideal $\tau$ is $\langle \{t_{\vec{v}}\}_{\vec{v} \in \N^{\times p}} \rangle$, and we have $I \subseteq \tau$, so $\Phi^{C_p} T \cong \Z[n].$
	Since $C_p$ acts on $T(C_p/e) = \Z[x_0,\ldots,x_{p-1}]$ by cyclic permutation of the variables, we have that
	$\Z[x_0,\ldots,x_{p-1}]^{C_p}$ is the subring of cyclic polynomials. When $p=2$ this coincides with the polynomial ring on elementary symmetric polynomials, and when $p > 2$ it is not polynomial (indeed, it is not even regular after base change to $\mathbb{C}$).
	Thus $\ghost(T)$ has the form
	\[
		\begin{tikzcd}[row sep=huge]
			\Z[x_0,\ldots,x_{p-1}]^{C_p} \times \Z[n]
			\ar[d, "\res" description]
			\\
			\Z[x_0,\ldots,x_{p-1}]
			\ar[u, bend right=50, "\nm"{right}]
			\ar[u, bend left=50, "\tr"{left}]
			\arrow[from=2-1, to=2-1, loop, in=300, out=240, distance=5mm, "\gamma"']
		\end{tikzcd}
	\]
	where $\res(g(x_0,\ldots,x_{p-1}), h(n)) = g(x_0,\ldots,x_{p-1})$, and $\tr$ and $\nm$ are given by
	\begin{align*}
		\tr(f(x_0,\dots,x_{p-1})) & = \left(\sum_{i=0}^{p-1}f(x_i,x_{i+1},\dots,x_{p-1},x_0,\dots,x_{i-1}),0\right)                \\
		\nm(f(x_0,\dots,x_{p-1})) & = \left(\prod_{i=0}^{p-1}f(x_i,x_{i+1},\dots,x_{p-1},x_0,\dots,x_{i-1}),f(n,n,\dots,n)\right).
	\end{align*}
\end{example}

\begin{example}\label{example: ghost of fixed point polynomial functor}
	We describe $\ghost(S)$ where $S = \uA[C_p/C_p]$ is the free Tambara functor on a fixed generator. By \cref{Ex:AxGG}, we have
	$$S(C_p/e) = \Z[x],$$
	$$S(C_p/C_p) = \Z[x,t,n]/\langle t^2-pt, tn-tx^p\rangle.$$
	The transfer ideal $\tau$ is given by $\langle t \rangle$, so $\Phi^{C_p} S \cong \Z[x,n].$
	Since $C_p$ acts on $S(C_p/e) = \Z[x]$ trivially, we have $S(C_p/e)^{C_p} = \Z[x]$.
	Thus $\ghost(S)$ has the form
	\[
		\begin{tikzcd}[row sep=huge]
			\Z[x] \times \Z[x,n]
			\ar[d, "{\res}" description]
			\\
			\Z[x]
			\ar[u, bend right=50, "\nm"{right}]
			\ar[u, bend left=50, "\tr"{left}]
			\arrow[from=2-1, to=2-1, loop, in=300, out=240, distance=5mm, "\gamma"']
		\end{tikzcd}
	\]
	where $\res$ is projection onto the first component, and $\tr$ and $\nm$ are given by
	\[
		\tr(f(x)) = (pf(x),0),\quad \mathrm{and}\quad \nm(f(x)) = (f(x)^p,f(n)).
	\]
\end{example}

\begin{remark}\label{remark:transversal}
	The ghost functor we consider here is similar in spirit 
	to the \emph{transversal Tambara functors} considered by Sulyma \cite{Sulyma2023}. A $C_p$-Tambara functor $T$ is \emph{transversal} when $T(C_p/C_p)$ is the pullback of $\Phi^{C_p}T$ and $T(C_p/e)^{C_p}$ along the maps to the cokernel of
	\[
		\res \circ \tr \colon T(C_p/e) \to T(C_p/e)^{C_p}.
	\]
	See \cite[Example 1.4]{Sulyma2023} for details. Transversal Tambara functors are designed so that it is easy to check the Tambara reciprocity relations (\cref{proposition: Tambara reciprocity}) and build a Tambara functor out of commutative rings and homomorphisms.

	One can check that the ghost functors we consider are transversal, and therefore enjoy the nice properties of transversal ones. In fact, the transversalization $\operatorname{\pitchfork} T$ provides an epi-mono factorization of the ghost map
	\[
		T \twoheadrightarrow \operatorname{\pitchfork}(T) \hookrightarrow \ghost(T)
	\]
	for any $C_p$-Tambara functor $T$. Note that the transversalization $\pitchfork$ is an idempotent operation, while the ghost $\ghost$ is not.
\end{remark}

\subsection{A Character Homomorphism}
The Burnside ring of $G$ is related to the ghost ring by means of a character homomorphism
\[
	\chi \colon A(G) \to \prod_{[H]} \Z,
\]
sending a finite $G$-set $X$ to the tuple $(|X^H|)_{[H]\leq G}$ of cardinalities of $H$-fixed points, as $H$ ranges over conjugacy classes of subgroups of $G$. Importantly, the character map is injective, and becomes an isomorphism after inverting $|G|$ on both the left and the right. These facts allow us to compute in the ghost ring which has a much simpler multiplicative structure than the Burnside ring.

\cref{example: Burnside ghost} turns the above into the prototypical example of the ghost of a $C_p$-Tambara functor. Below we define a character map $T \to \ghost(T)$ for any $C_p$-Tambara functor $T$, and show it is similarly injective when $T$ is $p$-torsion free, and an isomorphism when $p$ is invertible in $T$.

\begin{definition}
	The \emph{ghost map} is the map of Tambara functors
	\[
		\ghostmap_T \colon T\to \ghost(T)
	\]
	which is the identity on the underlying level, and given on the fixed level by
	\[
		(\ghostmap_T)_{C_p/C_p} = (\res,q)\colon T(C_p/C_p)\to T(C_p/e)^{C_p}\times \Phi^{C_p} T
	\]
	where $q \colon T(C_p/C_p) \to \Phi^{C_p} T$ is the quotient map.
\end{definition}

\begin{example}
	For the Burnside Tambara functor $\uA$, the map
	\[
		(\ghostmap_{\uA})_{C_p/C_p} \colon A(C_p) \cong \uA(C_p/C_p) \to \uA(C_p/e)^{C_p} \times \Phi^{C_p} \uA \cong \Z \times \Z
	\]
	from $A(C_p)$ to its ghost ring is the character map described above.  This map, called the ghost map, is used by Dress in \cite{Dress1971} to study $\Spec(A(G))$ for finite groups $G$, and is an essential tool in the study of Burnside rings. Calle and Ginnet use ghost maps to compute $\Spec(\uA)$ for $G$ a finite cyclic group in \cite{CG2023}.
\end{example}

The more general ghost map we consider here is studied in work of Th{\'e}venaz \cite{Thevenaz:SomeRemarks}, albeit in the context of Green functors instead of Tambara functors.  Notably, Thevenaz shows that the kernel and cokernel of this map consist entirely of $p$-torsion \cite[Theorem 3.2]{Thevenaz:SomeRemarks}, and moreover that the kernel is levelwise nilpotent; similar results hold for $C_p$-Tambara functors. The proofs extend to $C_p$-Tambara functors to give the following two lemmas:

\begin{lemma}\label{lemma: kernel of ghost is p-torsion}
	Let $T$ be a $C_p$-Tambara functor, and let $K = \ker (\ghostmap_T)$. Then $pK = 0$ (levelwise).
\end{lemma}

\begin{lemma}
	Let $T$ be a $C_p$-Tambara functor, and suppose $p$ is invertible in $T(C_p/C_p)$. Then $\ghostmap_T$ is an isomorphism.
\end{lemma}

Combining the previous two statements, we see that the ghost map is a reasonable notion of a character map for $C_p$-Tambara functors.

\begin{proposition}[cf. {\cite[Theorem 3.2]{Thevenaz:SomeRemarks}}]
	\label{cor: ghost injective}
	Let $T$ be a $C_p$-Tambara functor.
	\begin{enumerate}[(a)]
		\item If $T(C_p/C_p)$ is $p$-torsion free, then the ghost map is injective.
		\item If $p$ is invertible in $T(C_p/C_p)$, then the ghost map is an isomorphism.
	\end{enumerate}
\end{proposition}

We close this section with a fact about the kernel of the ghost map.

\begin{lemma}
	Let $T$ be a $C_p$-Tambara functor, and let $K = \ker (\ghostmap_T)$. Then $K^2 = 0$ (levelwise).
\end{lemma}
\begin{proof}
	Since $\ghostmap_T$ is an isomorphism on the underlying level, we have $K(C_p/e) = 0$. Thus, it suffices to show that $K(C_p/C_p)^2 = 0$. So, let $a,b \in K(C_p/C_p)$ be arbitrary. By definition of $\ghostmap_T$, we have that $\res(a) = 0$ and $b = \tr(c)$ for some $c \in T(C_p/e)$. By Frobenius reciprocity, we obtain
	\[ab = a\tr(c) = \tr(\res(a) c) = \tr(0) = 0,\]
	as desired.
\end{proof}

Since Tambara prime ideals are levelwise radical by \Cref{thm: Tambara prime are radical}, we conclude:

\begin{corollary}\label{cor: kernel of ghost is nil}
	Let $T$ be a $C_p$-Tambara functor, and let $K = \ker (\ghostmap_T)$. Then $K$ is contained in every prime ideal of $T$.
\end{corollary}

\subsection{Functoriality}

We make the observation in this subsection that the ghost construction is a functor $\Tamb_G \to \Tamb_G$, and the ghost map is a natural transformation $\id \to \ghost$. 

\begin{proposition}
	The ghost construction extends to a functor $\ghost \colon \Tamb_{C_p} \to \Tamb_{C_p}$. Given a morphism $\varphi \colon S \to T$ of $C_p$-Tambara functors, the underlying level $\ghost(\varphi)(C_p/e)$ is given by $\varphi(C_p/e)$ and the fixed level $\ghost(\varphi)(C_p/C_p)$ is given by $\varphi(C_p/C_p) \times \Phi^{C_p}(\varphi)$.
\end{proposition}

\begin{proposition}
	Let $T$ be a $C_p$-Tambara functor. The ghost map $\ghostmap_{T} \colon T \to \ghost(T)$ is natural in $T$.
\end{proposition}
\begin{proof}
	Let $\varphi \colon S \to T$ be a morphism of Tambara functors. We wish to show that
	\[
		\begin{tikzcd}
			S
			\ar[r, "{\ghostmap_{S}}"]
			\ar[d, "\varphi"']
			&
			\ghost(S)
			\ar[d, "{\ghost(\varphi)}"]
			\\
			T
			\ar[r, "{\ghostmap_{T}}"']
			&
			\ghost(T)
		\end{tikzcd}
	\]
	commutes. At the underlying level, $\ghostmap_{S}$ and $\ghostmap_{T}$ are identities, and $\ghost(\varphi)$ agrees with $\varphi$ so the diagram commutes. At the fixed level, we have
	\[
		\begin{tikzcd}[sep = large]
			S(C_p/C_p)
			\ar[r]
			\ar[d, "{\varphi_{C_p/C_p}}"]
			&
			S(C_p/e)^{C_p} \times S(C_p/C_p)/\tau
			\ar[d, "{\varphi_{C_p/e} \times \Phi^{C_p} \varphi}"]
			\\
			T(C_p/C_p)
			\ar[r]
			&
			T(C_p/e)^{C_p} \times T(C_p/C_p)/\tau
		\end{tikzcd}
	\]
	Going around the top-right gives
	\[x \mapsto (\res(x), x+\tau) \mapsto (\varphi_{C_p/e}(\res(x)), \varphi_{C_p/C_p}(x)+\tau),\]
	and going around the bottom-left gives
	\[x \mapsto \varphi_{C_p/C_p}(x) \mapsto (\res(\varphi_{C_p/C_p}(x)), \varphi_{C_p/C_p}(x) + \tau).\]
	These maps agree because $\varphi_{C_p/e} \circ \res = \res \circ \varphi_{C_p/C_p}$.
\end{proof}

\subsection{The Ghost and the Spectrum}
\label{section: the ghost and the spectrum}

For the purpose of computing $\Spec(T)$ we are interested in the function
\[
	\ghostmap^*\colon \Spec(\ghost(T))\to \Spec(T)
\]
induced on $\Spec$ by the ghost map.  It turns out this map is \emph{always} surjective.

\begin{thm}\label{Thm:GhostSurj}
	If  $T$ is any $C_p$-Tambara functor then the ghost map
	\[
		\ghostmap_{T}\colon T\to \ghost(T)
	\]
	induces a surjection $\Spec(\ghost(T))\to \Spec(T)$.
\end{thm}
\begin{proof}
	By \Cref{cor: levelwise nil integral map gives surjection on spec} and \Cref{cor: kernel of ghost is nil}, it suffices to prove that the ghost map is levelwise integral.  Since the map on the underlying is the identity, the only thing to check is that the ring map
	\[
		T(C_p/C_p)\xrightarrow{(\res,q)} T(C_p/e)^{C_p}\times \Phi^{C_p} T
	\]
	is integral.  Note that the quotient map $q$ is surjective and is therefore integral.  Furthermore, the map $\res \colon T(C_p/C_p)\to T(C_p/e)^{C_p}$ is integral because any fixed point $a$ satisfies the equation
	\[
		a^p- \res(\nm(a))=0.
	\]

	We are then done because the product of integral ring maps is always integral, which we will demonstrate now.  If $f = (f_1, f_2) \colon R\to A\times B$ is the product of integral ring maps and $(a,b)\in A\times B$ then pick a monic polynomial $p_2(x)$ with coefficients in $R$ such that $p_2(b)=0$.  Then we have $p_2(a,b) = (a',0)$ for some $a'\in A$.  Since $f_1$ is also integral, we pick a further monic polynomial $p_1(x)$ so that $p_1(a')=0$ and $p_1(0) = 0$. Then we have $(p_1\circ p_2)(a,b)=0$, as desired.
\end{proof}

In the proof above, we showed that $\ghostmap_T$ is levelwise integral.  This, combined with \cref{theorem: going up} and \cref{cor: lying over}, gives us the following corollary.

\begin{corollary}\label{cor: ghost has going up and lying over}
	For any $C_p$-Tambara functor $T$, $\ghostmap_T$ satisfies both going up and lying over.
\end{corollary}

In light of the theorem above we will be interested in understanding the $C_p$-Tambara prime ideals of ghosts.  Suppose that $\mf{a}\subseteq T(C_p/e)$ is a $C_p$-prime ideal and suppose that $\mf{b}\subseteq \Phi^{C_p} T$ is either a prime ideal or the entire ring.  Assume further that we have
\[
	\nm(\mf{a})+\tau\subseteq \mf{b}
\]
where $\nm$ is the norm operation of $T$. It is then easy to check that $\ghost(T)$ has a proper ideal $(\mf{a};\mf{b})$ given levelwise by
\begin{align*}
	(\mf{a};\mf{b})(C_p/e)   & = \mf{a}                     \\
	(\mf{a};\mf{b})(C_p/C_p) & = \mf{a}^{C_p}\times \mf{b}.
\end{align*}
This ideal is not prime in general, but we can classify exactly when it is a prime ideal.

\begin{lemma}\label{lemma: when ideal pairs are prime in ghost}
	Let $T$ be a $C_p$-Tambara functor. The ideal $(\mf{a};\mf{b})$ is a prime ideal of $\ghost(T)$ if and only if precisely one of the following holds:
	\begin{enumerate}[(a)]
		\item $\mf{b} = \Phi^{C_p} T$, or
		\item $\mf{a} = \nm^{-1}(\mf{b})$, where $\nm\colon T(C_p/e)\to \Phi^{C_p} T$ is the ring map induced by the norm.
	\end{enumerate}
\end{lemma}
\begin{proof}
	If both (a) and (b) hold, then $(\mf{a};\mf{b})$ is the whole Tambara functor $\ghost(T)$, which is not prime by definition.

	If neither of the two conditions above holds, pick an element $x\in \nm^{-1}(\mf{b})\setminus \mf{a}$.  Consider the elements $(0,1)\in \ghost(T)(C_p/C_p)$ and $x\in \ghost(T)(C_p/e)$. Neither of these elements is in $(\mf{a};\mf{b})$.  However, a straightforward computation shows that $\Q((\mf{a};\mf{b}),(0,1),(x,\nm(x))$ is satisfied and thus the ideal $(\mf{a};\mf{b})$ is not prime.  Thus we have shown that $(\mf{a};\mf{b})$ being prime implies that precisely one of (a) or (b) holds.

	To prove that all ideals of the form $(\mf{a}; \Phi^{C_p} T)$ are prime, note that there is a map of Tambara functors $\epsilon\colon \ghost(T)\to \FP(T(C_p/e))$ which is the identity at the underlying level and the restriction map at the fixed level.  The ideal $$(\mf{a}; \Phi^{C_p} T)$$ is the $\epsilon$-preimage of the ideal of $\FP(T(C_p/e))$ which corresponds to the ideal $\mf{a}$.  Since the preimage of a prime ideal is prime, we see that condition (a) implies that $(\mf{a};\mf{b})$ is prime.

	It remains to show that ideals of the form $(\nm^{-1}(\mf{b});\mf{b})$ are prime when $\mf{b}\subset T(C_p/C_p)/\tau$ is a proper prime ideal.   It is enough to show that for any $x$ and $y$ in $T$ the proposition $\Q((\nm^{-1}(\mf{b});\mf{b}),x,y)$ being satisfied implies that $x$ or $y$ is in $(\nm^{-1}(\mf{b});\mf{b})$.

	Suppose first that $x,y\in \ghost(T)(C_p/e)$ and $\Q((\nm^{-1}(\mf{b});\mf{b}),x,y)$ is true.  Then $xy\in \nm^{-1}(\mf{b})$, and because $\mf{b}$ is a prime ideal we must have either $x$ or $y$ is in $\nm^{-1}(\mf{b})$ which shows that $x$ or $y$ is in $(\nm^{-1}(\mf{b});\mf{b})$.

	Now suppose that $x\in \ghost(T)(C_p/e)$, $y\in \ghost(T)(C_p/C_p)$ and $\Q((\nm^{-1}(\mf{b});\mf{b}),x,y)$ is true.  Let us write $y = (y_1,y_2)$ where $y_1\in T(C_p/e)^{C_p}$ and $y_2\in \Phi^{C_p}(T)$.  If $x\in \nm^{-1}(\mf{b})$ then there is nothing to show, so assume that $x\notin\nm^{-1}(\mf{b})$.  We have
	\[
		\nm(x)\cdot y = \bigg(\prod\limits_{g\in G} gx,\nm(x) \bigg)\cdot (y_1,y_2) = \bigg( y_1\cdot \prod\limits_{g\in G}gx,y_2\nm(x) \bigg)
	\]
	is an element $(\nm^{-1}(\mf{b});\mf{b})$.  This implies that $\nm(x)y_2\in \mf{b}$, and by assumption we have $\nm(x)\notin \mf{b}$ so $y_2\in \mf{b}$.  On the other hand, we have
	\[
		\res(y_1,y_2)x = y_1x
	\]
	is in an element $(\nm^{-1}(\mf{b});\mf{b})$.  This implies that $y_1x\in \nm^{-1}(\mf{b})$ which is a prime ideal.  By assumption we have $x\notin \nm^{-1}(\mf{b})$ so $y_1\in \nm^{-1}(\mf{b})$.  Thus, either $x\in (\nm^{-1}(\mf{b});\mf{b})$ or $y\in (\nm^{-1}(\mf{b});\mf{b})$.

	Finally, suppose $x,y \in \ghost(T)(C_p/C_p)$.  Let use write $x = (x_1,x_2)$ and $y = (y_1,y_2)$ and suppose that $\Q((\nm^{-1}(\mf{b});\mf{b}),x,y)$ is true.  Since $xy\in (\nm^{-1}(\mf{b});\mf{b})$ we see that $x_2y_2\in \mf{b}$, and $x_1y_1\in \nm^{-1}(\mf{b})^{C_p}$.  Since $\mf{b}$ is prime, it must be that either $x_1$ or $y_1$ is in $\nm^{-1}(\mf{b})^{C_p}$ and either $x_2$ or $y_2$ is in $\mf{b}$.  Without loss of generality, assume that $x_1\in \nm^{-1}(\mf{b})^{C_p}$.  If $x_2\in \mf{b}$ we are done, as $(x_1,x_2)\in (\nm^{-1}(\mf{b});\mf{b})$.  If $x_2\notin \mf{b}$, then
	\[
		(x_1,x_2)\cdot \nm\res(y_1,y_2) = (x_1,x_2)\cdot(y_1^2,\nm(y_1)) = (x_1y_1^2,x_2\nm(y_1))
	\]
	is an element in $(\nm^{-1}(\mf{b});\mf{b})$.  Since $\mf{b}$ is prime, and $x_2\not\in \mf{b}$, we must have $\nm(y_1)\in \mf{b}$, so $y_1\in \nm^{-1}(\mf{b})^{C_p}$.  Since $x_2\notin \mf{b}$ implies $y_2\in \mf{b}$ we see that $y\in (\nm^{-1}(\mf{b});\mf{b})$.
\end{proof}

This gives us an explicit computation of all prime ideals of a $C_p$-Tambara functor.

\begin{proposition}\label{prop: primes of ghost}
	Let $T$ be a $C_p$-Tambara functor. Every prime ideal of $\ghost(T)$ is of one of the following two forms:
	\begin{enumerate}[(a)]
		\item $(\mf{a};\Phi^{C_p} T)$ for some $C_p$-prime ideal $\mf{a}\subset T(C_p/e)$, or
		\item $(\nm^{-1}(\mf{b});\mf{b})$ for some prime ideal $\mf{b}\in \Spec(\Phi^{C_p} T)$.
	\end{enumerate}
	In particular, there is a bijection
	\[
		\Spec(\ghost(T))\leftrightarrow \Spec(T(C_p/e)^{C_p})\amalg \Spec(\Phi^{C_p} T).
	\]
\end{proposition}
\begin{proof}
	Pairs of ideals as in the statement are indeed prime ideals of $\ghost(T)$ by \cref{lemma: when ideal pairs are prime in ghost}, so it suffices to show that all prime ideals have this form.  Let $\mf{p}\subseteq \ghost(T)$ be an arbitrary prime ideal.  Then $\mf{a} = \mf{p}(C_p/e)$ is a $C_p$-prime ideal by \cref{lem: bottom is G-prime}.  We know that $\mf{p}(C_p/C_p) = \mf{c}\times \mf{b}$ for some ideals $\mf{c}\subseteq T(C_p/e)^{C_p}$ and $\mf{b}\subseteq  \Phi^{C_p} T$ so we are done as soon as we show that $\mf{c} = \mf{a}^{C_p}$ and $\mf{b}$ is either a prime ideal or the entire ring.  The equation $\nm(\mf{a})+\tau\subseteq \mf{b}$ is automatically satisfied because $\mf{p}$ is an ideal of $\ghost(T)$.

	To see that $\mf{c} = \mf{a}^{C_p}$ we first note that $\mf{c}\ = \pi_1(\mf{c}\times \mf{b}) = \res(\mf{c} \times \mf{b}) \subseteq \mf{a}^{C_p}$.  On the other hand, if $x\in \mf{a}^{C_p}$ then $(x^p,\nm(x)+\tau)  = \nm(x)$ is an element in $\mf{c} \times \mf{b}$ since $\mf{p}$ is an ideal.  By \cref{thm: Tambara prime are radical}, $\mf{c}\times \mf{b}$ is a radical ideal, hence $\mf{c}$ is a radical ideal and so $x\in \mf{c}$ which proves $\mf{c} = \mf{a}^{C_p}$.

	We now show that $\mf{b}$ is a prime ideal or not proper.  Suppose that $x,y\in \Phi^{C_p} T$ and $xy\in \mf{b}$.  Then $(0,x)\cdot(0,y) = (0,xy)\in \mf{p}$.  Since all non-trivial multiplicative translates of $(0,x)$ and $(0,y)$ are zero we see that $\Q(\mf{p},(0,x),(0,y))$ is satisfied so $x$ or $y$ is in $\mf{b}$ by the primeness of $\mf{p}$.  Thus $\mf{b}$ is prime so long as it is a proper ideal.

	The claimed bijection follows from \cref{theorem: all the equispecs are the same}.
\end{proof}

\begin{remark}
	The bijection in the proposition is \emph{not} a homeomorphism in general.  Indeed, it is usually not even a poset isomorphism, as we have inclusions $(\nm^{-1}(\mf{b});\mf{b})\subseteq (\mf{a};\Phi^{C_p} T)$ whenever $\nm^{-1}(\mf{b})\subseteq \mf{a}$.
\end{remark}

\begin{example}\label{example: Burnside ghost picture}
	Using the previous results, we may recover Nakaoka's calculation of the spectrum of the Burnside Tambara functor over $C_p$ \cite{Nak2014a}; see also \cite{CG2023,6A25}. \Cref{example: ghost not inj on spec} described the Lewis diagram of $\ghost(\uA)$. The primes of $\ghost(\uA)$ are then easy to classify; the poset $\Spec(\ghost(\uA))$ is two copies of $\Spec \Z$, with one living above the other (\Cref{figure: Spec ghost}).
	\begin{figure}[h]
		\begin{tikzpicture}[scale=0.75]
			\draw (0,0) circle (0.075);
			\draw[fill] (1.5,0.5) circle (0.075);
			\draw[fill] (2,0.7) circle (0.075);
			\draw[fill] (2.5,0.85) circle (0.075);
			\draw (0,3) circle (0.075) node[above=0.1]{$0$};
			\draw[fill] (1.5,2.5) circle (0.075) node[above=0.1]{$2$};
			\draw[fill] (2,2.3) circle (0.075) node[above=0.1]{$3$};
			\draw[fill] (2.5,2.15) circle (0.075) node[above=0.1]{$5$};
			\draw[fill] (5,1.5) circle (0.075) node[above=0.1]{$p$};
			\draw plot [smooth] coordinates {(0.25,3) (5,1.5) (9,3)};
			\draw plot [smooth] coordinates {(0.25,0) (5,1.5) (9,0)};
			\draw[->] (0,0.25) -- (0,2.75);
			\draw[->] (1.5,0.75) -- (1.5,2.25);
			\draw[->] (2,0.95) -- (2,2.05);
			\draw[->] (2.5,1.1) -- (2.5,1.9);
			\node at (5,-1) {$\Spec(\uA)$};
		\end{tikzpicture}
		\hspace{2cm}
		\begin{tikzpicture}[scale=0.75]
			\draw (0,0) circle (0.075);
			\draw (0.25, 0) -- (9, 0);
			\draw[fill] (1.5,0) circle (0.075);
			\draw[fill] (2,0) circle (0.075);
			\draw[fill] (2.5,0) circle (0.075);
			\draw[fill] (5,0) circle (0.075);
			\draw (0,3) circle (0.075) node[above=0.1]{$0$};
			\draw (0.25, 3) -- (9, 3);
			\draw[fill] (1.5,3) circle (0.075) node[above=0.1]{$2$};
			\draw[fill] (2,3) circle (0.075) node[above=0.1]{$3$};
			\draw[fill] (2.5,3) circle (0.075) node[above=0.1]{$5$};
			\draw[fill] (5,3) circle (0.075) node[above=0.1]{$p$};
			\draw[->] (0,0.25) -- (0,2.75);
			\draw[->] (1.5,0.25) -- (1.5,2.75);
			\draw[->] (2,0.25) -- (2,2.75);
			\draw[->] (2.5,0.25) -- (2.5,2.75);
			\draw[->] (5,0.25) -- (5,2.75);
			\node at (5,-1) {$\Spec(\Gamma(\uA))$};
		\end{tikzpicture}
		\caption{$\Spec(\uA)$ and $\Spec(\protect\ghost(\uA))$}
		\label{figure: Spec ghost}
	\end{figure}

	Now applying $\ghostmap^*$ yields $\Spec(\uA)$, and here $\ghostmap^*$ is almost injective; it only collapses the two primes lying over $p$. Thus, we get the poset $\Spec(\uA) = \ghostmap^* \Spec(\ghost(\uA))$ (\Cref{figure: Spec ghost}), which looks like two copies of $\Spec(\Z)$ glued together at the point $(p)$ (with a different poset structure). We note here that $\uA \to \ghost(\uA)$ is an example of a levelwise-integral extension of $C_p$-Tambara functors which does not satisfy the incomparability property from commutative algebra.
\end{example}

We have from \Cref{theorem: all the equispecs are the same} that $\mf{a} \mapsto \mf{a}^{C_p}$ defines an isomorphism from the poset of $C_p$-prime ideals of $T(C_p/e)$ to the poset of prime ideals of $T(C_p/e)^{C_p}$ for any $C_p$-Tambara functor $T$. Thus, the primes of $\ghost(T)$ are determined by their $C_p/C_p$ levels, and we can equivalently describe the spectrum of $\ghost(T)$ without reference to $C_p$-primes.

\begin{corollary}\label{prop:alt classification of primes of ghost}
	Let $T$ be a $C_p$-Tambara functor. There is an isomorphism of posets between $\Spec(\ghost(T))$ and the set of pairs $(\mf{a}; \mf{b})$ where either:
	\begin{itemize}
		\item $\mf{a}$ is a prime ideal of $T(C_p/e)^{C_p}$ and $\mf{b} = \Phi^{C_p} T$;
		\item $\mf{b}$ is a prime ideal of $\Phi^{C_p} T$ and $\mf{a} = (\nm^{-1} \mf{b})^{C_p}$.
	\end{itemize}
	The ordering on such pairs is given by $(\mf{a}, \mf{b}) \leq (\mf{a}', \mf{b}')$ if and only if $\mf{a} \subseteq \mf{a}'$ and $\mf{b} \subseteq \mf{b}'$.
\end{corollary}

The classification of ideals in the ghost Tambara functor provides an efficient means of checking when a given Tambara functor is a domain, which is typically very difficult to check.

\begin{proposition}\label{prop: domain or not domain}
	Let $T$ be a nonzero $C_p$-Tambara functor.
	\begin{enumerate}[(a)]
		\item If $\nm\colon T(C_p/e) \to \Phi^{C_p} T$ is injective, $\Phi^{C_p} T$ is a domain, and $T(C_p/C_p)$ has no $p$-torsion, then $T$ is a domain.
		\item If $\nm\colon T(C_p/e) \to \Phi^{C_p} T$ is not injective and $\res\colon T(C_p/C_p) \to T(C_p/e)$ is not injective, then $T$ is not a domain.
	\end{enumerate}
\end{proposition}

\begin{proof}
	To prove $(a)$, we need to show that $0 \subseteq T$ is a prime ideal. By \Cref{cor: ghost injective}, our assumption that $T(C_p/C_p)$ has no $p$-torsion implies that $\ghostmap_T$ is injective, so it suffices to show that the zero ideal is prime in $\ghost(T)$. Since $\Phi^{C_p} T$ is a domain, we have that $(\nm^{-1}(0); 0)$ is a prime ideal of $\ghost(T)$, and since $\nm$ is injective, this ideal is $(0; 0) = 0$.

	Since $\ghostmap^*\colon \Spec(\ghost(T)) \to \Spec(T)$ is surjective, to prove $(b)$ it suffices to show that $0 \subset T$ is not in the image of $\ghostmap^*$. Assume, for the sake of contradiction, that $\ghostmap^*((\mf{p}; \mf{q})) = 0$. Any prime ideal $(\mf{p}; \mf{q})$ of $\ghost(T)$ has the form $(\mf{a}; \Phi^{C_p} T)$ or $(\nm^{-1}(\mf{b});\mf{b})$ where $\mf{a} \subseteq T(C_p/e)^{C_p}$ and $\mf{b} \subseteq \Phi^{C_p} T$ are prime ideals. Since $\ghostmap(C_p/e)$ is the identity, we must have $\mf{p}=0$, and since $\nm$ is not injective, we must have $(\mf{p};\mf{q}) = (0;\Phi^{C_p} T)$. We then have $s \in \ghostmap^{-1}(C_p/C_p)(\mf{p};\mf{q})$ if and only if $s \in \ker(\res)$ and $s \in T(C_p/C_p)$; the second condition is always satisfied, and since we assumed $\res$ is not injective, there exists some $s \neq 0 \in \ghostmap^{-1}(C_p/C_p)(\mf{p};\mf{q})$, in contradiction to the assumption that $\chi^*((\mf{p};\mf{q})) = 0$. Therefore $0 \subset T$ is not in the image of $\ghostmap^*$, so it is not prime.
\end{proof}

\begin{corollary}
	The functors $\uA$ and $\uA[C_p/C_p]$ are $C_p$-Tambara domains, and $\uA[C_p/e]$ is not a $C_p$-Tambara domain.
\end{corollary}
It is already known (\cite[Theorem 4.40]{Nak2012}) that $\uA$ is a $C_p$-Tambara domain; we provide another proof here as an immediate consequence of the prior proposition.
\begin{proof}
	For $\uA$, we have that $\nm \colon \uA(C_p/e) \to \Phi^{C_p} \uA$ is an isomorphism, and we note that $\uA(C_p/C_p) \cong \Z[t]/\langle t^2-pt \rangle$ has no $p$-torsion. Moreover, $\Phi^{C_p} \uA \cong \Z$ is a domain. By \Cref{prop: domain or not domain} we conclude that $\uA$ is a domain.

	For $\uA[C_p/C_p]$, we have that $\nm \colon \uA[C_p/C_p](C_p/e) \to \Phi^{C_p} \uA[C_p/C_p]$ is isomorphic to $x \mapsto n \colon \Z[x] \to \Z[x,n]$, which is injective. Moreover, \[\uA[C_p/C_p](C_p/C_p) \cong \Z[t,x,n]/\langle t^2-pt, tn-tx^2 \rangle\] has no $p$-torsion, and $\Phi^{C_p} \uA[C_p/C_p] = \Z[x,n]$ is a domain. By \Cref{prop: domain or not domain} we conclude that $\uA[C_p/C_p]$ is a domain.

	For $\uA[C_p/e]$, we see that $\nm \colon \uA[C_p/e](C_p/e) \to \Phi^{C_p} \uA[C_p/e]$ is not injective (e.g. $x_e - x_\gamma \mapsto 0$) and $\res \colon \uA[C_p/e](C_p/C_p) \to \uA[C_p/e](C_p/e)$ is not injective (e.g. $t-p \mapsto 0$). By \Cref{prop: domain or not domain} we conclude that $\uA[C_p/e]$ is not a domain.
\end{proof}

\section{Krull Dimension}
\label{sec: Krull}

In this section we use the ghost construction to compute the Krull dimension of several Tambara functors which appear throughout the paper.  Nonequivariantly, if $R\to S$ is an integral extension of commutative rings then $R$ and $S$ necessarily have the same Krull dimension.  One might hope that the same is true for Tambara functors, however examples show that this is false in general.  In particular, even though the map $T\to \ghost(T)$ is often a levelwise integral extension, it is possible for the dimensions to differ.  The main result of this section is that we can use the ghost map to provide an upper bound on the Krull dimension of a Tambara functor $T$.  In most examples of interest this upper bound is an equality.

\begin{definition}
	\label{def: krull dim}
	The \emph{Krull dimension} $\dim(T)$ of a nonzero Tambara functor $T$ is the supremum over all integers $k$ such that there is a chain
	\[
		\fp_0 \subset \fp_1\subset \fp_2\subset\dots\subset \fp_k
	\]
	of prime ideals in $T$. We declare by convention that $\dim(0) = -\infty$.
\end{definition}

We next record some consequences of lying over, and the going up theorem (\cref{theorem: going up}). If $T \to S$ is a morphism of Tambara functors which satisfies both going up and lying over, then a quick induction implies that any chain
\[
	\mf{p}_0\subset \mf{p}_1\subset\dots\subset \mf{p}_n
\]
of prime ideals in $T$ can be lifted to a chain
\[
	\mf{q}_0\subset \mf{q}_1\subset\dots\subset \mf{q}_n
\]
of prime ideals in $S$. We know that integral extensions satisfy both going up and lying over, and so we obtain:

\begin{proposition}\label{prop: going up and krull dim}
	Let $T$ and $S$ be Tambara functors. If $T \to S$ is a levelwise integral extension, then $\dim(T) \leq \dim(S)$.
\end{proposition}

In commutative algebra, integral extensions also satisfy the incomparability property, which forces the inequality above to be an equality \cite[Theorem 48]{Kaplansky}. Integral extensions of Tambara functors do not necessarily satisfy the incomparability property (\Cref{example: Burnside ghost picture}), and indeed there exist integral extensions of Tambara functors where this inequality is strict (\Cref{example: strict inequality of dimensions} below).

For the remainder of this section we restrict to the case $G = C_p$ and study the relationship between the dimension of $T$ and the dimension of $\ghost(T)$. In this setting, we saw (\Cref{cor: ghost has going up and lying over}) that $\ghostmap_T \colon T \to \ghost(T)$ satisfies both going up and lying over.

\begin{corollary}\label{corollary: Krull dim bound with ghost}
	For any Tambara functor $T$ we have
	\[
		\dim(T)\leq \dim(\ghost(T)).
	\]
\end{corollary}

The following example shows this bound can be strict.

\begin{example}\label{example: strict inequality of dimensions}
	Let $T = \FP(\F_p)$ be the fixed point $C_p$-Tambara functor of $\F_p$ with trivial $C_p$-action. The Tambara functor $T$ is a Tambara field \cite[Definition 4.28(2)]{Nak2012}, therefore $\Spec(T)$ is a singleton set. In particular, $\dim(T)=0$. On the other hand, \cref{prop: primes of ghost} shows that $\Spec(\ghost(T))$ has two points: $(0;0) \subset (0;\mathbb{F}_p)$ and thus $\dim(\Spec(\ghost(T)))=1$.
\end{example}

Our next goal is to find a way to compute the dimension of $\ghost(T)$. To this end, recall from \cref{prop: primes of ghost} that every prime ideal of $\ghost(T)$ can be written as either $(\mf{a};\Phi^{C_p}(T))$ or $(\nm^{-1}(\mf{b});\mf{b})$ where $\mf{a}\subseteq T(C_p/e)$ is a $C_p$-prime, $\mf{b}\subseteq \Phi^{C_p}(T)$ is prime, and $\nm\colon T(C_p/e)\to \Phi^{C_p}(T)$ is the ring map induced by applying the norm and modding out by transfers.  There are inclusions
\[
	(\nm^{-1}(\mf{b});\mf{b}) \subset (\mf{a};\Phi^{C_p}(T))
\]
whenever $\nm^{-1}(\mf{b})\subseteq \mf{a}$.

\begin{lemma}
	Let $T$ be a $C_p$-Tambara functor, and let $J$ be any $C_p$-invariant ideal of $T(C_p/e)$ such that $\ker(\nm) \leq J$.  Then $J$ is $C_p$-prime if and only if it is prime.
\end{lemma}
\begin{proof}
	Since every prime ideal is $C_p$-prime it suffices to prove the reverse implication. Suppose that $J$ is $C_p$-prime. Equivalently, $0$ is $C_p$-prime in $T(C_p/e)/J$. Since $\nm(x) = \nm(g \cdot x)$ for all $g \in C_p$ and all $x \in T(C_p/e)$, the $C_p$ action on $T(C_p/e)/J$ is trivial. Thus, $0$ is prime in $T(C_p/e)/J$, i.e. $J$ is prime.
\end{proof}

\begin{corollary}\label{corollary: a is often actually prime}
	Let $T$ be a $C_p$-Tambara functor. If $(\nm^{-1}(\mf{b});\mf{b}) \subseteq (\mf{a};\Phi^{C_p}(T))$ then $\mf{a}$ is a prime ideal of $T(C_p/e)$.
\end{corollary}
\begin{proof}
	This is immediate from the lemma because $\ker(\nm)\subseteq \nm^{-1}(\mf{b})$.
\end{proof}

Before going further we recall the notion of \emph{height} and \emph{coheight} of an ideal in a commutative ring.

\begin{definition}
	Let $R$ be a commutative ring and let $J \trianglelefteq R$ be a prime ideal. The \emph{height} of $J$, denoted by $\height_R(J)$, is the supremum over integers $k$ such that there is a chain of ideals
	\[
		J\supset \mf{q}_1 \supset \dots \supset \mf{q}_{k}
	\]
	where each $\mf{q}_i$ is prime. 
\end{definition}
\begin{definition}
	Let $R$ be a commutative ring and let $J \trianglelefteq R$ be a prime ideal. The \emph{coheight} of $J$ in $R$, denoted $\coht_R(J)$, is the supremum over integers $k$ such that there is a chain of ideals
	\[
		J\subset \mf{q}_1\subset\dots\subset \mf{q}_{k}\subset R 
	\]
	where each $\mf{q}_i$ is prime.
\end{definition}
An immediate consequence of the definitions is that for any ring $R$ and any prime ideal $J$ we have
\[
	\height_R(J)+\coht_R(J)\leq \dim(R).
\]
We also note that $\coht_R(J) = \dim(R/J)$ and $\height_R(J) = \dim(R_J)$, where $R_J$ is the localization of $R$ at $J$.

Later, we will need the following lemma.
\begin{lemma}[{\cite[Theorems 44 and 47]{Kaplansky}}]\label{lemma: integral coheight equality}
	Let $R$ and $S$ be commutative rings. If $f\colon R\to S$ is an integral ring map then for any prime ideal $\mf{p}$ of $S$ we have $\coht_S(\mf{p}) = \coht_{R}(f^*(\mf{p}))$.
\end{lemma}
\begin{proof}
	Here $R/f^* \mf{p} \to S/\mf{p}$ is an integral extension of rings, and thus \[\coht_R(f^* \mf{p}) = \dim(R/f^* \mf{p}) = \dim(S/\mf{p}) = \coht_S(\mf{p}).\qedhere\]
\end{proof}

\medskip

We now turn to finding an upper bound on the Krull dimension of $\ghost(T)$ for some $C_p$-Tambara functor $T$.

\begin{proposition}\label{proposition: the height coheight bound for Krull dimension}
	For a $C_p$-Tambara functor $T$ let $M(T)$ denote the supremum of the integers
	\[
		1 + \height_{\Phi^{C_p}(T)}(\mf{b})+ \coht_{T(C_p/e)}(\nm^{-1}(\mf{b}))
	\]
	where $\mf{b}$ runs over all prime ideals of $\Phi^{C_p}(T)$.\footnote{If $\Phi^{C_p}(T) = 0$, we have by convention that $M(T) = -\infty$.} Then
	\[
		\dim(\ghost(T))\leq \mathrm{max}(\dim T(C_p/e),M(T)).
	\]
\end{proposition}
\begin{proof}
	Suppose that
	\[
		\mf{p}_0\subset \mf{p}_1\subset \dots\subset \mf{p}_k
	\]
	is a chain of prime ideals in $\ghost(T)$. There are three cases, depending on the structure of the ideals $\mf{p}_i \subseteq \Gamma(T)$.

	For the first case, if all $\mf{p}_i$ are of the form $(\mf{a};\Phi^{C_p}(T))$, then (by \cref{theorem: all the equispecs are the same}) $k$ is no larger than $\dim(T(C_p/e)^{C_p})$. Moreover, since the inclusion of fixed points is an integral extension we have
	\[
		\dim(T(C_p/e)^{C_p})=\dim(T(C_p/e)),
	\]
	so $k\leq \dim(T(C_p/e))$.

	For the second case, if all the ideals are of the form $(\nm^{-1}(\mf{b});\mf{b})$, then we can extend the chain: if $\mf{p}_k = (\nm^{-1}(\mf{b}_k);\mf{b}_k)$ then we have $\mf{p}_k\subseteq (\nm^{-1}(\mf{b}_k);\Phi^{C_p}(T))$.  In particular, such a chain is never maximal.

	For the final case, suppose we have a $0\leq r<k$ such that $\mf{p}_i = (\nm^{-1}(\mf{b}_i);\mf{b}_i)$ for all $0\leq i\leq r$ and $\mf{p}_i = (\mf{a}_i;\Phi^{C_p}(T))$ for all $i>r$.  Since we are interested in understanding chains of maximum possible length, and
	\[
		(\nm^{-1}(\mf{b}_r);\mf{b}_r)\subseteq (\nm^{-1}(\mf{b}_r);\Phi^{C_p}(T))\subseteq (\mf{a}_{r+1}; \Phi^{C_p}(T)),
	\]
	we may safely assume $\mf{a}_{r+1} = \nm^{-1}(\mf{b}_r)$ and thus our chain has the form
	\[
		(\nm^{-1}(\mf{b}_0);\mf{b}_0)\subset\dots\subset (\nm^{-1}(\mf{b}_r);\mf{b}_r)\subset (\nm^{-1}(\mf{b}_r);\Phi^{C_p}(T))\subset \dots\subset (\mf{a}_{k};\Phi^{C_p}(T)).
	\]
	Note that $r\leq \height_{\Phi^{C_p}(T)}(\mf{b}_r)$, and that
	\[
		(\nm^{-1}(\mf{b}_r);\Phi^{C_p}(T))\subset \dots\subset (\mf{a}_{k};\Phi^{C_p}(T))
	\]
	is a length-$(k-r-1)$ chain of $C_p$-prime ideals of $T(C_p/e)$ containing $\mf{a}_{r+1} = \nm^{-1}(\mf{b}_r)$. By \cref{corollary: a is often actually prime}, all $C_p$-prime ideals in such a chain are actually prime, so we have
	\begin{equation}\label{equation: coheight inequality.}
		k\leq 1+r+ \coht_{T(C_p/e)}(\mf{a}_{r+1})
		\leq
		1+  \height_{\Phi^{C_p}(T)}(\mf{b}_r)+ \coht_{T(C_p/e)}(\nm^{-1}(\mf{b}_r))
		\leq M(T).
	\end{equation}
	\noindent Combining this with the bound $k\leq \dim(T(C_p/e))$ from above, we obtain the claimed bound on the dimension of $\ghost(T)$.
%
\end{proof}

We can improve this bound with some additional assumptions on the norm $\nm\colon T(C_p/e)\to \Phi^{C_p}(T)$.

\begin{thm}\label{theorem: krull dimension bounds for integral norm map}
	Let $T$ be a $C_p$-Tambara functor such that $\nm\colon T(C_p/e)\to \Phi^{C_p}(T)$ is an integral map.  Then
	\[
		\dim(T)\leq\dim(\ghost(T))\leq  \mathrm{max}(\dim(T(C_p/e)),\dim(\Phi^{C_p}(T))+1).
	\]
	If, in addition, the map $\nm$ is injective then $\dim(T(C_p/e)) = \dim(\Phi^{C_p}(T))$ and we obtain $\dim(T)\leq\dim(T(C_p/e))+1$.
\end{thm}
\begin{proof}
	The first inequality is immediate from \cref{corollary: Krull dim bound with ghost} so we only need to establish the bound on the dimension of $\ghost(T)$. For any $\mf{b}\in \Spec(\Phi^{C_p}(T))$,  \cref{lemma: integral coheight equality} tells us that
	\[
		\coht_{T(C_p/e)}(\nm^{-1} \mf{b}) = \coht_{\Phi^{C_p}(T)}(\mf{b}),
	\]
	and thus for any $\mf{b}\in \Spec(\Phi^{C_p}(T))$ we have
	\begin{align*}
		1+  \height_{\Phi^{C_p}(T)}(\mf{b})+ \coht_{T(C_p/e)}(\nm^{-1}(\mf{b}))
		 & =  1+  \height_{\Phi^{C_p}(T)}(\mf{b})+ \coht_{\Phi^{C_p}(T)}(\mf{b}) \\
		 & \leq 1+\dim(\Phi^{C_p}(T))
	\end{align*}
	and  thus $M(T)\leq 1+\dim(\Phi^{C_p}(T))$.  This gives us the inequality from the statement.  The final claim follows from the fact that if $R\to S$ is an integral extension of commutative rings then $\dim(R) = \dim(S)$.
\end{proof}

\begin{corollary}\label{corollary: Krull dimension computer}
	Let $T$ be a $C_p$-Tambara functor. If $\nm\colon T(C_p/e)\to \Phi^{C_p}(T)$ is integral and $\dim(T(C_p/e))\geq \dim(\Phi^{C_p}(T))+1$ then $\dim(T) = \dim(\ghost(T)) = \dim(T(C_p/e))$.
\end{corollary}
\begin{proof}
	Let $D = \dim(T(C_p/e))$. It is clear from the theorem that we have $\dim(T) \leq \dim(\ghost(T)) \leq D$, so it suffices to show that $D \leq \dim(T)$. Pick a chain
	\[
		\mf{a}_0\subset \mf{a}_{1}\subset\dots\subset \mf{a}_{D}
	\]
	of $C_p$-prime ideals in $T(C_p/e)$.  This gives a chain
	\[
		(\mf{a}_0;\Phi^{C_p}(T))\subset (\mf{a}_{1};\Phi^{C_p}(T))\subset\dots\subset (\mf{a}_{D};\Phi^{C_p}(T))
	\]
	in $\ghost(T)$.  The preimages of these ideals under the map $\ghostmap_T\colon T\to \ghost(T)$ are distinct because the bottom level of each is $\mf{a}_i$.  Thus this is a chain of prime ideals of $T$ with length $D$.
\end{proof}

We now turn to computing the Krull dimension in several examples.  In the first three of these examples the map $\nm\colon T(C_p/e)\to \Phi^{C_p}(T)$ is integral so we apply \cref{theorem: krull dimension bounds for integral norm map}.  In the last example we calculate using the bound of \cref{proposition: the height coheight bound for Krull dimension}.

\begin{example}\label{example: Krull dim of polynomial free generator}
	We consider the free Tambara functor $T = \uA[C_p/e]$.  In this case we have $\Phi^{C_p}(T) = \Z[n]$, and $T(C_p/e) = \Z[x_0,\dots,x_{p-1}]$.  The norm map is determined by $\nm(x_i)=n$ for all $i$.  This map is surjective and hence integral.  Since $\dim(\Z[n]) +1= 3\leq p+1 = \dim(\mathbb{Z}[x_0,\dots,x_{p-1}])$, \cref{corollary: Krull dimension computer} tells us that $\dim(\uA[C_p/e]) = p+1$ for all $p$.
\end{example}

\begin{example}
	We consider the $C_p$-Burnside Tambara functor $\uA$. We have $\Phi^{C_p}(\uA) \cong \uA(C_p/e)\cong \Z$.  The norm map is an isomorphism.  Thus \cref{theorem: krull dimension bounds for integral norm map} tells us that
	\[
		\dim(\uA) \leq \dim(\Z)+1 =2.
	\]
	In fact we have $\dim(\uA)=2$.  A chain in $\ghost(\uA)$ of length $2$, which restricts to a chain of length $2$ in $\uA$ is given by
	\[
		(0;0)\subset (q;q)\subset (q;\mathbb{Z})
	\]
	where $q$ is any prime except $p$.
\end{example}
\begin{example}\label{example: Krull dim of RU}
	Let $\RU$ be the complex representation $C_p$-Tambara functor given by the Lewis diagram
	\[
		\begin{tikzcd}[row sep=huge]
			\Z[C_p]
			\ar[d, "\res" description]
			\\
			\Z
			\ar[u, bend right=50, "\nm"{right}]
			\ar[u, bend left=50, "\tr"{left}]
			\arrow[from=2-1, to=2-1, loop, in=300, out=240, distance=5mm, "{\mathrm{trivial}}"']
		\end{tikzcd}
	\]
	where the transfer, norm, and restriction are determined by
	\[
		\tr(k) = k(1+\gamma+\dots+\gamma^{p-1}),\quad \nm(k) = k+\left(\frac{k^p-k}{p}\right)\cdot(1+\gamma+\dots+\gamma^{p-1}), \quad \text{and} \quad \res(\gamma)=1.
	\]

	At the level $C_p/C_p$, the element $1$ represents the trivial complex representation of $C_p$, and the element $\gamma$ represents the $1$-dimensional $C_p$-representation given by multiplication by $e^{2\pi i/p}$.

	We have $\Phi^{C_p}(RU) = \Z[\xi]$ where $\xi$ is a $p$-th root of unity and $RU(C_p/e)=\mathbb{Z}$.  The map $\nm\colon \Z\to \Z[\xi]$ is the unique ring map and is an integral extension because $\Z[\xi]$ is the ring of integers in the field extension $\QQ\to \QQ[\xi]$.  Since the norm map is also injective we have
	\[
		\dim(\RU)\leq \dim(\Z)+1=2.
	\]

	In fact we have $\dim(\RU)=2$: let $q$ be a prime which is not $p$ and let $\mf{b}\subset \Z[\xi]$ be a prime ideal which lies over $q$, i.e.\ $\nm^{-1}(\mf{b}) = \mf{b}\cap \Z = \langle q \rangle$.  The chain of ideals
	\[
		(0;0)\subset (\langle q \rangle ;\mf{b})\subset (\langle q \rangle;\mathbb{Z}[\xi])
	\]
	is a proper chain of length $2$ in $\Spec(\ghost(\RU))$.  We need to show that the preimage of this chain is a proper chain in $\RU$.  It is clear that $(0;0)\subset (\langle q \rangle ;\mf{b})$ remains a proper inclusion.  To see that the other inclusion remains proper it suffices to find an element $z\in \RU(C_p/C_p)$ so that $\res(z)\in \langle q \rangle$ but $z+ \tau\notin \mf{b}$.

	Pick a prime $r$ equal to neither $p$ nor $q$ and pick integers $a$ and $b$ so that $ap+bq=r$. Let $z = r-a\tr(1)\in \RU(C_p/C_p)$.  The restriction of $z$ is
	\[
		\res(z) = r-ap =bq
	\]
	which is in $\langle q \rangle$, and $z+\tau=r+\tau$. Since $r \neq q$, we have $r \notin \langle q \rangle = \nm^{-1}(\mf{b})$, so $r+\tau \notin \mf{b}$, so $z + \tau \notin \mf{b}$.
\end{example}

The final example we consider is $T=\uA[C_p/C_p]$.  In this case, $\Phi^{C_p}(T) \cong \mathbb{Z}[n,x]$, $T(C_p/e)\cong \Z[n]$, and the norm map $\nm\colon \Z[n]\to \Z[x,n]$ is the inclusion.  Since the norm map is not an integral map we need to treat this using \cref{proposition: the height coheight bound for Krull dimension}.

The next lemma follows from \cite[Theorem 4.19]{Eisenbud} and the fact that $\Z$ is an Jacobson ring.

\begin{lemma}\label{lem: finite type over Jacobson => pres closed points}
	If $\mf{m}$ is a maximal ideal in $\Z[x,n]$ then $\mf{m}\cap \Z[n]$ is a maximal ideal of $\Z[n]$.
\end{lemma}


\begin{example}\label{example: Krull dim polynomial fixed generator}
	Let $T = \uA[C_p/C_p]$. We will show that $\dim(T)=4$.  We begin by establishing that $5$ is an upper bound using \cref{proposition: the height coheight bound for Krull dimension}.  We have $\dim(T(C_p/e))=2$, so we will have to evaluate $M(T)$.  Suppose that $\mf{b}\in \Spec(\Z[x,n])$ is a prime ideal. We will show that
	\[
		1+\height_{\Z[x,n]}(\mf{b})+\coht_{\Z[n]}(\nm^{-1}(\mf{b}))\leq 4
	\]
	for any choice of $\mf{b}$.  We divide into different cases depending on the height of $\mf{b}$.

	If $\height_{\Z[x,n]}(\mf{b})\leq 1$ we have
	\[
		1+\height_{\Z[x,n]}(\mf{b})+\coht_{\Z[n]}(\nm^{-1}(\mf{b})) \leq  1+1+\dim(\Z[n]) = 4.
	\]

	If $\height_{\Z[x,n]}(\mf{b})=2$, then \cite[Theorem 37]{Kaplansky} implies that $\height_{\Z[n]}(\nm^{-1} \mf{b}) \geq 1$, so $\coht_{\Z[n]}(\nm^{-1}(\mf{b}))\leq \dim(\Z[n])-1 =1$. Then we have
	\[
		1+\height_{\Z[x,n]}(\mf{b})+\coht_{\Z[n]}(\nm^{-1}(\mf{b})) \leq  1+2+1 = 4.
	\]

	Finally, if $\height_{\Z[x,n]}(\mf{b})=3$ then $\mf{b}$ is a maximal ideal so \Cref{lem: finite type over Jacobson => pres closed points} implies that $\nm^{-1}(\mf{b})$ is also maximal (i.e. $\coht_{\Z[n]}(\nm^{-1}(\mf{b}))=0$), so we have
	\[
		1+\height_{\Z[x,n]}(\mf{b})+\coht_{\Z[n]}(\nm^{-1}(\mf{b})) \leq  1+3+0 = 4.
	\]

	Thus, we have shown that $M(T)\leq 4$. Finally, note that for any prime $q\neq p$ there is a proper chain of prime ideals in $\ghost(T)$ given by
	\[
		(0;0)\subset (\langle q \rangle;\langle q \rangle)\subset (\langle q,x \rangle; \langle q, x\rangle )\subset (\langle q,x \rangle; \langle q,x,n \rangle)\subset (\langle q,x \rangle;\Z[x,n])
	\]
	which has length $4$. Call the preimages of these ideals under $\ghostmap_{T}$ $\mf{p}_0,\dots,\mf{p}_4$. In $T(C_p/C_p)$, consider the elements of   $a_1=q$, $a_2=x$, $a_3=n$, and $a_4=q-p+t$.  One checks that $a_i\in \mf{p}_i\setminus \mf{p}_{i-1}$ for all $i$ hence the $\mf{p}_i$ give a proper chain of prime ideals of length $4$ in $T$, whence $\dim(T)=4$.

\end{example}

\section{Examples}

With all of the necessary theory in place,  we now turn to computing examples of Nakaoka spectra.

\subsection{Complex Representation Ring}

Throughout this section, let $\gamma$ be the generator of $C_p$. Recall the complex representation ring Tambara functor $\RU$ from \cref{example: Krull dim of RU}.

To compute $\ghost(\RU)$ we observe that the ring
\[
	\Z[C_p]/\langle 1+\gamma+\dots+\gamma^{p-1} \rangle
\]
is isomorphic to $\Z[\xi]$ where $\xi$ is a primitive $p$-th root of unity.  Thus the ghost Tambara functor $\ghost(\RU)$ is given by the Lewis diagram
\[
	\begin{tikzcd}[row sep=huge]
		\Z\times \Z[\xi]
		\ar[d, "\res" description]
		\\
		\Z
		\ar[u, bend right=50, "\nm"{right}]
		\ar[u, bend left=50, "\tr"{left}]
		\arrow[from=2-1, to=2-1, loop, in=300, out=240, distance=5mm, "{\mathrm{trivial}}"']
	\end{tikzcd}
\]
where the transfer, norm, and restriction are determined by
\[
	\tr(k) = (kp,0),\quad \nm(k) = (k^p,k), \quad \text{and} \quad \res(k,f)=k.
\]

We now turn to the determination of the primes of $\ghost(\RU)$.  Recall that by \cref{prop: primes of ghost} the Tambara prime ideals are pairs of the form $(\mf{a};\Z[\xi])$ or $(\nm^{-1}(\mf{b});\mf{b})$ where $\mf{a}\subset \Z$ and $\mf{b}\subset \Z[\xi]$ are prime ideals. Any such $\mf{a}$ has the form $0$ or $\langle q \rangle$ where $q$ is a prime number. To describe $\Spec(\Z[\xi])$, we recall the following:

\begin{thm}[{\cite[Chapter 6]{milneANT}}]
	Let $q$ be a prime number.
	\begin{enumerate}
		\item If $q=p$ then there is exactly one prime $\mf{b}\subset \Z[\xi]$ such that $\mf{b} \cap \Z = \langle p \rangle$. In particular, $\mf{b} = \langle \xi -1 \rangle$ and $\mf{b}^{p-1} = \langle p \rangle$.
		\item If $q\neq p$ then let $f$ be the smallest positive integer so that $q^f\equiv 1 \mod p$ and let $e = \frac{p-1}{f}$.  Then there are exactly $e$ prime ideals $\mf{b}_1,\dots,\mf{b}_e$ so that $\mf{b}_i\cap \Z = \langle q \rangle$ and for each $i$ we have $\mf{b}^f = \langle q \rangle$.
	\end{enumerate}
\end{thm}

\begin{thm}\label{thm:RU}
	Let $\RU$ denote the complex representation ring $C_p$-Tambara functor. There is a bijection
	$$\Spec(\RU) \cong \dfrac{\Spec(\Z) \times \{s, \ell\}}{(\langle p\rangle ,s) \sim (\langle p\rangle,\ell)}.$$
	Under this bijection, $\Spec(\RU)$ has a subbasis of closed sets
	\[\{V_e(n) \colon n \in \Z\} \cup \{V_{C_p}(f) \colon f \in \Z[x]/(x^p-1)\},\]
	where $V_e(n)$ has points determined by
	\[
		V_e(n) =
		\big\{ [(\mf{p},k)]
		\colon
		n \in \mf{p},\ k \in \{s, \ell\}
		\big\}
	\]
	and $V_{C_p}(f)$ has points determined by
	\[
		V_{C_p}(f) =
		\big\{
		[(\mf{p},\ell)]
		:
		f \in (1 - x) + \mf{p}\cdot\Z[x]/(x^p-1)
		\big\}
		\sqcup
		\big\{
		[(\mf{p},s)]
		:
		\mf{p} \neq \langle p \rangle,\  f \in \mf{p}\cdot\Z[x]/(x^p-1)
		\big\}.
	\]
\end{thm}

\begin{proof}
	Throughout this proof, $q$ denotes a prime not equal to $p$. Since $\nm(k) = k$ in $\ghost(\RU)$, every prime ideal of $\RU$ is of one of the following forms:
	\begin{align*}
		\ghostmap^*(0;0), \quad \quad                            & \ghostmap^*(0;\Z),      \\
		\ghostmap^*(p;\langle p,\xi-1 \rangle), \quad \quad      & \ghostmap^*(q;\Z[\xi]), \\
		\ghostmap^*(q;\langle q,\xi-\alpha \rangle), \quad \quad & \ghostmap^*(q;\Z),
	\end{align*}
	where $\alpha$ is a primitive $p$-th root of unity in $\F_q$. Using the definition of the ghost map, it is straightforward to show that
	\begin{align*}
		\ghostmap^*(0,0)                             & =\begin{pmatrix} 0 \\ 0 \end{pmatrix},           & \ghostmap^*(0,\Z[\xi]) & =\begin{pmatrix} 1-\gamma \\ 0 \end{pmatrix},     \\
		\ghostmap^*(p; \langle p,\xi-1 \rangle)      & =\begin{pmatrix} p, 1-\gamma \\ p \end{pmatrix}, & \ghostmap^*(p,\Z[\xi]) & =\begin{pmatrix} p, 1-\gamma \\ p \end{pmatrix}, \\
		\ghostmap^*(q; \langle q,\xi-\alpha \rangle) & \supseteq \begin{pmatrix} q \\ q \end{pmatrix},  & \ghostmap^*(q,\Z[\xi]) & =\begin{pmatrix} q, 1-\gamma \\ q \end{pmatrix},
	\end{align*}
	where $\binom{a}{b}$ denotes the ideal $\mf{p}$ with $\mf{p}(C_p/C_p)= \langle a\rangle$ and $\mf{p}(C_p/e) = \langle b \rangle$.
	Evidently, $\ghostmap^*(q; \langle q,\xi-\alpha \rangle)(C_p/e) = q$. In fact, we claim that $\ghostmap^*(q; \langle q,\xi-\alpha \rangle)(C_p/C_p) = \langle q \rangle$, so the inclusion in the bottom left is actually an equality. To see this, observe that we have a sequence of inclusions of prime ideals
	$$ \begin{pmatrix} 0 \\ 0 \end{pmatrix} \subset \begin{pmatrix} q \\ q \end{pmatrix} \subseteq \ghostmap^*(q; \langle q,\xi-\alpha \rangle) \subset \begin{pmatrix} q, 1-\gamma \\ q \end{pmatrix}$$
	where we claim that the outer inclusions are proper. Proper containment on the left is clear. The right-hand containment is proper since $1-\gamma \not\in \ghostmap^*(q;\langle q,\xi-\alpha \rangle)(C_p/C_p)$. Indeed, $\ghostmap(1-\gamma) = (0,1-\xi)$, and $1 - \xi \not\in \langle q,\xi-\alpha \rangle$, since if it were, we would have $\Z[\xi]/\langle q,\xi-\alpha \rangle = 0 \neq \F_q$. Since the Krull dimension of $\RU$ is $2$ by \cref{example: Krull dim of RU}, we obtain the claimed description of $\ghostmap^*(q; \langle q,\xi-\alpha \rangle)$. This gives the claimed set-theoretic description of $\Spec(\RU)$, where $\ghostmap^*(\mf{p};\mf{q})$ corresponds to $(\mf{p},s)$ and $\ghostmap^*(\mf{p}; \Z[\xi])$ corresponds to $(\mf{p},\ell)$.

	To describe the topology on $\Spec(\RU)$, we recall that it has a subbasis of closed sets determined by principal ideals. Principal ideals of $\RU$ are either generated by an element $n \in \RU(C_p/e) = \Z$, and we denote the corresponding  closed set by $V_e(n)$, or by an element $f \in \RU(C_p/C_p) = \Z[C_p] \cong \Z[x]/(x^p-1)$, and we denote the corresponding closed set by $V_{C_p}(f)$. We then use the levelwise description of the prime ideals to determine the elements of these  closed sets in the subbasis.
\end{proof}
\begin{remark}
	In fact, the linearization map $\uA \to \RU$ of $C_p$-Tambara functors induces a homeomorphism $\Spec(\RU) \cong \Spec(\uA)$. The linearization map is an isomorphism of $C_p$-Tambara functors if $p = 2$, but not if $p$ is an odd prime. Thus when $p$ is odd, we obtain examples of non-isomorphic $C_p$-Tambara functors with homeomorphic Nakaoka spectra.
\end{remark}


\subsection{Free \texorpdfstring{$C_p$}{Cₚ}-Tambara Functors}
\label{sec: main example}
Explicit descriptions of $\Spec(\uA[C_p/e])$ and $\Spec(\uA[C_p/C_p])$ are beyond reach, essentially because affine $n$-space $\A^n_{\mathbb{Z}} = \Spec(\Z[x_1,\ldots,x_n])$ does not admit an explicit description for $n>1$. However, the results above allow us to say some things about these spectra.

The strategy here is to study the ideals of the ghost $C_p$-Tambara functors $\ghost(\uA[C_p/C_p])$ and $\ghost(\uA[C_p/e])$, and use our understanding of the ideals of the ghost (\cref{lemma: when ideal pairs are prime in ghost} and \cref{Thm:GhostSurj}) to study their preimages under the ghost map.

Recall the terms ``underlying" and ``fixed" from \cref{notation:underlying/fixed}. Since we only deal with $C_p$-Tambara functors here, we will use these terms freely throughout this section.

\subsubsection{The Prime Ideal Spectrum of \texorpdfstring{$\uA[C_p/C_p]$}{A[Cₚ/Cₚ]}}

Let $T = \uA[C_p/C_p]$, and let $\ghost(T)$ be the ghost functor. Recall the structure of $T$ from \cref{Ex:AxGG} and the structure of $\ghost(T)$ from \cref{example: ghost of fixed point polynomial functor}.

By \cref{prop: primes of ghost}, the prime ideals of $\ghost(T)$ come in two forms:
\begin{enumerate}[(1)]
	\item $(\fa; \Z[x,n])$, where $\fa \subseteq \Z[x]$ is a $C_p$-prime.
	\item $(\nm^{-1}(\fb); \fb)$, where $\fb \subseteq \Z[x,n]$ is prime.
\end{enumerate}
Given that the $C_p$-action on $T(C_p/e) = \Z[x]$ is trivial, we replace (1) with:
\begin{enumerate}[($1'$)]
	\item $(\fa; \Z[x,n])$, where $\fa \subseteq \Z[x]$ is prime.
\end{enumerate}
Explicitly, ideals of the form ($1'$) look like
\[
	\begin{tikzcd}[row sep = large]
		\fa \times \Z[x,n] \ar[d,"\res" description]\\
		\fa \ar[u, bend left=50, "\tr"]
		\ar[u, bend right=50, "\nm"']
		\arrow[from=2-1, to=2-1, loop, in=300, out=240, distance=5mm, "{\mathrm{trivial}}"']
	\end{tikzcd}
\]
where $\fa = \langle q, f \rangle$ for some prime $q \in \Z$ and $f \in \Z[x]$ is a monic polynomial that is irreducible modulo $q$, or $\fa = \langle f \rangle$ for some irreducible $f \in \Z[x]$.

Ideals of the form (2) look like
\[
	\begin{tikzcd}[row sep = large]
		\nm^{-1}(\fb) \times \fb
		\ar[d,"\res" description]
		\\
		\nm^{-1}(\fb)
		\ar[u, bend left=50, "\tr"]
		\ar[u, bend right=50, "\nm"']
		\arrow[from=2-1, to=2-1, loop, in=300, out=240, distance=5mm, "{\mathrm{trivial}}"']
	\end{tikzcd}
\]
where $\fb \subseteq \Z[x,n]$ is a prime ideal. Making ideals of this form more explicit would require a complete understanding of $\Spec(\Z[x,n])$, which we do not have.

\begin{remark}\label{remark: A2 is hard to compute}
	It is worth taking some time to explain why one cannot hope to produce an explicit description of the points in $\Spec(\Z[x,n])$ (or higher-dimensional affine spaces over $\Z$), whereas explicit descriptions of $\Spec(\Z)$ and $\Spec(\Z[x])$ are well-known.

	To understand $\Spec(\Z[x])$, we consider the map $\Spec(\Z[x]) \to \Spec(\Z)$ coming from the inclusion $\Z \to \Z[x]$. Given a point $\langle p \rangle \in \Spec(\Z)$, we have a pullback square of topological spaces
	\[
		\begin{tikzcd}
			\Spec(\Z[x] \otimes \kappa(p))
			\ar[r]
			\ar[d]
			\ar[dr, "\lrcorner" description, very near start, phantom]
			&
			\{\langle p \rangle\}
			\ar[d, hook]
			\\
			\Spec(\Z[x])
			\ar[r]
			&
			\Spec(\Z)
		\end{tikzcd}
	\]

	where $\kappa(p)$ denotes the residue field of the point $\langle p \rangle \in \Spec(\Z)$, i.e.\ the field of fractions of $\Z/p$. This residue field is either $\F_p$ (if $p$ is a prime number) or $\QQ$ (if $p$ is $0$). In either case, we have that $\Spec(\Z[x] \otimes \kappa(p)) \cong \Spec(\kappa(p)[x])$ has one point for every monic irreducible polynomial in $\kappa(p)[x]$. Moreover, because $\Z/p$ is always a UFD, these monic irreducible polynomials in $\kappa(p)[x]$ have well-behaved lifts to $\Z[x]$. This gives a concrete description of the fiber of $\Spec(\Z[x]) \to \Spec(\Z)$ above each point in $\Spec(\Z)$, and thus we obtain the standard description of $\Spec(\Z[x])$.

	It is this last ``lifting'' step that breaks down when we attempt to study $\Spec(\Z[x,n])$ by comparison to $\Spec(\Z[x])$. The attempt would look like this: we take a prime ideal of $\Z[x]$, for example $\langle x^2+5 \rangle$. Then we obtain a pullback square of topological spaces
	\[
		\begin{tikzcd}
			\Spec(\Z[x,n] \otimes_{\Z[x]} \kappa(x^2+5))
			\ar[r]
			\ar[d]
			\ar[dr, "\lrcorner" description, very near start, phantom]
			&
			\{\langle x^2+5 \rangle\}
			\ar[d, hook]
			\\
			\Spec(\Z[x,n])
			\ar[r]
			&
			\Spec(\Z[x])
		\end{tikzcd}
	\]
	Once again, $\Spec(\Z[x,n] \otimes_{\Z[x]} \kappa(x^2+5)) \cong \Spec(\kappa(x^2-5)[n])$ has one point for each monic irreducible polynomial in $\kappa(x^2-5)[n]$. However, the ring $\Z[x]/\langle x^2+5 \rangle$ is famously \emph{not} a UFD! Thus, it is not clear that there is any canonical way to lift these monic polynomials to $\Z[x,n]$. In order to clear this obstruction and obtain some uniform description of the prime ideals of $\Z[x,n]$, we would need some global control over the class groups of all number fields, which is beyond hope.
\end{remark}

Even though explicitly calculating all ideals of the form (2) is beyond our reach, we can still understand preimages of ideals of the form ($1'$) under the ghost map $T \to \ghost(T)$. Let $(\fa; \Z[x,n])$ be a prime ideal of the form ($1'$) in $\ghost(T)$.

On the underlying level, the ghost map is the identity. So the preimage of $(\fa; \Z[x,n])$ on the underlying level is just $\fa \subseteq \Z[x] = T(C_p/e)$.

On the fixed level, things are a little more complicated. Given any $g(t,x,n) \in T(C_p/C_p)$, we may use the relations $t^2 = pt$ and $tn = tx^p$ in $T(C_p/C_p)$ to rewrite
\[
	g(t,x,n) = g_0(x,n) + tg_1(x)
\]
for some $g_0 \in \Z[x,n]$ and $g_1 \in \Z[x]$. If this is to be in the preimage of $(\fa, \Z[x,n])$ under the ghost map, then we must have
\[
	\res(g(t,x,n)) = g(p,x,x^p) \in \fa
	\qquad \text{ and } \qquad
	g(t,x,n) + \tau = g_0(x,n) \in \Z[x,n].
\]
This second condition is vacuously true, so we focus on the first, which becomes
\begin{align*}
	\res(g(t,x,n)) \in \fa & \iff g(p,x,x^p) = g_0(x,x^p) + pg_1(x) \in \big(f(x)\big).
\end{align*}
In this situation, the preimage of $(\fa; \Z[x,n])$ at the top level is those polynomials $g(t,x,n)$ such that $g_0(x,x^p) + p g_1(x) \in \fa$. This condition does not simplify substantially for either of the two forms of prime ideal in $\Z[x]$: either $\fa = \langle q, f(x) \rangle$ or $\fa = \langle f(x) \rangle$.

\begin{example}\label{example: coincidence1}
	Consider the ideal $(\langle p, x \rangle; \Z[x,n])$ of $\ghost(T)$ for $T$ the $C_p$-Tambara functor $\uA[C_p/C_p]$. On the underlying level, the preimage is $\langle p,x \rangle$ as above. On the fixed level, the preimage is generated by polynomials
	\[
		g(t,x,n) = g_0(x,n) + tg_1(x)
	\]
	which restrict (in $\ghost(T)$) to $g_0(x,x^p) + pg_1(x) \in \langle p,x \rangle$. Equivalently, such a polynomial must have $g_0(x,x^p) \in \langle p,x \rangle$, which means that the constant term of $g_0(x,n)$ must be divisible by $p$. So the preimage is a single ideal of $T$, which looks like
	\[
		\begin{tikzcd}[row sep = large]
			\big\{ 	g_0(x,n) + tg_1(x) \in T(C_p/C_p)
			\mid
			p \text{ divides } g_0(0,0)
			\big\}
			\ar[d,"\res" description]
			\\
			\langle p,x \rangle.
			\ar[u, bend left=50, "\tr"]
			\ar[u, bend right=50, "\nm"']
			\arrow[from=2-1, to=2-1, loop, in=300, out=240, distance=5mm, "{\mathrm{trivial}}"']
		\end{tikzcd}
	\]
\end{example}

\begin{example}\label{example: coincidence2}
	Consider the prime ideal $\fb = \langle p,x,n \rangle$ of $\Z[x,n]$. This corresponds to a prime ideal of $\ghost(T)$ of type (2), namely $(\nm^{-1}\langle p,x,n \rangle; \langle p,x,n \rangle)$. In this case, the norm in $\ghost(T)$ is determined by $f(x) \mapsto f(n)$, so
	\[
		\nm^{-1}\langle p,x,n \rangle = \langle p,x \rangle.
	\]
	What is the preimage of this prime under the ghost map?

	On the underlying level, it is again $\langle p,x \rangle$.

	On the fixed level, a polynomial
	\[
		g(t,x,n) = g_0(x,n) + tg_1(x) \in T(C_p/C_p)
	\]
	is in the preimage of the ghost map if
	\[
		\res(g(t,x,n)) = g_0(x,x^p) + pg_1(x) \in \langle p,x\rangle
		\qquad \text{ and } \qquad
		g(t,x,n) + \tau = g_0(x,n) \in \langle p,x,n \rangle.
	\]
	Both of these conditions are met if and only if the constant term of $g_0$ is divisible by $p$. So we see that this prime $(\langle p,x \rangle;\langle p,x,n \rangle)$ and the prime $(\langle p,x \rangle;\Z[x,n])$ of the previous example are both primes of $\ghost(T)$ lying over the same prime in $T$.
\end{example}

We have described the prime ideals of $\uA[C_p/C_p]$, but not uniquely; a priori, primes of the form $\ghostmap^*(\mf{a}; \Z[x,n])$ may coincide with primes of the form $\ghostmap^*(\nm^{-1} \mf{b}; \mf{b})$, for example. Additionally, in order to understand the topology on $\Spec(\uA[C_p/C_p])$ (in light of \Cref{thm: topology in noetherian case}), we must describe its poset structure. Thus, it is also important to determine the containments between primes of the forms (1) and (2). There are four possible types of containments, which we address in the following four propositions. In what follows, let $\mf{a}$ and $\mf{a}'$ be $C_p$-prime ideals of $\uA[C_p/C_p](C_p/e) = \Z[x]$, and let $\mf{b}$ and $\mf{b}'$ be prime ideals of $\Phi^{C_p}\uA[C_p/C_p] = \Z[x,n]$.

\begin{proposition}
	We have a containment $\ghostmap^*(\mf{a}; \Z[x,n]) \subseteq \ghostmap^*(\mf{a}'; \Z[x,n])$ if and only if $\mf{a} \subseteq \mf{a}'$.
\end{proposition}
\begin{proof}
	The underlying level of $\ghostmap^*(\mf{a}; \Z[x,n])$ is $\mf{a}$, so $\ghostmap^*(\mf{a}; \Z[x,n]) \subseteq \ghostmap^*(\mf{a}'; \Z[x,n])$ implies $\mf{a} \subseteq \mf{a}'$. The reverse implication is trivial.
\end{proof}

\begin{proposition}
	If $\mf{b} \subseteq \mf{b}'$ then $\ghostmap^*(\nm^{-1} \mf{b}; \mf{b}) \subseteq \ghostmap^*(\nm^{-1} \mf{b}'; \mf{b}')$. If  $\ghostmap^*(\nm^{-1} \mf{b}; \mf{b}) \subseteq \ghostmap^*(\nm^{-1} \mf{b}'; \mf{b}')$ then $\mf{b} \subseteq \mf{b}'$ or $\langle p,n-x^p\rangle \subseteq \mf{b}'$.
\end{proposition}
\begin{proof}
	The first implication is trivial. For the second implication, suppose $\ghostmap^*(\nm^{-1} \mf{b}; \mf{b}) \subseteq \ghostmap^*(\nm^{-1} \mf{b}'; \mf{b}')$ and $\mf{b} \nsubseteq \mf{b}'$. Pick some $g(x,n) \in \mf{b} \setminus \mf{b}'$, and set $f_1(t,x,n) \coloneqq (n-x^p)g(x,n)$. We have $f_1(p,n,n^p) = 0 \in \mf{b}$ and $f_1(0,x,n) = (n-x^p)g(x,n) \in \mf{b}$, so $f_1 \in \ghostmap^*(\nm^{-1} \mf{b}; \mf{b})(C_p/C_p)$. Thus, $f_1 \in \ghostmap^*(\nm^{-1} \mf{b}'; \mf{b}')(C_p/C_p)$, telling us that $(n-x^p)g(x,n) = f_1(0,x,n) \in \mf{b}'$. Since $g \notin \mf{b}'$ and $\mf{b}'$ is a prime ideal, we conclude that $n-x^p \in \mf{b}'$. Similarly, we define $f_2(t,x,n) \coloneqq (t-p)g(x,n)$. Then $f_2 \in \ghostmap^*(\nm^{-1} \mf{b}; \mf{b})(C_p/C_p)$, so $f_2 \in \ghostmap^*(\nm^{-1} \mf{b}'; \mf{b}')(C_p/C_p)$, so $(-pg(x,n)) = f_2(0,x,n) \in \mf{b}'$. Since $g \notin \mf{b}'$ and $\mf{b}'$ is prime, we conclude that $p \in \mf{b}'$.
\end{proof}

\begin{corollary}
	We have $\ghostmap^*(\nm^{-1} \mf{b}; \mf{b}) = \ghostmap^*(\nm^{-1} \mf{b}'; \mf{b}')$ if and only if $\mf{b} = \mf{b}'$.
\end{corollary}
\begin{proof}
	The reverse implication is trivial. For the forward implication, suppose for contradiction that $\ghostmap^*(\nm^{-1} \mf{b}; \mf{b}) = \ghostmap^*(\nm^{-1} \mf{b}'; \mf{b}')$ and $\mf{b} \neq \mf{b}'$. Without loss of generality, $\mf{b} \nsubseteq \mf{b}'$, so $\langle p,n-x^p \rangle \subseteq \mf{b}'$. Then
	\[
		\langle p,n-x^p \rangle \subseteq \ghostmap^*(\nm^{-1} \mf{b}'; \mf{b}')(C_p/C_p)
		= \ghostmap^*(\nm^{-1} \mf{b}; \mf{b})(C_p/C_p),
	\]
	so $\langle p,n-x^p \rangle \subseteq \mf{b}$. Now $\mf{b}$ is a prime ideal of $\Z[x,n]$ which contains $\langle p,n-x^p \rangle$; such ideals correspond bijectively with prime ideals of $\Z[x,n]/\langle p,n-x^p \rangle \cong \F_p[x]$. We know $\F_p[x]$ is a PID, with its nonzero prime ideals being generated by monic irreducible polynomials. Any monic polynomial in $\F_p[x]$ can be lifted to a monic polynomial in $\Z[x]$, so we conclude that $\mf{b} = \langle p,n-x^p \rangle$ or $\mf{b} = \langle p,n-x^p,f \rangle$ for some monic polynomial $f \in \Z[x]$ such that $f$ is irreducible modulo $p$. The former is not possible, since we have assumed $\mf{b} \nsubseteq \mf{b}'$. Thus, $\mf{b} = \langle p,n-x^p,f \rangle$ for some monic polynomial $f \in \Z[x]$ such that $f$ is irreducible modulo $p$. Now $\nm^{-1} \mf{b}$ is the kernel of
	\[
		\Z[x] \xrightarrow{x \mapsto n} \Z[x,n] \xrightarrow{x \mapsto x, n \mapsto x^p} \F_p[x] \xrightarrow{x \mapsto \alpha} \F_{p^d}
	\]
	where $d$ is the degree of $f$ and $\alpha$ is a generator of $\F_{p^d}$ with minimal polynomial $f$. In particular, we see that $f \in \nm^{-1} \mf{b}$, because
	\[
		f(\alpha) = 0 \implies f(\alpha^p) = 0
	\]
	since $\alpha$ and $\alpha^p$ are Galois conjugates in $\F_{p^d}$. Thus, $f \in \ghostmap^*(\nm^{-1} \mf{b}; \mf{b})(C_p/C_p) = \ghostmap^*(\nm^{-1} \mf{b}'; \mf{b}')(C_p/C_p)$, so $f \in \mf{b}'$. But now we have $\mf{b} = \langle p,n-x^p,f \rangle \subseteq \mf{b}'$, a contradiction.
\end{proof}

\begin{proposition}
	We have a containment $\ghostmap^*(\mf{a}; \Z[x,n]) \subseteq \ghostmap^*(\nm^{-1} \mf{b}; \mf{b})$ if and only if one of the following holds:
	\begin{enumerate}[(a)]
		\item $\mf{b} = \langle p,n-x^p \rangle$ and $\mf{a} \subseteq \langle p \rangle$, or
		\item $\mf{b} = \langle p,n-x^p,f \rangle$ and $\mf{a} \subseteq \langle p,f \rangle$ for some monic polynomial $f \in \Z[x]$ such that $f$ is irreducible mod $p$.
	\end{enumerate}
\end{proposition}
\begin{proof}
	First, suppose $\ghostmap^*(\mf{a}; \Z[x,n]) \subseteq \ghostmap^*(\nm^{-1} \mf{b}; \mf{b})$. We have $t-p \in \ghostmap^*(\mf{a}; \Z[x,n])(C_p/C_p)$, so $t-p \in \ghostmap^*(\nm^{-1} \mf{b}; \mf{b})(C_p/C_p)$, which implies $p \in \mf{b}$. Likewise, $n-x^p \in \ghostmap^*(\mf{a}; \Z[x,n])(C_p/C_p)$, so $n-x^p \in \mf{b}$. Now $\mf{b}$ is a prime ideal of $\Z[x,n]$ which contains $\langle p,n-x^p \rangle$; as in the previous corollary, we conclude that $\mf{b} = \langle p,n-x^p\rangle$ or $\mf{b} = \langle p,n-x^p,f\rangle$ for some monic polynomial $f \in \Z[x]$ such that $f$ is irreducible modulo $p$. In the first case, we have $\nm^{-1} \mf{b} = \nm^{-1} \langle p,n-x^p \rangle =  \langle p \rangle$, whence $\mf{a} \subseteq \langle p \rangle$, exactly as desired for (1). In the second case, we have as before that $p \in \nm^{-1} \mf{b}$ and that $f \in \nm^{-1} \mf{b}$. Now since $\langle p,f \rangle$ is a maximal ideal of $\Z[x]$, we conclude that $\nm^{-1} \mf{b} = \langle p,f \rangle$. Thus, $\mf{a} \subseteq \langle p,f \rangle$, exactly as desired for (2).

	In the other direction, suppose first that $\mf{b} = \langle p,n-x^p \rangle$ and $\mf{a} \subseteq \langle p \rangle$. It suffices to show that $\ghostmap^*(\langle p \rangle; \Z[x,n]) \subseteq \ghostmap^*(\nm^{-1} \mf{b}; \mf{b})$, which can be checked directly.

	Likewise, suppose that $\mf{b} = \langle p,n-x^p,f \rangle$ and $\mf{a} \subseteq \langle p,f \rangle$ for some monic polynomial $f \in \Z[x]$ such that $f$ is irreducible modulo $p$. It suffices to show that $\ghostmap^*(\langle p,f \rangle; \Z[x,n]) \subseteq \ghostmap^*(\nm^{-1} \mf{b}; \mf{b})$, which can be checked directly.
\end{proof}

\begin{proposition}
	We have a containment $\ghostmap^*(\nm^{-1} \mf{b}; \mf{b}) \subseteq \ghostmap^*(\mf{a}; \Z[x,n])$ if and only if $f(x) \in \mf{a}$ for all $f(n) \in \mf{b} \cap \Z[n]$.
\end{proposition}
\begin{proof}
	First, suppose we have $\ghostmap^*(\nm^{-1} \mf{b}; \mf{b}) \subseteq \ghostmap^*(\mf{a}; \Z[x,n])$. At the underlying level, this says precisely that $f(x) \in \mf{a}$ for all $f(n) \in \mf{b} \cap \Z[n]$.

	In the other direction, suppose $f(x) \in \mf{a}$ for all $f(n) \in \mf{b} \cap \Z[n]$. This tells us that $\ghostmap^*(\nm^{-1} \mf{b}; \mf{b})(C_p/e) \subseteq \ghostmap^*(\mf{a}; \Z[x,n])(C_p/e)$, and only containment at the fixed level remains to be shown. So, let $f(t,x,n) \in \ghostmap^*(\nm^{-1} \mf{b}; \mf{b})(C_p/C_p)$ be arbitrary. Then, in particular, we have $f(p,n,n^p) \in \mf{b}$. Setting $h(x) \coloneqq f(p,x,x^p)$, this says that $h(n) \in \mf{b}$, so $h(x) \in \mf{a}$. Thus, $f(p,x,x^p) \in \mf{a}$, which tells us that $f(t,x,n) \in \ghostmap^*(\mf{a}; \Z[x,n])(C_p/C_p)$.
\end{proof}

\begin{corollary}
	We have equality $\ghostmap^*(\mf{a}; \Z[x,n]) = \ghostmap^*(\nm^{-1} \mf{b}; \mf{b})$ if and only if one of the following holds:
	\begin{enumerate}[(a)]
		\item $\mf{b} = \langle p,n-x^p \rangle$ and $\mf{a} = \langle p\rangle$, or
		\item $\mf{b} = \langle p,n-x^p,f \rangle$ and $\mf{a} = \langle p,f \rangle$ for some monic polynomial $f \in \Z[x]$ that is irreducible modulo $p$.
	\end{enumerate}
\end{corollary}
\begin{proof}
	Combine the prior two propositions with the fact that $\nm^{-1}\langle p,n-x^p,f\rangle = \langle p,f \rangle$.
\end{proof}

Note that \Cref{example: coincidence1,example: coincidence2} demonstrate this coincidence in the case $\mf{b} = \langle p,n-x^p,x \rangle = \langle p,x,n \rangle$. Overall, we see that:

\begin{thm}[{The structure of $\Spec(\uA[C_p/C_p])$ as a set}]\label{thm:set structure C_p/C_p}
	Let $\mf{a}$ be a $C_p$-prime ideal of $\uA[C_p/C_p](C_p/e) = \Z[x]$, and let $\mf{b}$ be a prime ideal of $\Phi^{C_p}\uA[C_p/C_p] = \Z[x,n]$.
	\begin{enumerate}[(a)]
		\item Distinct primes in $\ghost(\uA[C_p/C_p])$ of the form $(\mf{a}; \Z[x,n])$ remain distinct after applying $\ghostmap^*$.
		\item Distinct primes in $\ghost(\uA[C_p/C_p])$ of the form $(\nm^{-1} \mf{b}; \mf{b})$ remain distinct after applying $\ghostmap^*$.
		\item There is a coincidence $\ghostmap^*(\langle p \rangle; \Z[x,n]) = \ghostmap^*(\nm^{-1} \langle p,n-x^p \rangle; \langle p,n-x^p \rangle)$.
		\item There is a coincidence $\ghostmap^*(\langle p,f \rangle; \Z[x,n]) = \ghostmap^*(\nm^{-1} \langle p,n-x^p,f \rangle; \langle p,n-x^p,f \rangle)$ for each monic polynomial $f \in \Z[x]$ such that $f$ is irreducible modulo $p$.
		\item There are no other coincidences.
	\end{enumerate}
\end{thm}

\subsubsection{The Prime Ideal Spectrum of \texorpdfstring{$\uA[C_p/e]$}{A[Cₚ/e]}} Recall the structure of $\uA[C_p/e]$ from \Cref{Ex:AxGe} and the structure of $\ghost(\uA[C_p/e])$ from \Cref{Ex:Ghost of AxGe}.

By \cref{prop: primes of ghost}, the prime ideals of $\ghost(\uA[C_p/e])$ come in two forms:
\begin{enumerate}[(1)]
	\item $(\fa; \Z[n])$, where $\fa \subseteq \Z[x_0,x_1, \ldots, x_{p-1}]$ is a $C_p$-prime ideal.
	\item $(\nm^{-1}(\fb); \fb)$, where $\fb \subseteq \Z[n]$ is a prime ideal.
\end{enumerate}

Explicitly, ideals of the form $(1)$ look like

\[
	\begin{tikzcd}[row sep=large]
		(\fa)^{C_p} \times \Z[n]
		\ar[d, "\res" description]
		\\
		\fa
		\ar[u, bend left=50, "\tr"]
		\ar[u, bend right=50, "\nm"']
		\arrow[from=2-1, to=2-1, loop, in=300, out=240, distance=5mm, "\gamma"']
	\end{tikzcd}
\]

\noindent where $\fa$ is a $C_p$-prime ideal of $\Z[x_0,x_1, \ldots, x_{p-1}]$. As noted in \Cref{remark: A2 is hard to compute}, knowing the $C_p$-prime ideals of $\Z[x_0,x_1, \ldots, x_{p-1}]$ (i.e. the prime ideals of $\Z[x_0,x_1, \ldots, x_{p-1}]^{C_p}$) is quite difficult.

Explicitly, ideals of the form $(2)$ look like

\[
	\begin{tikzcd}[row sep=large]
		(\nm^{-1}(\fb))^{C_p} \times \fb
		\ar[d, "\res" description]
		\\
		\nm^{-1}(\fb)
		\ar[u, bend left=50, "\tr"]
		\ar[u, bend right=50, "\nm"']
		\arrow[from=2-1, to=2-1, loop, in=300, out=240, distance=5mm, "\gamma"']
	\end{tikzcd}
\]

\noindent for $\fb = \langle q,f \rangle$ for some prime $q \in \Z$ and $f(n) \in \Z[n]$ a monic polynomial that is irreducible modulo $q$, or $\fb = \langle f \rangle$ for $f(n) \in \Z[n]$ an irreducible polynomial.

Applying $\ghostmap^*$ to these ideals produces all of the prime ideals of $\uA[C_p/e]$.
Again, we must classify the containments between these prime ideals. In what follows, let $\mf{a}$ and $\mf{a}'$ be $C_p$-prime ideals of $\uA[C_p/e](C_p/e) = \Z[x_0, \ldots, x_{p-1}]$, and let $\mf{b}$ and $\mf{b}'$ be prime ideals of $\Phi^{C_p}\uA[C_p/e] = \Z[n]$.

\begin{proposition}
	We have a containment $\ghostmap^*(\mf{a}; \Z[n]) \subseteq \ghostmap^*(\mf{a}'; \Z[n])$ if and only if $\mf{a} \subseteq \mf{a}'$.
\end{proposition}
\begin{proof}
	The underlying level of $\ghostmap^*(\mf{a}; \Z[x,n])$ is $\mf{a}$, so $\ghostmap^*(\mf{a}; \Z[n]) \subseteq \ghostmap^*(\mf{a}'; \Z[n])$ implies $\mf{a} \subseteq \mf{a}'$. The reverse implication is trivial.
\end{proof}

\begin{proposition}
	If $\mf{b} \subseteq \mf{b}'$ then $\ghostmap^*(\nm^{-1} \mf{b}; \mf{b}) \subseteq \ghostmap^*(\nm^{-1} \mf{b}'; \mf{b}')$. If  $\ghostmap^*(\nm^{-1} \mf{b}; \mf{b}) \subseteq \ghostmap^*(\nm^{-1} \mf{b}'; \mf{b}')$ then $\mf{b} \subseteq \mf{b}'$ or $p \in \mf{b}'$.
\end{proposition}
\begin{proof}
	The first implication is trivial. For the second implication, suppose $\ghostmap^*(\nm^{-1} \mf{b}; \mf{b}) \subseteq \ghostmap^*(\nm^{-1} \mf{b}'; \mf{b}')$ and $\mf{b} \nsubseteq \mf{b}'$. Pick some $g(n) \in \mf{b} \setminus \mf{b}'$, which we also view as an element of $\uA[C_p/e](C_p/C_p)$. Then $\res((t_0 - p)g(n)) = 0$ and $-pg(n) \in \mf{b}$, so $(t_0-p)g(n) \in \ghostmap^*(\nm^{-1} \mf{b}; \mf{b})(C_p/C_p)$. Thus, $(t_0-p)g(n) \in \ghostmap^*(\nm^{-1} \mf{b}'; \mf{b}')(C_p/C_p)$, so $-pg(n) \in \mf{b}'$. Since $\mf{b}'$ is prime and $g(n) \notin \mf{b}'$, we conclude that $p \in \mf{b}'$.
\end{proof}

\begin{corollary}
	We have an equality $\ghostmap^*(\nm^{-1} \mf{b}; \mf{b}) = \ghostmap^*(\nm^{-1} \mf{b}'; \mf{b}')$ if and only if $\mf{b} = \mf{b}'$.
\end{corollary}
\begin{proof}
	The reverse implication is trivial. For the forward implication, suppose for contradiction that $\ghostmap^*(\nm^{-1} \mf{b}; \mf{b}) = \ghostmap^*(\nm^{-1} \mf{b}'; \mf{b}')$ and $\mf{b} \neq \mf{b}'$. Without loss of generality, $\mf{b} \nsubseteq \mf{b}'$, so $p \in \mf{b}'$. Then we also have $p \in \nm^{-1} \mf{b}' = \nm^{-1}(\mf{b})$, so $p = \nm(p) \in \mf{b}$.


	Since $\mf{b}$ is a prime ideal of $\Z[n]$, we conclude that $\mf{b} = \langle p \rangle$ or $\mf{b} = \langle p,f(n) \rangle$ for some monic polynomial $f(n) \in \Z[n]$ which is irreducible modulo $p$. The former would imply $\mf{b} \subseteq \mf{b}'$, which we have assumed is not the case; thus $\mf{b} = \langle p,f(n) \rangle$. Now we claim that $\res(f(n)) \in \nm^{-1} \mf{b}$; to see this, consider the composite
	\[
		\Z[n]
		\xrightarrow{n \mapsto x_0 \cdots x_{p-1}}
		\Z[x_0, \dots, x_{p-1}]
		\xrightarrow{x_i \mapsto n}
		\F_p[n]
		\xrightarrow{n \mapsto \alpha}
		\F_{p^d}
	\]
	where $d$ is the degree of $f$ and $\alpha$ is a generator of $\F_{p^d}$ with minimal polynomial $f$. This composite map sends $f(n)$ to $f(\alpha^p)$, which equals to $0$ because $\alpha$ and $\alpha^p$ are Galois conjugates in $\F_{p^d}$. Unwinding definitions, this says precisely that $\res(f) \in \nm^{-1} \mf{b}$. Thus, $f(n) \in \ghostmap^*(\nm^{-1} \mf{b}; \mf{b}) = \ghostmap^*(\nm^{-1} \mf{b}'; \mf{b}')$, and so $f(n) \in \mf{b}'$. But now $\mf{b} \subseteq \mf{b}'$, contradicting our prior assumption.
\end{proof}

\begin{proposition}
	We have a containment $\ghostmap^*(\mf{a}; \Z[n]) \subseteq \ghostmap^*(\nm^{-1} \mf{b}; \mf{b})$ if and only if one of the following holds:
	\begin{enumerate}[(a)]
		\item $\mf{b} = \langle p \rangle$ and $\mf{a} \subseteq \langle p \rangle + \varepsilon$, or
		\item $\mf{b} = \langle p,f(n) \rangle$ and $\mf{a} \subseteq \langle p,f(x_0) \rangle + \varepsilon$ for some monic polynomial $f(n) \in \Z[n]$ which is irreducible modulo $p$,
	\end{enumerate}
	where $\varepsilon$ is the ideal $\langle x_0 - x_1, x_1 - x_2, \dots, x_{p-2} - x_{p-1} \rangle$ of $\Z[x_0, \dots, x_{p-1}]$.
\end{proposition}
\begin{proof}
	We note that $t_0-p \in \ghostmap^*(\mf{a}; \Z[n])(C_p/C_p) \subseteq \ghostmap^*(\nm^{-1} \mf{b}; \mf{b})$, so $p \in \mf{b}$. Now since $\mf{b}$ is a prime ideal of $\Z[n]$, we have either $\mf{b} = \langle p \rangle$ or $\mf{b} = \langle p,f(n) \rangle$ for some monic polynomial $f(n) \in \Z[n]$ which is irreducible modulo $p$. In the first case, we have
	\[\mf{a} = \ghostmap^*(\mf{a}; \Z[n])(C_p/e) \subseteq \ghostmap^*(\nm^{-1} \mf{b}; \mf{b})(C_p/e) = \nm^{-1}\mf{b} = \langle p \rangle + \varepsilon,\]
	as desired. In the second case, we have
	\[\mf{a} = \ghostmap^*(\mf{a}; \Z[n])(C_p/e) \subseteq \ghostmap^*(\nm^{-1} \mf{b}; \mf{b})(C_p/e) = \nm^{-1}\mf{b} = \langle p, f(x_0) \rangle + \varepsilon,\]
	which is again exactly what we wanted.

	In the other direction, first suppose that $\mf{b} = \langle p \rangle$ and $\mf{a} \subseteq \langle p \rangle + \varepsilon$. We wish to show that $\ghostmap^*(\mf{a}; \Z[n]) \subseteq \ghostmap^*(\nm^{-1} \mf{b}; \mf{b})$, and we have by assumption that $\mf{a} \subseteq \nm^{-1} \mf{b}$, so we must only show containment at the fixed level. Thus, let $g \in \ghostmap^*(\mf{a}; \Z[n])(C_p/C_p) = \res^{-1} \mf{a}$ be arbitrary. We know that $\res(g) \in \nm^{-1} \mf{b}$, so we wish to show that $g + \tau \in \mf{b} = \langle p \rangle$. Equivalently, we want to show that $g$ is sent to a multiple of $p$ once we set each $t_{\vec v}$ to $0$. In particular, we may assume that $g=g(n)$ does not contain any of the $t_{\vec{v}}$.  Supposing this, we have that
	\[
		\nm(\res(g(n))) = \nm(g(x_0\cdots x_{p-1})) = g(n^p) + \tau
	\]
	is in $\mf{b} = \langle p \rangle$ and it follows that $g(n)$ was divisible by $p$.

	Finally, suppose that $\mf{b} = \langle p,f(n) \rangle$ and $\mf{a} \subseteq \langle p,f(x_0) \rangle + \varepsilon$. We wish to show that $\ghostmap^*(\mf{a}; \Z[n]) \subseteq \ghostmap^*(\nm^{-1} \mf{b}; \mf{b})$, and we have by assumption that $\mf{a} \subseteq \nm^{-1} \mf{b}$, so we must only show containment at the fixed level. Thus, let $g \in \ghostmap^*(\mf{a}; \Z[n])(C_p/C_p) = \res^{-1} \mf{a}$ be arbitrary. We know that $\res(g) \in \nm^{-1} \mf{b}$, so we wish to show that $g + \tau \in \mf{b} = \langle p, f(n) \rangle$. Equivalently, we want to show that $g$ is sent to zero by
	\[\Z[n][t_{\vec v} | \vec v \in \N^p]/I \xrightarrow{t_{\vec v} \mapsto 0} \F_p[n] \xrightarrow{n \mapsto \alpha} \F_{p^d}\]
	where $d$ is the degree of $f$ and $\alpha$ is a generator of $\F_{p^d}$ with minimal polynomial $f$. What we know is that $g$ is sent to zero by
	\[\Z[n][t_{\vec v} | \vec v \in \N^p] \xrightarrow{\res} \Z[x_0, \dots, x_{p-1}] \xrightarrow{x_i \mapsto x_0} \F_p[x_0]\xrightarrow{x_0 \mapsto \alpha} \F_{p^d}.\]
	The composite of the first two maps has the effect $n \mapsto x_0^p$, $t_{\vec v} \mapsto px_0^{\sum \vec v} = 0$. Thus, we may assume that $g$ does not contain any variables $t_{\vec v}$, i.e. we have a polynomial $g(n)$ in the single variable $n$. We have that $g(\alpha^p) = 0$, and since $\alpha$ and $\alpha^p$ are Galois conjugates in $\F_{p^d}$, we conclude that $g(\alpha) = 0$.
\end{proof}

\begin{proposition}
	We have a containment $\ghostmap^*(\nm^{-1} \mf{b}; \mf{b}) \subseteq \ghostmap^*(\mf{a}; \Z[n])$ if and only if $\nm^{-1} \mf{b} \subseteq \mf{a}$.
\end{proposition}
\begin{proof}
	The forward implication is trivial by taking $C_p/e$ levels. For the reverse implication, suppose $\nm^{-1} \mf{b} \subseteq \mf{a}$. We wish to show $\ghostmap^*(\nm^{-1} \mf{b}; \mf{b}) \subseteq \ghostmap^*(\mf{a}; \Z[n])$, and we already have (by assumption) the containment at the $C_p/e$ level, so it suffices to show that $\ghostmap^*(\nm^{-1} \mf{b}; \mf{b})(C_p/C_p) \subseteq \ghostmap^*(\mf{a}; \Z[n])(C_p/C_p)$. Thus, let $g \in \ghostmap^*(\nm^{-1} \mf{b}; \mf{b})(C_p/C_p)$ be arbitrary. Then $\res(g) \in \nm^{-1} \mf{b} \subseteq \mf{a}$, so $g \in \ghostmap^*(\mf{a}; \Z[n])(C_p/C_p)$, as desired.
\end{proof}

\begin{corollary}
	We have an equality $\ghostmap^*(\mf{a}; \Z[n]) = \ghostmap^*(\nm^{-1} \mf{b}; \mf{b})$ if and only if one of the following holds:
	\begin{enumerate}[(a)]
		\item $\mf{b} = \langle p \rangle$ and $\mf{a} = \langle p \rangle + \varepsilon$, or
		\item $\mf{b} = \langle p,f(n) \rangle$ and $\mf{a} = \langle p,f(x_0) \rangle + \varepsilon$ for some monic polynomial $f(n) \in \Z[n]$ which is irreducible modulo $p$.
	\end{enumerate}
\end{corollary}
\begin{proof}
	Combine the prior two propositions with the computations $\nm^{-1}\langle p \rangle = \langle p \rangle + \varepsilon$ and $\nm^{-1}\langle p,f(n) \rangle = \langle p,f(x_0) \rangle + \varepsilon$, noting that $\langle p \rangle + \varepsilon$ and $\langle p,f(x_0) \rangle + \varepsilon$ are $C_p$-prime ideals of $\Z[x_0, \dots, x_{p-1}]$ because they are $C_p$-invariant and prime.
\end{proof}

Overall, we see that:

\begin{thm}[{The structure of $\Spec(\uA[C_p/e])$ as a set}]\label{thm:set structure C_p/e}
	Let $\mf{a}$ be a $C_p$-prime ideal of $\uA[C_p/e](C_p/e) = \Z[x_0, \ldots, x_{p-1}]$ and let $\mf{b}$ be a prime ideal of $\Phi^{C_p}\uA[C_p/e] = \Z[n]$.
	\begin{enumerate}[(a)]
		\item Distinct primes in $\ghost(\uA[C_p/e])$ of the form $(\mf{a}; \Z[n])$ remain distinct after applying $\ghostmap^*$.
		\item Distinct primes in $\ghost(\uA[C_p/e])$ of the form $(\nm^{-1} \mf{b}; \mf{b})$ remain distinct after applying $\ghostmap^*$.
		\item There is a coincidence $\ghostmap^*(\langle p \rangle + \varepsilon; \Z[n]) = \ghostmap^*(\nm^{-1} \langle p \rangle; \langle p \rangle)$.
		\item There is a coincidence $\ghostmap^*(\langle p,f(x_0)\rangle  + \varepsilon; \Z[n]) = \ghostmap^*(\nm^{-1} \langle p,f(n)\rangle ; \langle p,f(n)\rangle)$ for each monic polynomial $f(n) \in \Z[n]$ such that $f(n)$ is irreducible modulo $p$.
		\item There are no other coincidences.
	\end{enumerate}
	In the above, $\varepsilon$ denotes the ideal $\langle x_0 - x_1, x_1 - x_2, \dots, x_{p-2} - x_{p-1} \rangle$ of $\Z[x_0, \dots, x_{p-1}]$.
\end{thm}

\subsection{The Tambara Affine Line}\label{SS:A1}

In this section we complete the computation of the $C_p$-Tambara affine line by computing the behavior on $\Spec$ of the structure maps $\cores$, $\conm$, $\cotr$, and $\coc$ of the co-Tambara object $\uA[{-}]$.  To simplify notation we will write $\Poly_{C_p}$ as $\Poly$ in this section. To recall, the co-Tambara structure on $\uA[{-}]$ consists of corestriction, conorm, and cotransfer morphisms
\begin{align*}
	\cores \colon & \uA[C_p/e] \to \uA[C_p/C_p], \\
	\conm \colon  & \uA[C_p/C_p] \to \uA[C_p/e], \\
	\cotr \colon  & \uA[C_p/C_p] \to \uA[C_p/e],
\end{align*}
and a coconjugation automorphism
\[
	\coc_\gamma \colon \uA[C_p/e] \to \uA[C_p/e],
\]
where we fix a generator $\gamma$ of $C_p$. The structure of the co-Tambara object $A[-]$ can be visualized with the following diagram, dual to a Lewis diagram:
\[
	\uA[-] \colon \quad
	\begin{tikzcd}[row sep = large]
		\uA[C_p/C_p]
		\ar[d, "\cotr", bend left=50]
		\ar[d, "\conm"', bend right=50]
		\\
		\uA[C_p/e]
		\ar[u, "\cores" description]
		\arrow[from=2-1, to=2-1, loop, in=240, out=300, distance=5mm, "\coc_\gamma"]
	\end{tikzcd}
\]

These co-Tambara operations corepresent (by the universal property of free Tambara functors) the operations of $C_p$-Tambara functors. Since these free Tambara functors are corepresentable up to group completion, i.e.
\[\uA[C_p/H] = \Poly(C_p/H, {-})^+,\]
we have by the Yoneda lemma that the co-Tambara operations are given by precomposition with distinguished morphisms in $\Poly$. For example, letting $q \colon C_p/e \to C_p/C_p$ be the unique map of $G$-sets,
\[\cores \colon \uA[C_p/e] \to \uA[C_p/C_p]\]
is induced by
\[(R_q)^* \colon \Poly(C_p/e, {-}) \to \Poly(C_p/C_p, {-}),\]
where we recall the distinguished bispans $T, N, R$ from \eqref{eq: distinguished bispans}.
Likewise, $\cotr$ is induced by $(T_q)^*$, and $\conm$ is induced by $(N_q)^*$. Finally, $\coc_\gamma$ is induced by $(c_\gamma)^*$, where $\conj_\gamma \colon C_p/e \to C_p/e$ is the morphism of $C_p$-sets sending the trivial coset $\{e\}$ to $\{\gamma\}$.

The Tambara affine line $\underline{\A}^{\!1}$ is the result of applying the Nakaoka spectrum functor to the co-Tambara object $\uA[-]$: it is the data of the Nakaoka spectra of $\uA[C_p/C_p]$ and $\uA[C_p/e]$, and the maps between these induced by cotransfer, conorm, corestriction, and coconjugation. These data may be organized into a Lewis diagram valued in topological spaces:
\[
	\underline{\A}^{\!1} \colon \quad
	\begin{tikzcd}[row sep = large]
		\Spec(\uA[C_p/C_p])
		\ar[d, "\cores^*" description]
		\\
		\Spec(\uA[C_p/e]).
		\ar[u, "\cotr^*", bend left=50]
		\ar[u, "\conm^*"', bend right=50]
		\arrow[from=2-1, to=2-1, loop, in=300, out=240, distance=5mm, "\coc_\gamma^*"']
	\end{tikzcd}
\]

In \cref{sec: main example}, we described the Nakaoka spectra of the above diagram. It remains to describe the structure morphisms. Our descriptions of $\Spec(\uA[C_p/e])$ and $\Spec(\uA[C_p/C_p])$ came through our theorems on the ghost construction, which reduced the problem to understanding the rings at each level of the free $C_p$-Tambara functors. We understood these rings by giving presentations for them in \Cref{Ex:AxGe,Ex:AxGG}. In order to now understand the action of $\cores$, $\conm$, $\cotr$, $\coc$ on prime ideals, we must translate between these presentations and the description of elements of a free Tambara functor as (formal differences) of bispans in $G\mhyphen\set$. We record in \Cref{table:dictionary} a dictionary between the generators of the presentations of the levels of the free $C_p$-Tambara functors and bispans in $C_p$-set.

\begin{table}[h!]
	\[\begin{array}{c|ccc}
			\text{Ring}           & \text{Generator} &                 & \text{Corresponding Bispan} \\ \hline \\
			\uA[C_p/C_p](C_p/C_p) & t                & \leftrightarrow &
			\begin{tikzcd}[ampersand replacement=\&]
				\& \varnothing \& {C_p/e} \\
				{C_p/C_p} \&\&\& {C_p/C_p}
				\arrow[from=1-2, to=1-3]
				\arrow[from=1-2, to=2-1]
				\arrow[from=1-3, to=2-4]
			\end{tikzcd}                                                  \\ \\
			\uA[C_p/C_p](C_p/C_p) & x                & \leftrightarrow &
			\begin{tikzcd}[ampersand replacement=\&]
				\& {C_p/C_p} \& {C_p/C_p} \\
				{C_p/C_p} \&\&\& {C_p/C_p}
				\arrow[from=1-2, to=1-3]
				\arrow[from=1-2, to=2-1]
				\arrow[from=1-3, to=2-4]
			\end{tikzcd}                                                  \\ \\
			\uA[C_p/C_p](C_p/C_p) & n                & \leftrightarrow &
			\begin{tikzcd}[ampersand replacement=\&]
				\& {C_p/e} \& {C_p/C_p} \\
				{C_p/C_p} \&\&\& {C_p/C_p}
				\arrow[from=1-2, to=1-3]
				\arrow[from=1-2, to=2-1]
				\arrow[from=1-3, to=2-4]
			\end{tikzcd}                                                  \\ \\
			\uA[C_p/C_p](C_p/e)   & x                & \leftrightarrow &
			\begin{tikzcd}[ampersand replacement=\&]
				\& {C_p/e} \& {C_p/e} \\
				{C_p/C_p} \&\&\& {C_p/e}
				\arrow[from=1-2, to=1-3]
				\arrow[from=1-2, to=2-1]
				\arrow[from=1-3, to=2-4]
			\end{tikzcd}                                                  \\ \\
			\uA[C_p/e](C_p/C_p)   & n                & \leftrightarrow &
			\begin{tikzcd}[ampersand replacement=\&]
				\& {C_p/e} \& {C_p/C_p} \\
				{C_p/e} \&\&\& {C_p/C_p}
				\arrow[from=1-2, to=1-3]
				\arrow[from=1-2, to=2-1]
				\arrow[from=1-3, to=2-4]
			\end{tikzcd}                                                  \\ \\
			\uA[C_p/e](C_p/C_p)   & t_{\vec v}       & \leftrightarrow &
			\begin{tikzcd}[ampersand replacement=\&]
				\& {\coprod_{i=0}^{p-1} (C_p/e)^{\amalg v_i}} \&\& {C_p/e} \\
				{C_p/e} \&\&\&\& {C_p/C_p}
				\arrow["{(s_{\vec v, 0}, \dots, s_{\vec v, p-1})}", from=1-2, to=1-4]
				\arrow[from=1-2, to=2-1]
				\arrow[from=1-4, to=2-5]
			\end{tikzcd}                     \\ \\
			\uA[C_p/e](C_p/e)     & x_i              & \leftrightarrow &
			\begin{tikzcd}[ampersand replacement=\&]
				\& {C_p/e} \& {C_p/e} \\
				{C_p/e} \&\&\& {C_p/e}
				\arrow["{c_{\gamma^i}}", from=1-2, to=1-3]
				\arrow[from=1-2, to=2-1]
				\arrow[from=1-3, to=2-4]
			\end{tikzcd}
		\end{array}\]

	\[\text{where} \quad s_{\vec v, i} = (\underbrace{c_{\gamma^i}, \dots, c_{\gamma^i}}_{v_i}).\]
	\caption{Dictionary between generators of our ring presentations for the free $C_p$-Tambara functors and bispans.}
	\label{table:dictionary}
\end{table}

From these descriptions and the fact that $\cores$, $\conm$, $\cotr$, and $\coc$ are morphisms of Tambara functors, we compute directly the behaviour of $\cores$, $\conm$, $\cotr$, and $\coc$ on the generators (\Cref{fig:co-ops on generators}).

\begin{figure}[h!]
	\[\begin{array}{ccc|ccc}
			\Poly(C_p/e,C_p/C_p)   & \xrightarrow{\cores} & \Poly(C_p/C_p,C_p/C_p)       & \Poly(C_p/C_p,C_p/C_p) & \xrightarrow{\conm}       & \Poly(C_p/e,C_p/C_p)     \\
			n                      & \mapsto              & n                            & t                      & \mapsto                   & \tr(1) = t_{(0,\dots,0)} \\
			t_{\vec v}             & \mapsto              & t x^{(\sum_{i=0}^{p-1} v_i)} & x                      & \mapsto                   & n                        \\
			                       &                      &                              & n                      & \mapsto                   & n^p                      \\
			\Poly(C_p/e,C_p/e)     & \xrightarrow{\cores} & \Poly(C_p/C_p,C_p/e)                                                                                         \\
			x_i                    & \mapsto              & x                            & \Poly(C_p/C_p,C_p/e)   & \xrightarrow{\conm}       & \Poly(C_p/e,C_p/e)       \\
			                       &                      &                              & x                      & \mapsto                   & \prod_{i=0}^{p-1} x_i    \\

			\Poly(C_p/C_p,C_p/C_p) & \xrightarrow{\cotr}  & \Poly(C_p/e,C_p/C_p)                                                                                         \\
			t                      & \mapsto              & \tr(1) = t_{(0,\dots,0)}     & \Poly(C_p/e,C_p/e)     & \xrightarrow{\coc_\gamma} & \Poly(C_p/e,C_p/e)       \\
			x                      & \mapsto              & t_{(1,0,\dots,0)} = \tr(x_0) & x_i                    & \mapsto                   & x_{(i+1\!\!\pmod{p})}    \\
			n                      & \mapsto              & \nm(\sum_{i=0}^{p-1} x_i)                                                                                    \\
			                       &                      &                              & \Poly(C_p/e,C_p/C_p)   & \xrightarrow{\coc_\gamma} & \Poly(C_p/e,C_p/C_p)     \\

			\Poly(C_p/C_p,C_p/e)   & \xrightarrow{\cotr}  & \Poly(C_p/e,C_p/e)           & n                      & \mapsto                   & n                        \\
			x                      & \mapsto              & \sum_i x_i                   & t_{\vec v}             & \mapsto                   & t_{\vec v}
		\end{array}\]
	\caption{The action of the co-Tambara structure morphisms on the levelwise generators of the free $C_p$-Tambara functors.}
	\label{fig:co-ops on generators}
\end{figure}

This allows us to compute the actions of $\ghost(\cores)$, $\ghost(\cotr)$, $\ghost(\conm)$, and $\ghost(\coc_\gamma)$ on primes. In the following descriptions, recall that $\A^1$ is a ring object in schemes and that $\A^p \cong (\A^1)^p$.  In particular, there are addition and multiplication maps $\A^p\to \A^1$, which we denote by $\Sigma$ and $\Pi$, respectively.

\begin{proposition}
	\label{prop: corestriction}
	The map $(\ghost(\cores))^* \colon \Spec \ghost (\uA[C_p/C_p]) \to \Spec \ghost(\uA[C_p/e])$ sends
	\begin{align*}
		(\mf{a}; \Z[x,n])         & \mapsto (\Delta \mf{a}; \Z[n])                           \\
		(\nm^{-1} \mf{b}; \mf{b}) & \mapsto (\nm^{-1} (\Z[n]\cap \mf{b}); \Z[n] \cap \mf{b})
	\end{align*}
	where $\Delta \colon \Spec(\Z[x]) \cong \A^1 \to \A^p \cong \Spec(\Z[x_0, \dots, x_{p-1}])$ is the diagonal embedding.
\end{proposition}

\begin{proposition}
	\label{prop: cotransfer}
	The map $(\ghost(\cotr))^* \colon \Spec \ghost(\uA[C_p/e]) \to \Spec \ghost(\uA[C_p/C_p])$ sends
	\begin{align*}
		(\mf{a}; \Z[n])           & \mapsto (\Sigma \mf{a}; \Z[x,n])                    \\
		(\nm^{-1} \mf{b}; \mf{b}) & \mapsto (\nm^{-1} \alpha_p \mf{b}; \alpha_p \mf{b})
	\end{align*}
	where $\alpha_p \colon \Spec(\Z[n]) \to \Spec(\Z[x,n])$ is induced by the map $\Z[x,n]\to \Z[n]$ which sends $x$ to $0$ and $n$ to $pn$.
\end{proposition}

\begin{proposition}
	\label{prop: conorm}
	The map  $(\ghost(\conm))^* \colon \Spec \ghost(\uA[C_p/e]) \to \Spec \ghost(\uA[C_p/C_p])$ sends
	\begin{align*}
		(\mf{a}; \Z[n])           & \mapsto (\Pi \mf{a}; \Z[x,n])                       \\
		(\nm^{-1} \mf{b}; \mf{b}) & \mapsto (\nm^{-1}\delta_p(\mf{b}); \delta_p \mf{b})
	\end{align*}
	where $\delta_p \colon \Spec(\Z[n]) \to \Spec(\Z[x,n])$ is induced by the map $\Z[x,n]\to \Z[n]$ which sends $x$ to $n$ and $n$ to $n^p$.
\end{proposition}

\begin{proposition}
	\label{prop: coconjugation}
	$(\ghost(\coc_\gamma))^* \colon \Spec \ghost(\uA[C_p/e]) \to \Spec \ghost(\uA[C_p/e])$ is the identity.
\end{proposition}

The above propositions effectively compute the effect of the co-Tambara operations on the spaces $\Spec(\uA[-])$.  For example, consider the operation $\cores^*\colon \Spec(\uA[C_p/C_p])\to \Spec(\uA[C_p/e])$.  There is a commutative square
\[
	\begin{tikzcd}
		\Spec(\uA[C_p/C_p])
		\ar[r, "\cores^*"]
		&
		\Spec(\uA[C_p/e])
		\\
		\Spec(\ghost(\uA[C_p/C_p]))
		\ar[r, "\ghost(\cores)^*"]
		\ar[u, "\ghostmap^*", two heads]
		&
		\Spec(\ghost(\uA[C_p/e]))
		\ar[u, "\ghostmap^*", two heads]
	\end{tikzcd}
\]
which, because the vertical maps are surjections, tells us that $\cores^*$ can be computed by first pulling back to the spectrum of the ghost, using \cref{prop: corestriction}, and the applying $\ghostmap^*$ again.  For instance, if $(\mf{a};\Z[x,n])\in \Spec(\ghost(\uA[C_p/C_p]))$ is a prime ideal then
\[
	\cores^*(\ghostmap^*(\mf{a};\Z[x,n])) = \ghostmap^*(\Delta\mf{a};\Z[n]).
\]

Completely analogous computations can be done for $\cotr^*$ and $\conm^*$.  The action of $\coc_{\gamma}^*$ is, of course, the identity. This completes our description of the Tambara affine line.

\bibliographystyle{alpha}
\bibliography{references.bib}

@article{6A19,
  author     = {Barthel, Tobias and Hausmann, Markus and Naumann, Niko and
                Nikolaus, Thomas and Noel, Justin and Stapleton, Nathaniel},
  title      = {The {B}almer spectrum of the equivariant homotopy category of
                a finite abelian group},
  journal    = {Invent. Math.},
  fjournal   = {Inventiones Mathematicae},
  volume     = {216},
  year       = {2019},
  number     = {1},
  pages      = {215--240},
  issn       = {0020-9910,1432-1297},
  doi        = {10.1007/s00222-018-0846-5},
  url        = {https://doi.org/10.1007/s00222-018-0846-5},
  note		= {\url{https://doi.org/10.1007/s00222-018-0846-5}}
}

@article{6A25,
 author = {Calle, Maxine E. and Chan, David and Mehrle, David and Quigley, J. D. and Spitz, Ben and Van Niel, Danika},
 title = {The spectrum of the {Burnside} {Tambara} functor},
journal = {Int. Math. Res. Notices (IMRN), to appear},
fjournal = {International Mathematics Research Notices (IMRN), to appear},
 year = {2025}
}

@article{Bal05,
  author     = {Balmer, Paul},
  title      = {The spectrum of prime ideals in tensor triangulated categories},
  journal    = {J. Reine Angew. Math.},
  fjournal   = {Journal f\"ur die Reine und Angewandte Mathematik. [Crelle's
                Journal]},
  volume     = {588},
  year       = {2005},
  pages      = {149--168},
  issn       = {0075-4102,1435-5345},
  doi        = {10.1515/crll.2005.2005.588.149},
  note        = {\url{https://doi.org/10.1515/crll.2005.2005.588.149}}
}

@inproceedings{Bal10,
  author     = {Balmer, Paul},
  title      = {Tensor triangular geometry},
  booktitle  = {Proceedings of the {I}nternational {C}ongress of {M}athematicians. {V}olume {II}},
  pages      = {85--112},
  publisher  = {Hindustan Book Agency, New Delhi},
  year       = {2010},
  note		= {\url{https://www.mathunion.org/fileadmin/ICM/Proceedings/ICM2010.2/ICM2010.2.pdf}},
  isbn       = {978-81-85931-08-3; 978-981-4324-32-8; 981-4324-32-9},
}

@article{BGH20,
  author     = {Barthel, Tobias and Greenlees, J. P. C. and Hausmann, Markus},
  title      = {On the {B}almer spectrum for compact {L}ie groups},
  journal    = {Compos. Math.},
  fjournal   = {Compositio Mathematica},
  volume     = {156},
  year       = {2020},
  number     = {1},
  pages      = {39--76},
  issn       = {0010-437X,1570-5846},
  mrclass    = {55P42 (55P91)},
  mrnumber   = {4036448},
  mrreviewer = {Samik\ Basu},
  doi        = {10.1112/s0010437x19007656},
  note        = {\url{https://doi.org/10.1112/s0010437x19007656}}
}

@article{BGHL2019,
  author     = {Blumberg, Andrew J. and Gerhardt, Teena and Hill, Michael A. and Lawson, Tyler},
  title      = {The {W}itt vectors for {G}reen functors},
  journal    = {J. Algebra},
  fjournal   = {Journal of Algebra},
  volume     = {537},
  year       = {2019},
  pages      = {197--244},
  issn       = {0021-8693,1090-266X},
  doi        = {10.1016/j.jalgebra.2019.07.014},
  note        = {\url{https://doi.org/10.1016/j.jalgebra.2019.07.014}}
}

@article{BH2015,
  author     = {Blumberg, Andrew J. and Hill, Michael A.},
  title      = {Operadic multiplications in equivariant spectra, norms, and transfers},
  journal    = {Adv. Math.},
  fjournal   = {Advances in Mathematics},
  volume     = {285},
  year       = {2015},
  pages      = {658--708},
  issn       = {0001-8708,1090-2082},
  doi        = {10.1016/j.aim.2015.07.013},
  note        = {\url{https://doi.org/10.1016/j.aim.2015.07.013}}
}

@incollection{BH2019,
  author     = {Blumberg, Andrew J. and Hill, Michael A.},
  booktitle  = {Homotopy theory: tools and applications},
  publisher  = {Amer. Math. Soc.},
  title      = {The right adjoint to the equivariant operadic forgetful functor on incomplete {T}ambara functors},
  year       = {2019},
  pages      = {75--92},
  series     = {Contemp. Math.},
  volume     = {729},
  doi        = {10.1090/conm/729/14691},
  note        = {\url{https://doi.org/10.1090/conm/729/14691}}
}

@article{Bru2005,
  author   = {Brun, Morten},
  title    = {Witt vectors and {T}ambara functors},
  journal  = {Adv. Math.},
  fjournal = {Advances in Mathematics},
  volume   = {193},
  year     = {2005},
  number   = {2},
  pages    = {233--256},
  issn     = {0001-8708,1090-2082},
  doi      = {10.1016/j.aim.2004.05.002},
  note      = {\url{https://doi.org/10.1016/j.aim.2004.05.002}}
}

@article{BS17,
  author     = {Balmer, Paul and Sanders, Beren},
  title      = {The spectrum of the equivariant stable homotopy category of a
                finite group},
  journal    = {Invent. Math.},
  fjournal   = {Inventiones Mathematicae},
  volume     = {208},
  year       = {2017},
  number     = {1},
  pages      = {283--326},
  issn       = {0020-9910,1432-1297},
  mrclass    = {18E30 (55P42 55U35)},
  mrnumber   = {3621837},
  mrreviewer = {Geoffrey\ M. L. Powell},
  doi        = {10.1007/s00222-016-0691-3},
  note        = {\url{https://doi.org/10.1007/s00222-016-0691-3}}
}

@article{burnside_table_of_marks,
  title   = {On the Representation of a Group of Finite Order as a Permutation Group, and on the Composition of Permutation Groups},
  volume  = {s1-34},
  rights  = {© 1901 London Mathematical Society},
  issn    = {1460-244X},
  note     = {\url{https://doi.org/10.1112/plms/s1-34.1.159}},
  doi     = {10.1112/plms/s1-34.1.159},
  pages   = {159--168},
  number  = {1},
  journal = {Proceedings of the London Mathematical Society},
  author  = {Burnside, W.},
  urldate = {2024-09-28},
  year    = {1901},
  langid  = {english}
}

@article{CG2023,
  title    = {The spectrum of the {B}urnside {T}ambara functor of a cyclic group},
  volume   = {227},
  issn     = {0022-4049},
  note      = {\url{https://doi.org/10.1016/j.jpaa.2023.107344}},
  doi      = {https://doi.org/10.1016/j.jpaa.2023.107344},
  pages    = {107344},
  number   = {8},
  journal  = {Journal of Pure and Applied Algebra},
  author   = {Calle, Maxine and Ginnett, Sam},
  year     = {2023}
}

@article{DHS,
  author     = {Devinatz, Ethan S. and Hopkins, Michael J. and Smith, Jeffrey H.},
  title      = {Nilpotence and stable homotopy theory. {I}},
  journal    = {Ann. of Math. (2)},
  fjournal   = {Annals of Mathematics. Second Series},
  volume     = {128},
  year       = {1988},
  number     = {2},
  pages      = {207--241},
  issn       = {0003-486X,1939-8980},
  doi        = {10.2307/1971440},
  note        = {\url{https://doi.org/10.2307/1971440}}
}

@book{Dress1971,
  author    = {Dress, Andreas W. M.},
  title     = {Notes on the theory of representations of finite groups. {P}art {I}: {T}he {B}urnside ring of a finite group and some {AGN}-applications},
  publisher = {Universit\"{a}t Bielefeld, Fakult\"{a}t f\"{u}r Mathematik, Bielefeld},
  year      = {1971},
  pages     = {iv+158+A28+B31 pp. (loose errata)}
}

@book{Eisenbud,
  author    = {Eisenbud, David},
  title     = {Commutative algebra with a view toward algebraic geometry},
  series    = {Graduate Texts in Mathematics},
  volume    = {150},
  publisher = {Springer-Verlag, New York},
  year      = {1995},
  pages     = {xvi+785},
  isbn      = {0-387-94268-8; 0-387-94269-6},
  doi       = {10.1007/978-1-4612-5350-1},
  note       = {\url{https://doi.org/10.1007/978-1-4612-5350-1}}
}

@article{ElmantoHaugseng,
  author     = {Elmanto, Elden and Haugseng, Rune},
  title      = {On distributivity in higher algebra {I}: the universal property of bispans},
  journal    = {Compos. Math.},
  fjournal   = {Compositio Mathematica},
  volume     = {159},
  year       = {2023},
  number     = {11},
  pages      = {2326--2415},
  issn       = {0010-437X,1570-5846},
  doi        = {10.1112/s0010437x23007388},
  note        = {\url{https://doi.org/10.1112/s0010437x23007388}}
}

@article{HH16,
  author       = {Hill, Michael A. and Hopkins, Michael J.},
  title        = {Equivariant symmetric monoidal structures},
  journal = {arXiv preprint, arXiv:1610.03114},
  year         = {2016},
  note          = {\url{https://arxiv.org/abs/1610.03114}},
  arxiv        = {arXiv:1610.03114}
}

@article{HM2019,
  title    = {An equivariant tensor product on {M}ackey functors},
  journal  = {Journal of Pure and Applied Algebra},
  volume   = {223},
  number   = {12},
  pages    = {5310-5345},
  year     = {2019},
  issn     = {0022-4049},
  doi      = {https://doi.org/10.1016/j.jpaa.2019.04.001},
  note 		= {\url{https://doi.org/10.1016/j.jpaa.2019.04.001}},
  author   = {Michael A. Hill and Kristen Mazur}
}

@article{HS98,
  author     = {Hopkins, Michael J. and Smith, Jeffrey H.},
  title      = {Nilpotence and stable homotopy theory. {II}},
  journal    = {Ann. of Math. (2)},
  fjournal   = {Annals of Mathematics. Second Series},
  volume     = {148},
  year       = {1998},
  number     = {1},
  pages      = {1--49},
  issn       = {0003-486X,1939-8980},
  mrclass    = {55P42 (55N20 55Q10)},
  mrnumber   = {1652975},
  mrreviewer = {David\ A.\ Blanc},
  doi        = {10.2307/120991},
  note        = {\url{https://doi.org/10.2307/120991}}
}

@book{Kaplansky,
  author    = {Kaplansky, Irving},
  title     = {Commutative rings},
  edition   = {Revised},
  publisher = {University of Chicago Press, Chicago, Ill.-London},
  year      = {1974},
  pages     = {ix+182},
}

@article{MandellMay,
  author     = {Mandell, M. A. and May, J. P.},
  title      = {Equivariant orthogonal spectra and {$S$}-modules},
  journal    = {Mem. Amer. Math. Soc.},
  fjournal   = {Memoirs of the American Mathematical Society},
  volume     = {159},
  year       = {2002},
  number     = {755},
  pages      = {x+108},
  issn       = {0065-9266,1947-6221},
  mrclass    = {55P91 (18E30 55P42 55P43 55P48)},
  mrnumber   = {1922205},
  mrreviewer = {J.\ P. C. Greenlees},
  doi        = {10.1090/memo/0755},
  url        = {https://doi.org/10.1090/memo/0755}
}

@phdthesis{Maz2013,
  title  = {On the {{Structure}} of {{Mackey Functors}} and {{Tambara Functors}}},
  author = {Mazur, Kristen},
  year   = {2013},
  month  = may,
  school = {University of Virginia}
}

@misc{milneANT,
  author = {Milne, James S.},
  title  = {Algebraic Number Theory (v3.08)},
  year   = {2020},
  note   = {\url{https://www.jmilne.org/math/}},
  pages  = {166}
}

@article{MQS2024,
  title         = {Koszul Resolutions over Free Incomplete {T}ambara Functors for Cyclic $p$-Groups},
  author        = {David Mehrle and J. D. Quigley and Michael Stahlhauer},
  year          = {2024},
  journal       = {arXiv preprint, arXiv:2407.18382},
  archiveprefix = {arXiv},
  primaryclass  = {math.AT},
  note          = {\url{https://arxiv.org/abs/2407.18382}}
}

@book{MumfordFogarty,
  author     = {Mumford, David and Fogarty, John},
  title      = {Geometric invariant theory},
  series     = {Ergebnisse der Mathematik und ihrer Grenzgebiete [Results in
                Mathematics and Related Areas]},
  volume     = {34},
  edition    = {Second},
  publisher  = {Springer-Verlag, Berlin},
  year       = {1982},
  pages      = {xii+220},
  isbn       = {3-540-11290-1},
}

@article{Nak2012,
  title    = {Ideals of {{Tambara}} Functors},
  author   = {Nakaoka, Hiroyuki},
  year     = {2012},
  month    = jul,
  volume   = {230},
  pages    = {2295--2331},
  issn     = {00018708},
  doi      = {10.1016/j.aim.2012.04.021},
  url      = {https://linkinghub.elsevier.com/retrieve/pii/S0001870812001612},
  journal  = {Advances in Mathematics},
  language = {en},
  number   = {4-6},
  note = {\url{https://doi.org/10.1016/j.aim.2012.04.021}}
}

@article {Nak2013,
    AUTHOR = {Nakaoka, Hiroyuki},
     TITLE = {A generalization of the {D}ress construction for a {T}ambara
              functor, and its relation to polynomial {T}ambara functors},
   JOURNAL = {Adv. Math.},
  FJOURNAL = {Advances in Mathematics},
    VOLUME = {235},
      YEAR = {2013},
     PAGES = {237--260},
      ISSN = {0001-8708,1090-2082},
   MRCLASS = {19A22},
  MRNUMBER = {3010058},
MRREVIEWER = {Charles\ Weibel},
       DOI = {10.1016/j.aim.2012.11.013},
       URL = {https://doi.org/10.1016/j.aim.2012.11.013},
       note = {\url{https://doi.org/10.1016/j.aim.2012.11.013}}
}

@article{Nak2014a,
  title    = {The Spectrum of the {{Burnside Tambara}} Functor on a Finite Cyclic p -Group},
  author   = {Nakaoka, Hiroyuki},
  year     = {2014},
  month    = jan,
  volume   = {398},
  pages    = {21--54},
  issn     = {00218693},
  doi      = {10.1016/j.jalgebra.2013.09.010},
  url      = {https://linkinghub.elsevier.com/retrieve/pii/S0021869313004997},
  urldate  = {2020-06-17},
  journal  = {Journal of Algebra},
  language = {en},
  note = {\url{10.1016/j.jalgebra.2013.09.010}}
}

@article{Str2012,
  title         = {Tambara Functors},
  author        = {Strickland, Neil},
  year          = {2012},
  note          = {\url{https://arxiv.org/abs/1205.2516}},
  urldate       = {2020-02-04},
  archiveprefix = {arXiv},
  eprint        = {1205.2516},
  eprinttype    = {arxiv},
  journal       = {arXiv preprint, arXiv:1205.2516},
  primaryclass  = {math}
}

@article{Sulyma2023,
  title         = {Prisms and {T}ambara functors {I}: Twisted powers, transversality, and the perfect sandwich},
  author        = {Yuri J. F. Sulyma},
  year          = {2023},
  journal       = {arXiv preprint, arXiv:2309.03181},
  eprint        = {2309.03181},
  archiveprefix = {arXiv},
  primaryclass  = {math.NT},
  url           = {https://arxiv.org/abs/2309.03181},
  note          = {\url{https://arxiv.org/abs/2309.03181}},
}

@article{Tam1993,
  title    = {On Multiplicative Transfer},
  author   = {Tambara, D.},
  year     = {1993},
  month    = jan,
  volume   = {21},
  pages    = {1393--1420},
  issn     = {0092-7872, 1532-4125},
  doi      = {10.1080/00927879308824627},
  url      = {http://www.tandfonline.com/doi/abs/10.1080/00927879308824627},
  note		= {\url{https://doi.org/10.1080/00927879308824627}},
  urldate  = {2020-09-28},
  journal  = {Communications in Algebra},
  language = {en},
  number   = {4}
}

@article{Thevenaz:SomeRemarks,
  author     = {Th{\'{e}}venaz, Jacques},
  title      = {Some remarks on {$G$}-functors and the {B}rauer morphism},
  journal    = {J. Reine Angew. Math.},
  fjournal   = {Journal f{\"{u}}r die Reine und Angewandte Mathematik. [Crelle's
                Journal]},
  volume     = {384},
  year       = {1988},
  pages      = {24--56},
  issn       = {0075-4102,1435-5345},
  mrclass    = {20C20},
  mrnumber   = {929977},
  mrreviewer = {David\ Benson},
  doi        = {10.1515/crll.1988.384.24},
  url        = {https://doi.org/10.1515/crll.1988.384.24},
  note        = {\url{https://doi.org/10.1515/crll.1988.384.24}}
}

\end{document}